# I-LAMM FOR SPARSE LEARNING: SIMULTANEOUS CONTROL OF ALGORITHMIC COMPLEXITY AND STATISTICAL ERROR*

By Jianqing Fan[†,‡], Han Liu[†], Qiang Sun[†], and Tong Zhang[§,‡]

*Princeton University[†], Fudan University[‡] and Tecent AI Lab[§]*

We propose a computational framework named iterative local adaptive majorize-minimization (I-LAMM) to simultaneously control algorithmic complexity and statistical error when fitting high dimensional models. I-LAMM is a two-stage algorithmic implementation of the local linear approximation to a family of folded concave penalized quasi-likelihood. The first stage solves a convex program with a crude precision tolerance to obtain a coarse initial estimator, which is further refined in the second stage by iteratively solving a sequence of convex programs with smaller precision tolerances. Theoretically, we establish a phase transition: the first stage has a sublinear iteration complexity, while the second stage achieves an improved linear rate of convergence. Though this framework is completely algorithmic, it provides solutions with optimal statistical performances and controlled algorithmic complexity for a large family of nonconvex optimization problems. The iteration effects on statistical errors are clearly demonstrated via a contraction property. Our theory relies on a localized version of the sparse/restricted eigenvalue condition, which allows us to analyze a large family of loss and penalty functions and provide optimality guarantees under very weak assumptions (For example, I-LAMM requires much weaker minimal signal strength than other procedures). Thorough numerical results are provided to support the obtained theory.

**1. Introduction.** Modern data acquisitions routinely measure massive amounts of variables, which can be much larger than the sample size, making statistical inference an ill-posed problem. For inferential tractability and interpretability, one common approach is to exploit the penalized M-estimator

$$(1.1) \qquad \widehat{\boldsymbol{\beta}} = \operatorname*{argmin}_{\boldsymbol{\beta}\in\mathbb{R}^d}\Big\{\mathcal{L}(\boldsymbol{\beta}) + \mathcal{R}_\lambda(\boldsymbol{\beta})\Big\},$$

*This research is supported in part by NIH Grants 5R01-GM072611-11, R01-GM100474-04, R01-MH102339, R01-GM083084 and R01-HG06841, NSF grants DMS-1206464-04, DMS-1308566, DMS-1454377, DMS-1206464, IIS-1408910, IIS-1332109, and Science and Technology Commission of Shanghai Municipality 16JC1402600.

*MSC 2010 subject classifications:* Primary 62J07; secondary 62C20, 62H35

where $\mathcal{L}(\cdot)$ is a smooth loss function, $\mathcal{R}_\lambda(\cdot)$ is a sparsity-inducing penalty with a regularization parameter $\lambda$. Our framework encompasses the square loss, logistic loss, Gaussian graphical model negative log-likelihood loss, Huber loss, and the family of folded concave penalties [9]. Finding optimal statistical procedures with controlled computational complexity characterizes the efforts of high-dimensional statistical learning in the last two decades. This paper makes an important leap toward this grand challenge by proposing a general algorithmic strategy for solving (1.1) even when $\mathcal{R}_\lambda(\boldsymbol{\beta})$ is nonconvex.

A popular choice of $\mathcal{R}_\lambda(\boldsymbol{\beta})$ is the Lasso penalty [25], a convex penalty. Though a large literature exists on understanding the theory of penalized M-estimators with convex penalties [8, 4, 26, 22], it has been well known [9, 33] that the convex penalties introduce non-negligible estimation biases. In addition, the algorithmic issues for finding global minimizer are rarely addressed. To eliminate the estimation bias, a family of folded-concave penalties was introduced by [9], which includes the smooth clipped absolute deviation (SCAD) [9], minimax concave penalty (MCP) [29], and capped $\ell_1$-penalty [32]. Compared to their convex counterparts, these nonconvex penalties eliminate the estimation bias and attain more refined statistical rates of convergence. However, it is more challenging to analyze the theoretical properties of the resulting estimators due to nonconvexity of the penalty functions. Existing work on nonconvex penalized M-estimators treats the statistical properties and practical algorithms separately. On one hand, statistical properties are established for the hypothetical global optimum (or some local minimum), which is usually unobtainable by any practical algorithm in polynomial time. For example, [9] showed that there exists a local solution that possesses an oracle property; [15] and [11] showed that the oracle estimator is a local minimizer with high probability. Later on, [16] and [30] proved that the global optimum achieves the oracle property under certain conditions. Nevertheless, none of these paper specifies an algorithm to find the desired solution. More recently, [20, 22, 1] develop a projected gradient algorithm with desired statistical guarantees. However, they need to modify the estimating procedures to include an additional $\ell_1$-ball constraint, $\|\boldsymbol{\beta}\|_1 \leq R$, which depends on the unknown true parameter. On the other hand, practitioners have developed numerous heuristic algorithms for nonconvex optimization problems, but without theoretical guarantees. One such example is the coordinate optimization strategy studied in [6] and [13].

So there is a gap between theory and practice: What is actually computed is not the same as what has been proved. To bridge this gap, we propose an iterative local adaptive majorize-minimization (I-LAMM) algorithm for fitting high dimensional statistical models. Unlike most existing methods, which are mainly motivated from a statistical perspective and ignore the computational consideration, I-LAMM is both algorithmic and statistical: it computes an estimator within polynomial time and achieves optimal statistical accuracy for this estimator. In particular, I-LAMM obtains estimators with the strongest statistical guarantees for a wide family of loss functions under the weakest possible assumptions. Moreover, the statistical properties are established for the estimators computed exactly by our algorithm, which is designed to control the cost of computing resources. Compared to existing works [20, 22, 1], our method does not impose any constraint that depends on the unknown true parameter.

Inspired by the local linear approximation to the folded concave penalty [34], we use I-LAMM to solve a sequence of convex programs up to a prefixed optimization precision

$$(1.2) \quad \min_{\boldsymbol{\beta} \in \mathbb{R}^d} \left\{ \mathcal{L}(\boldsymbol{\beta}) + \mathcal{R}(\boldsymbol{\lambda}^{(\ell-1)} \odot \boldsymbol{\beta}) \right\}, \text{ for } \ell = 1, \ldots, T,$$

where $\boldsymbol{\lambda}^{(\ell-1)} = \left(\lambda \mathrm{w}(|\widetilde{\beta}_1^{(\ell-1)}|), \ldots, \lambda \mathrm{w}(|\widetilde{\beta}_d^{(\ell-1)}|)\right)^{\mathrm{T}}$, $\widetilde{\boldsymbol{\beta}}^{(\ell)}$ is an approximate solution to the $\ell$th optimization problem in (1.2), $\mathrm{w}(\cdot)$ is a weighting function, $\mathcal{R}(\cdot)$ is a decomposable convex penalty function, and '$\odot$' denotes the Hadamard product. In this paper, we mainly consider $\mathcal{R}(\boldsymbol{\beta}) = \|\boldsymbol{\beta}\|_1$, though our theory is general. The weighting function corresponds to the derivative of the folded concave penalty in [9], [34] and [11].

In particular, the I-LAMM algorithm obtains a crude initial estimator $\widetilde{\boldsymbol{\beta}}^{(1)}$ and further solves the optimization problem (1.2) for $\ell \geq 2$ with established algorithmic and statistical properties. This provides theoretical insights on how fast the algorithm converges and how much computation is needed, as well as the desired statistical properties of the obtained estimator. The whole procedure consists of $T$ convex programs, each only needs to be solved approximately to control the computational cost. Under mild conditions, we show that only $\log(\lambda\sqrt{n})$ steps are needed to obtain the optimal statistical rate of convergence. Even though I-LAMM solves *approximately* a sequence of convex programs, the solution enjoys the same optimal statistical property of the unobtainable global optimum for the folded-concave penalized regression. The adaptive stopping rule for solving each convex program in (1.2) allows us to control both computational costs and statistical errors. Figure 1 provides a geometric illustration of the I-LAMM procedure. It contains a contraction stage and a tightening stage as described below.

∗ Contraction Stage: In this stage ($\ell = 1$), we approximately solve a convex optimization problem (1.2), starting from *any* initial value $\widetilde{\boldsymbol{\beta}}^{(0)}$, and

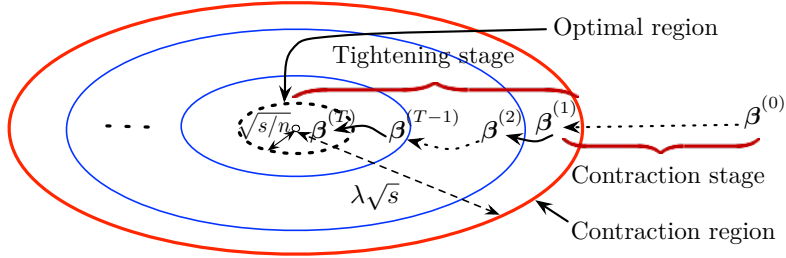

FIG 1. *Geometric illustration of the contraction property. The contraction stage produces an initial estimator, starting from* any *initial value $\widetilde{\boldsymbol{\beta}}^{(0)}$ that falls in the contraction region which secures the tightening stage to enjoy optimal statistical and computational rates of convergence. The tightening stage adaptively refines the contraction estimator till it enters the optimal region, which is stated in* (1.3). *Here $\lambda$ is a regularization parameter, $s$ the number of nonzero coefficients in $\boldsymbol{\beta}^*$ and $n$ the sample size.*

terminate the algorithm as long as the approximate solution enters a desired contraction region which will be characterized in Section 2.3. The obtained estimator is called the contraction estimator, which is very crude and only serves as initialization.

∗ Tightening Stage: This stage involves multiple tightening steps ($\ell \geq 2$). Specifically, we iteratively tighten the contraction estimator by solving a sequence of convex programs. Each step contracts its initial estimator towards the true parameter until it reaches the optimal region of convergence. At that region, further iteration does not improve statistical performance. See Figure 1. More precisely, we will show the following contraction property

$$(1.3) \qquad \|\widetilde{\boldsymbol{\beta}}^{(\ell)} - \boldsymbol{\beta}^*\|_2 \lesssim \sqrt{\frac{s}{n}} + \delta \cdot \|\widetilde{\boldsymbol{\beta}}^{(\ell-1)} - \boldsymbol{\beta}^*\|_2 \quad \text{for } \ell \geq 2,$$

where $\boldsymbol{\beta}^*$ is the true regression coefficient, $\delta \in (0,1)$ a prefixed contraction parameter and $\sqrt{s/n}$ the order of statistical error. Tightening helps improve the accuracy only when $\|\widetilde{\boldsymbol{\beta}}^{(\ell-1)} - \boldsymbol{\beta}^*\|_2$ dominates the statistical error. The iteration effect is clearly demonstrated. Since $\widetilde{\boldsymbol{\beta}}^{(\ell)}$ is only used to create an adaptive weight for $\widetilde{\boldsymbol{\beta}}^{(\ell+1)}$, we can control the iteration complexity by solving each subproblem in (1.2) approximately. What differs from the contraction stage is that the initial estimators in the tightening stage are already in the contraction region, making the optimization algorithm enjoy geometric rate of convergence. This allows us to rapidly solve (1.2) with small optimization error.

∗ (Phase Transition in Algorithmic Convergence) In the contraction stage ($\ell = 1$), the optimization problem is not strongly convex and therefore our

algorithm has only a sublinear convergence rate. Once the solution enters the contraction region, we will show that the feasible solutions are sparse and the objective function is essentially 'low'-dimensional and becomes (restricted) strongly convex and smooth in that region. Therefore, our algorithm has a linear convergence rate for $\ell > 1$. Indeed, this holds even for $\ell = 1$, which admits a sublinear rate until it enters into the contraction region and enjoys a linear rate of convergence after that. See Figure 2. But this estimator (for $\ell = 1$) is the estimator that corresponds to LASSO penalty, not folded concave penalty that we are looking for.

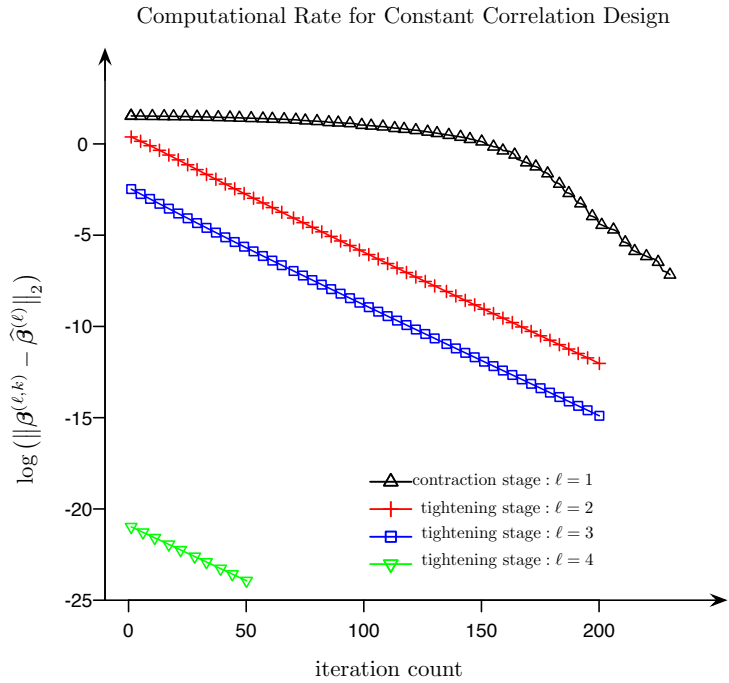

FIG 2. *Computational rate of convergence in each stage for the simulation experiment specified in case 2 in Example 6.1. The x-axis is the iteration count $k$ within the $\ell$th subproblem. The phase transition from sublinear rate to liner rate of algorithmic convergence is clearly seen once the iterations enter the contraction region. Here $\widehat{\boldsymbol{\beta}}^{(\ell)}$ is the global minimizer of the $\ell$th optimization problem in (1.2) and $\boldsymbol{\beta}^{(\ell,k)}$ is its kth iteration (see Figure 3). For $\ell = 1$, the initial estimation sequence has sublinear rate and once the solution sequence enters the contraction region, it becomes linear convergent. For $\ell \geq 2$, the algorithm achieves linear rate, since all estimators $\boldsymbol{\beta}^{(\ell,k-1)}$ are in the contraction region.*

This paper makes four major contributions. First, I-LAMM offers an algorithmic approach to obtain the optimal estimator with controlled computing resources. Second, compared to the existing literature, our method requires weaker conditions due to a novel localized analysis of sparse learning problems. Specifically, our method does not need the extra ball constraint as in [20] and [28], which is an artifact of their proofs. Third, our computational framework takes the approximate optimization error into analysis and provides theoretical guarantees for the estimator that is computed by the algorithm. Fourth, our method provides new theoretical insights about the adaptive lasso and folded-concave penalized regression. In particular, we bridge these two methodologies together using a unified framework. See Section 3.2 for more details.

The rest of this paper proceeds as follows. In Section 2, we introduce I-LAMM and its implementation. Section 3 is contributed to new insights into existing methods for high dimensional regression. In Section 4, we introduce both the localized sparse eigenvalue and localized restricted eigenvalue conditions. Statistical property and computational complexity are then presented. In Section 5, we outline the key proof strategies. Numerical simulations are provided to evaluate the proposed method in Section 6. We conclude by discussions in Section 7. All the proofs are postponed to the online supplement.

**Notation:** For $\mathbf{u} = (u_1, u_2, \ldots, u_d)^{\mathrm{T}} \in \mathbb{R}^d$, we define the $\ell_q$-norm of $\mathbf{u}$ by $\|\mathbf{u}\|_q = (\sum_{j=1}^d |u_j|^q)^{1/q}$, where $q \in [1, \infty)$. Let $\|\mathbf{u}\|_{\min} = \min\{u_j : 1 \leq j \leq d\}$. For a set $\mathcal{S}$, let $|\mathcal{S}|$ denote its cardinality. We define the $\ell_0$-pseudo norm of $\mathbf{u}$ as $\|\mathbf{u}\|_0 = |\mathrm{supp}(\mathbf{u})|$, where $\mathrm{supp}(\mathbf{u}) = \{j : u_j \neq 0\}$. For an index set $\mathcal{I} \subseteq \{1, \ldots, d\}$, $\mathbf{u}_\mathcal{I} \in \mathbb{R}^d$ is defined to be the vector whose $i$-th entry is equal to $u_i$ if $i \in \mathcal{I}$ and zero otherwise. Let $\mathbf{A} = [a_{i,j}] \in \mathbb{R}^{d \times d}$. For $q \geq 1$, we define $\|\mathbf{A}\|_q$ as the matrix operator $q$-norm of $\mathbf{A}$. For index sets $\mathcal{I}, \mathcal{J} \subseteq \{1, \ldots, d\}$, we define $\mathbf{A}_{\mathcal{I},\mathcal{J}} \in \mathbb{R}^{d \times d}$ to be the matrix whose $(i,j)$-th entry is equal to $a_{i,j}$ if $i \in \mathcal{I}$ and $j \in \mathcal{J}$, and zero otherwise. We use $\mathrm{sign}(x)$ to denote the sign of $x$: $\mathrm{sign}(x) = x/|x|$ if $x \neq 0$ and $\mathrm{sign}(x) = 0$ otherwise. For two functionals $f(n, d, s)$ and $g(n, d, s)$, we denote $f(n, d, s) \gtrsim g(n, d, s)$ if $f(n, d, s) \geq C g(n, d, s)$ for a constant $C$; $f(n, d, s) \lesssim g(n, d, s)$ otherwise.

**2. Methodology.** In this paper, we assume that the loss function $\mathcal{L}(\cdot) \in \mathcal{F}_\mathcal{L}$, a family of general convex loss functions specified in Appendix A.

2.1. *Local Adaptive Majorize-Minimization.* Recall that the estimators are obtained by solving a sequence of convex programs in (1.2). We require the function $\mathrm{w}(\cdot)$ used therein to be taken from the tightening function class

$\mathcal{T}$, defined as

$$(2.1) \quad \mathcal{T} = \Big\{ \mathrm{w}(\cdot) \in \mathcal{M} : \mathrm{w}(t_1) \leq \mathrm{w}(t_2) \text{ for all } t_1 \geq t_2 \geq 0,$$
$$0 \leq \mathrm{w}(t) \leq 1 \text{ if } t \geq 0, \mathrm{w}(t) = 0 \text{ if } t \leq 0 \Big\}.$$

To fix ideas, we take $\mathcal{R}_\lambda(\boldsymbol{\beta})$ in (1.1) to be $\sum_{j=1}^{d} p_\lambda(|\beta_j|)$, where $p_\lambda(\cdot)$ is a folded concave penalty [9] such as the SCAD or MCP. As discussed in [9], the penalized likelihood function in (1.1) is folded concave with respect to $\boldsymbol{\beta}$, making it difficult to be maximized. We propose to use the adaptive local linear approximation (adaptive LLA) to the penalty function [34, 12] and approximately solve

$$(2.2) \quad \operatorname*{argmin}_{\boldsymbol{\beta}} \Big\{ \mathcal{L}(\boldsymbol{\beta}) + \sum_{j=1}^{d} p'_\lambda(|\widetilde{\beta}_j^{(\ell-1)}|)|\beta_j| \Big\}, \text{ for } 1 \leq \ell \leq T,$$

where $\widetilde{\beta}_j^{(\ell-1)}$ is the $j$th component of $\widetilde{\boldsymbol{\beta}}^{(\ell-1)}$ and $\widetilde{\boldsymbol{\beta}}^{(0)}$ can be an arbitrary bad initial value: $\widetilde{\boldsymbol{\beta}}^{(0)} = \mathbf{0}$, for example. If we assume that $\mathrm{w}(\cdot) \equiv \lambda^{-1} p'_\lambda(\cdot) \in \mathcal{T}$, such as the SCAD or MCP, then the adaptive LLA algorithm can be regarded as a special case of our general formulation (1.2). Note that the LLA algorithm with $\ell_q$-penalty ($q < 1$) is not covered by our algorithm since its derivative is unbounded at the origin and thus $\lambda^{-1} p'_\lambda(\cdot) \notin \mathcal{T}$. The latter creates a zero-absorbing state: once a component is shrunk to zero, it will remain zero throughout the remaining iterations, as noted in [10]. Of course, we can truncate the loss derivative of the loss function to resolve this issue.

We now propose a local adaptive majorize-minimization (LAMM) principal, which will be repeatedly called to practically solve the optimization problem (2.2). We first review the majorize-minimization (MM) algorithm. To minimize a general function $f(\boldsymbol{\beta})$, at a given point $\boldsymbol{\beta}^{(k)}$, MM majorizes it by $g(\boldsymbol{\beta}|\boldsymbol{\beta}^{(k)})$, which satisfies

$$g(\boldsymbol{\beta}|\boldsymbol{\beta}^{(k)}) \geq f(\boldsymbol{\beta}) \quad \text{and} \quad g(\boldsymbol{\beta}^{(k)}|\boldsymbol{\beta}^{(k)}) = f(\boldsymbol{\beta}^{(k)})$$

and then compute $\boldsymbol{\beta}^{(k+1)} = \operatorname{argmin}_{\boldsymbol{\beta}} \{ g(\boldsymbol{\beta}|\boldsymbol{\beta}^{(k)}) \}$ [17, 14]. The objective value of such an algorithm is non-increasing in each step, since

$$(2.3) \quad f(\boldsymbol{\beta}^{(k+1)}) \overset{\text{major.}}{\leq} g(\boldsymbol{\beta}^{(k+1)} \,|\, \boldsymbol{\beta}^{(k)}) \overset{\text{min.}}{\leq} g(\boldsymbol{\beta}^{(k)} \,|\, \boldsymbol{\beta}^{(k)}) \overset{\text{init.}}{=} f(\boldsymbol{\beta}^{(k)}).$$

An inspection of the above arguments shows that the majorization requirement is not necessary. It requires only the local property:

$$(2.4) \quad f(\boldsymbol{\beta}^{(k+1)}) \leq g(\boldsymbol{\beta}^{(k+1)}|\boldsymbol{\beta}^{(k)}) \text{ and } g(\boldsymbol{\beta}^{(k)}|\boldsymbol{\beta}^{(k)}) = f(\boldsymbol{\beta}^{(k)})$$

for the inequalities in (2.3) to hold.

Inspired by the above observation, we locally majorize (2.2) at the $\ell$th step. It is similar to the iteration steps used in the (proximal) gradient method [5, 23]. Instead of computing and storing a large Hessian matrix as in [34], we majorize $\mathcal{L}(\boldsymbol{\beta})$ in (2.2) at $\widetilde{\boldsymbol{\beta}}^{(\ell-1)}$ by an isotropic quadratic function

$$\mathcal{L}(\widetilde{\boldsymbol{\beta}}^{(\ell-1)}) + \langle \nabla \mathcal{L}(\widetilde{\boldsymbol{\beta}}^{(\ell-1)}), \boldsymbol{\beta} - \widetilde{\boldsymbol{\beta}}^{(\ell-1)} \rangle + \frac{\phi}{2} \|\boldsymbol{\beta} - \widetilde{\boldsymbol{\beta}}^{(\ell-1)}\|_2^2,$$

where $\nabla$ is used to denote derivative. By Taylor's expansion, it suffices to take $\phi$ that is no smaller than the largest eigenvalue of $\nabla^2 \mathcal{L}(\widetilde{\boldsymbol{\beta}}^{(\ell-1)})$. More importantly, the isotropic form also allows a simple analytic solution to the subsequent majorized optimization problem:

$$(2.5) \quad \operatorname*{argmin}_{\boldsymbol{\beta} \in \mathbb{R}^d} \Big\{ \mathcal{L}(\widetilde{\boldsymbol{\beta}}^{(\ell-1)}) + \langle \nabla \mathcal{L}(\widetilde{\boldsymbol{\beta}}^{(\ell-1)}), \boldsymbol{\beta} - \widetilde{\boldsymbol{\beta}}^{(\ell-1)} \rangle$$
$$+ \frac{\phi}{2} \|\boldsymbol{\beta} - \widetilde{\boldsymbol{\beta}}^{(\ell-1)}\|_2^2 + \sum_{j=1}^d p'_\lambda(|\widetilde{\beta}_j^{(\ell-1)}|) |\beta_j| \Big\}.$$

With $\boldsymbol{\lambda}^{(\ell-1)} = \big(p'_\lambda(|\widetilde{\beta}_1^{(\ell-1)}|), \ldots, p'_\lambda(|\widetilde{\beta}_d^{(\ell-1)}|)\big)^{\mathrm{T}}$, it is easy to show that (2.5) is minimized at

$$\boldsymbol{\beta}^{(\ell,1)} = T_{\boldsymbol{\lambda}^{(\ell-1)},\phi}(\widetilde{\boldsymbol{\beta}}^{(\ell-1)}) \equiv S\big(\widetilde{\boldsymbol{\beta}}^{(\ell-1)} - \phi^{-1} \nabla \mathcal{L}(\widetilde{\boldsymbol{\beta}}^{(\ell-1)}), \phi^{-1} \boldsymbol{\lambda}^{(\ell-1)}\big),$$

where $S(\mathbf{x}, \boldsymbol{\lambda})$ is the soft-thresholding operator, defined by $S(\mathbf{x}, \boldsymbol{\lambda}) \equiv \big(\operatorname{sign}(x_j) \cdot \max\{|x_j| - \lambda_j, 0\}\big)$. The simplicity of this updating rule is due to the fact that (2.5) is an unconstrained optimization problem. This is not the case in [20] and [28].

However, finding the value of $\phi \geq \|\nabla^2 \mathcal{L}(\widetilde{\boldsymbol{\beta}}^{(\ell-1)})\|_2$ is not an easy task in computation. To avoid storing and computing the largest eigenvalue of a big matrix, we now state the LAMM algorithm, thanks to the local requirement (2.4). The basic idea of LAMM is to start from a very small isotropic parameter $\phi_0$ and then successively inflate $\phi$ by a factor $\gamma_u > 1$ (say, 2). If the solution satisfies (2.4), we stop this part of the algorithm, which will make the target value non-increasing. Since after the $k$th iteration, $\phi = \gamma_u^{k-1} \phi_0$, there always exists a $k$ such that it is no larger than $\|\nabla^2 \mathcal{L}(\widetilde{\boldsymbol{\beta}}^{(\ell-1)})\|_2$. In this manner, the LAMM algorithm will find a smallest iteration to make (2.4) hold.

Specifically, our proposed LAMM algorithm to solve (2.5) at $\widetilde{\boldsymbol{\beta}}^{(\ell-1)}$ begins with $\phi = \phi_0$, say $10^{-6}$, iteratively increases $\phi$ by a factor of $\gamma_u > 1$ inside the

**Algorithm 1** The LAMM algorithm in the $k$th iteration of the $\ell$th tightening subproblem.

1: **Algorithm**: $\{\boldsymbol{\beta}^{(\ell,k)}, \phi^{(\ell,k)}\} \leftarrow \text{LAMM}(\boldsymbol{\lambda}^{(\ell-1)}, \boldsymbol{\beta}^{(\ell,k-1)}, \phi_0, \phi^{(\ell,k-1)})$
2: **Input**: $\boldsymbol{\lambda}^{(\ell-1)}, \boldsymbol{\beta}^{(\ell,k-1)}, \phi_0, \phi^{(\ell,k-1)}$
3: **Initialize**: $\phi^{(\ell,k)} \leftarrow \max\{\phi_0, \gamma_u^{-1}\phi^{(\ell,k-1)}\}$
4: **Repeat**
5: $\quad \boldsymbol{\beta}^{(\ell,k)} \leftarrow T_{\boldsymbol{\lambda}^{(\ell-1)}, \phi^{(\ell,k)}}(\boldsymbol{\beta}^{(\ell,k-1)})$
6: $\quad$ **If** $F(\boldsymbol{\beta}^{(\ell,k)}, \boldsymbol{\lambda}^{(\ell-1)}) > \Psi_{\boldsymbol{\lambda}^{(\ell-1)}, \phi^{(\ell,k)}}(\boldsymbol{\beta}^{(\ell,k)}; \boldsymbol{\beta}^{(\ell,k-1)})$ **then** $\phi^{(\ell,k)} \leftarrow \gamma_u \phi^{(\ell,k)}$
7: **Until** $F(\boldsymbol{\beta}^{(\ell,k)}, \boldsymbol{\lambda}^{(\ell-1)}) \leq \Psi_{\boldsymbol{\lambda}^{(\ell-1)}, \phi^{(\ell,k)}}(\boldsymbol{\beta}^{(\ell,k)}; \boldsymbol{\beta}^{(\ell,k-1)})$
8: **Return** $\{\boldsymbol{\beta}^{(\ell,k)}, \phi^{(\ell,k)}\}$

$\ell$th step of optimization, and computes

$$\boldsymbol{\beta}^{(\ell,1)} = T_{\boldsymbol{\lambda}^{(\ell-1)}, \phi^{(\ell,k)}}(\boldsymbol{\beta}^{(\ell,0)}), \qquad \text{with } \phi^{(\ell,k)} = \gamma_u^{k-1}\phi_0, \quad \boldsymbol{\beta}^{(\ell,0)} = \widetilde{\boldsymbol{\beta}}^{(\ell-1)},$$

until the local property (2.4) holds. In our context, LAMM stops when

$$\Psi_{\boldsymbol{\lambda}^{(\ell-1)}, \phi^{(\ell,k)}}(\boldsymbol{\beta}^{(\ell,1)}, \boldsymbol{\beta}^{(\ell,0)}) \geq F(\boldsymbol{\beta}^{(\ell,1)}, \boldsymbol{\lambda}^{(\ell-1)}),$$

where $F(\boldsymbol{\beta}, \boldsymbol{\lambda}^{(\ell-1)}) \equiv \mathcal{L}(\boldsymbol{\beta}) + \sum_{j=1}^d \lambda_j^{(\ell-1)}|\beta_j|$ and

$$\Psi_{\boldsymbol{\lambda}^{(\ell-1)}, \phi^{(\ell,k)}}(\boldsymbol{\beta}, \boldsymbol{\beta}^{(\ell,0)}) \equiv \mathcal{L}(\boldsymbol{\beta}^{(\ell,0)}) + \langle \nabla \mathcal{L}(\boldsymbol{\beta}^{(\ell,0)}), \boldsymbol{\beta} - \boldsymbol{\beta}^{(\ell,0)} \rangle \\ + \frac{\phi^{(\ell,k)}}{2}\|\boldsymbol{\beta} - \boldsymbol{\beta}^{(\ell,0)}\|_2^2 + \sum_{j=1}^d \lambda_j^{(\ell-1)}|\beta_j|.$$

Inspired by [23], to accelerate LAMM within the next majorizing step, we keep track of the sequence $\{\phi^{(\ell,k)}\}_{\ell,k}$ and set $\phi^{(\ell,k)} = \max\{\phi_0, \gamma_u^{-1}\phi^{(\ell,k-1)}\}$, with the convention that $\phi_{\ell,0} = \widetilde{\phi}_{\ell-1}$ and $\widetilde{\phi}_0 = \phi_0$, in which $\widetilde{\phi}_{\ell-1}$ is the isotropic parameter corresponding to the solution $\widetilde{\boldsymbol{\beta}}^{(\ell-1)}$. This is summarized in Algorithm 1 with a generic initial value.

The LAMM algorithm solves only one local majorization step. It corresponds to moving one horizontal step in Figure 3. To solve (2.2), we need to use LAMM iteratively, which we shall call the iterative LAMM (I-LAMM) algorithm, and compute a sequence of solutions $\boldsymbol{\beta}^{(\ell,k)}$ using the initial value $\boldsymbol{\beta}^{(\ell,k-1)}$. Figure 3 depicts the schematics of our algorithm: the $\ell$th row corresponds to solving the $\ell$th subproblem in (2.2) approximately, beginning by computing the adaptive weight $\boldsymbol{\lambda}^{(\ell-1)}$. The number of iterations needed within each row will be discussed in the sequel.

2.2. *Stopping Criterion.* I-LAMM recognizes that the exact solutions to (2.2) can never be achieved in practice with algorithmic complexity control.

Instead, in the $\ell$th optimization subproblem, we compute the approximate solution, $\widetilde{\boldsymbol{\beta}}^{(\ell)}$, up to an optimization error $\varepsilon$, the choice of which will be discussed in next subsection. To calculate this approximate solution, starting from the initial value $\boldsymbol{\beta}^{(\ell,0)} = \widetilde{\boldsymbol{\beta}}^{(\ell-1)}$, the algorithm constructs a solution sequence $\{\boldsymbol{\beta}^{(\ell,k)}\}_{k=1,2,\cdots}$ using the introduced LAMM algorithm. See Figure 3.

We then introduce a stopping criterion for the I-LAMM algorithm. From optimization theory (Section 5.5 in [5]), we know that any exact solution $\widehat{\boldsymbol{\beta}}^{(\ell)}$ to the $\ell$th subproblem in (2.2) satisfies the first order optimality condition:

$$(2.6) \quad \nabla \mathcal{L}(\widehat{\boldsymbol{\beta}}^{(\ell)}) + \boldsymbol{\lambda}^{(\ell-1)} \odot \boldsymbol{\xi} = \mathbf{0}, \text{ for some } \boldsymbol{\xi} \in \partial \|\widehat{\boldsymbol{\beta}}^{(\ell)}\|_1 \in [-1,1]^d,$$

where $\partial$ is used to indicate the subgradient operator. The set of subgradients of a function $f : \mathbb{R}^d \to \mathbb{R}$ at a point $x_0$, denoted as $\partial f(x_0)$, is defined as the collection of vectors, $\boldsymbol{\xi}$, such that $f(x) - f(x_0) \geq \boldsymbol{\xi}^{\mathrm{T}}(x - x_0)$, for any $x$. Thus, a natural measure for suboptimality of $\boldsymbol{\beta}$ can be defined as

$$\omega_{\boldsymbol{\lambda}^{(\ell-1)}}(\boldsymbol{\beta}) = \min_{\boldsymbol{\xi} \in \partial \|\boldsymbol{\beta}\|_1} \{\|\nabla \mathcal{L}(\boldsymbol{\beta}) + \boldsymbol{\lambda} \odot \boldsymbol{\xi}\|_\infty\}.$$

For a prefixed optimization error $\varepsilon$, we stop the algorithm within the $\ell$th subproblem when $\omega_{\boldsymbol{\lambda}^{(\ell-1)}}(\boldsymbol{\beta}^{(\ell,k)}) \leq \varepsilon$. We call $\widetilde{\boldsymbol{\beta}}^{(\ell)} \equiv \boldsymbol{\beta}^{(\ell,k)}$ an $\varepsilon$-optimal solution. More details can be found in Algorithm 2.

---

**Algorithm 2** I-LAMM algorithm for each subproblem in (2.2).

1: **Algorithm:** $\{\widetilde{\boldsymbol{\beta}}^{(\ell)}\} \leftarrow$ I-LAMM$(\boldsymbol{\lambda}^{(\ell-1)}, \boldsymbol{\beta}^{(\ell,0)})$
2: **Input:** $\phi_0 > 0$
3: **for** $k = 0, 1, \cdots$ until $\omega_{\boldsymbol{\lambda}^{(\ell-1)}}(\boldsymbol{\beta}^{(\ell,k)}) \leq \varepsilon$ **do**
4:     $\{\boldsymbol{\beta}^{(\ell,k)}, \phi^{(\ell,k)}\} \leftarrow$ LAMM$(\boldsymbol{\lambda}^{(\ell-1)}, \boldsymbol{\beta}^{(\ell,k-1)}, \phi_0)$
5: **end for**
6: **Output:** $\widetilde{\boldsymbol{\beta}}^{(\ell)} = \boldsymbol{\beta}^{(\ell,k)}$

---

$$\begin{aligned}
\boldsymbol{\lambda}^{(0)}: & \quad \boldsymbol{\beta}^{(1,0)} = \mathbf{0} \xRightarrow{\text{LAMM}} \boldsymbol{\beta}^{(1,1)} \xRightarrow{\text{LAMM}} \cdots \xRightarrow{\text{LAMM}} \boldsymbol{\beta}^{(1,k_1)} = \widetilde{\boldsymbol{\beta}}^{(1)}, \ k_1 \lesssim \varepsilon_c^{-2}; \\
\boldsymbol{\lambda}^{(1)}: & \quad \boldsymbol{\beta}^{(2,0)} = \widetilde{\boldsymbol{\beta}}^{(1)} \xRightarrow{\text{LAMM}} \boldsymbol{\beta}^{(2,1)} \xRightarrow{\text{LAMM}} \cdots \xRightarrow{\text{LAMM}} \boldsymbol{\beta}^{(2,k)} = \widetilde{\boldsymbol{\beta}}^{(2)}, \ k_2 \lesssim \log(\varepsilon_t^{-1}); \\
& \quad \vdots \qquad\qquad\qquad\qquad \vdots \qquad\qquad\qquad\qquad \vdots \\
\boldsymbol{\lambda}^{(T-1)}: & \quad \boldsymbol{\beta}^{(T,0)} = \widetilde{\boldsymbol{\beta}}^{(T-1)} \xRightarrow{\text{LAMM}} \boldsymbol{\beta}^{(T,1)} \xRightarrow{\text{LAMM}} \cdots \xRightarrow{\text{LAMM}} \boldsymbol{\beta}^{(T,k)} = \widetilde{\boldsymbol{\beta}}^{(T)}, \ k_T \lesssim \log(\varepsilon_t^{-1}).
\end{aligned}$$

FIG 3. *Paradigm illustration of I-LAMM. $k_\ell, 1 \leq \ell \leq T$, is the iteration index for the $\ell$th optimization in (2.2). $\varepsilon_c$ and $\varepsilon_t$ are the precision parameters for the contraction and tightening stage respectively and will be described in Section 2.3 in detail.*

**Remark 2.1.** The I-LAMM algorithm is an early-stop variant of the ISTA algorithm to handle general loss functions and nonconvex penalties [2]. The LAMM principal serves as a novel perspective for the proximal gradient method.

2.3. *Tightening After Contraction.* From the computational perspective, optimization in (2.2) can be categorized into two stages: contraction ($\ell = 1$) and tightening ($2 \leq \ell \leq T$). In the contraction stage, we start from an arbitrary initial value, which can be quite remote from the underlying true parameter. We take $\varepsilon$ as $\varepsilon_c \asymp \lambda$, reflecting the precision needed to bring the initial solution to a contracting neighborhood of the global minimum. For instance, in linear model with sub-Gaussian errors, $\varepsilon_c$ can be taken in the order of $\sqrt{\log d/n}$. This stage aims to find a good initial estimator $\widetilde{\boldsymbol{\beta}}^{(1)}$ for the subsequent optimization subproblems in the tightening stage. Recall that $s = \|\boldsymbol{\beta}^*\|_0$ is the sparsity level. We will show in section 4.3 that with a properly chosen $\lambda$, the approximate solution $\widetilde{\boldsymbol{\beta}}^{(1)}$, produced by the early stopped I-LAMM algorithm, falls in the region of such good initials estimators
$$\big\{\boldsymbol{\beta} : \|\boldsymbol{\beta} - \boldsymbol{\beta}^*\|_2 \leq C\lambda\sqrt{s} \text{ and } \boldsymbol{\beta} \text{ is sparse}\big\}.$$
We call this region the contraction region.

However, the estimator $\widetilde{\boldsymbol{\beta}}^{(1)}$ suffers from a suboptimal statistical rate of convergence which is inferior to the refined one obtained by nonconvex regularization. A second stage to tighten this coarse contraction estimator into the optimal region of convergence is needed. This is achieved by the subsequent optimization ($\ell \geq 2$) and referred to as a tightening stage. Because the initial estimators are already good and sparse at each iteration of the tightening stage, the I-LAMM algorithm at this stage enjoys geometric rate of convergence, due to the sparse strong convexity. Therefore, the optimization error $\varepsilon = \varepsilon_t$ can be much smaller to simultaneously ensure statistical accuracy and control computational complexity. To achieve the oracle rate $\sqrt{s/n}$: $\varepsilon_t$ must be no larger than the order of $\sqrt{1/n}$. A graphical illustration of the full algorithm is presented in Figure 3. Theoretical justifications are provided in Section 4. From this perspective, we shall also call the psuedo-algorithm in (1.2) or (2.2), combined with LAMM, the tightening after contraction (TAC) algorithm.

## 3. New Insights into Existing Methods.

3.1. *Connection to One-step Local Linear Approximation.* In the low dimensional regime, [34] shows that the one-step LLA algorithm produces an oracle estimator if the maximum likelihood estimator (MLE) is used for initialization. They thus claim that the multi-step LLA is unnecessary. However, this is not the case in high dimensions, under which an unbiased initial estimator, such as the MLE, is not available. In this paper, we show that starting from a possibly arbitrary bad initial value (such as $\mathbf{0}$), the contraction stage can produce a sparse coarse estimator. Each tightening step then refines the estimator from previous step to the optimal region of convergence by

$$(3.1) \qquad \|\widetilde{\boldsymbol{\beta}}^{(\ell)} - \boldsymbol{\beta}^*\|_2 \lesssim \sqrt{\frac{s}{n}} + \delta \cdot \|\widetilde{\boldsymbol{\beta}}^{(\ell-1)} - \boldsymbol{\beta}^*\|_2, \text{ for } 2 \leq \ell \leq T,$$

where $\delta \in (0,1)$ is a prefixed contraction parameter. Unlike the one-step method in [12], the role of iteration is clearly evidenced in (3.1).

An important aspect of our algorithm (2.2) is that we use the solvable approximate solutions, $\widetilde{\boldsymbol{\beta}}^{(\ell)}$'s, rather than the exact ones, $\widehat{\boldsymbol{\beta}}^{(\ell)}$'s. In order to practically implement (2.2) for a general convex loss function, [34] propose to locally approximate $\mathcal{L}(\boldsymbol{\beta})$ by a quadratic function:

$$(3.2) \quad \mathcal{L}(\widehat{\boldsymbol{\beta}}^{(0)}) + \langle \nabla \mathcal{L}(\widehat{\boldsymbol{\beta}}^{(0)}), \boldsymbol{\beta} - \widehat{\boldsymbol{\beta}}^{(0)}\rangle + \frac{1}{2}(\boldsymbol{\beta} - \widehat{\boldsymbol{\beta}}^{(0)})^{\mathrm{T}} \nabla^2 \mathcal{L}(\widehat{\boldsymbol{\beta}}^{(0)})(\boldsymbol{\beta} - \widehat{\boldsymbol{\beta}}^{(0)}),$$

where $\widehat{\boldsymbol{\beta}}^{(0)}$ is a 'good' initial estimator of $\boldsymbol{\beta}^*$ and $\nabla^2 \mathcal{L}(\widehat{\boldsymbol{\beta}}^{(0)})$ is the Hessian evaluated at $\widehat{\boldsymbol{\beta}}^{(0)}$. However, in high dimensions, evaluating the $d \times d$ Hessian is not only computationally intensive but also requires a large storage cost. In addition, the optimization problem (2.2) can not be solved analytically with approximation (3.2). We resolve these issues by proposing the isotropic quadratic approximation, see Section 2.

3.2. *New Insight into Folded-concave Regularization and Adaptive Lasso.* The adaptive local linear approximation (2.2) provides new insight into folded-concave regularization and adaptive Lasso. To correct the Lasso's estimation bias, folded-concave regularization [9] and its one-step implementation, adaptive Lasso [33, 34, 12], have drawn much research interest due to their attractive statistical properties. For a general loss function $\mathcal{L}(\boldsymbol{\beta})$, the adaptive Lasso solves

$$\widehat{\boldsymbol{\beta}}_{\mathrm{adapt}} = \operatorname*{argmin}_{\boldsymbol{\beta}} \left\{ \mathcal{L}(\boldsymbol{\beta}) + \lambda \sum_{j=1}^{d} \mathrm{w}(\beta_{\mathrm{init},j})|\beta_j| \right\},$$

where $\beta_{\mathrm{init},j}$ is an initial estimator of $\beta_j$. We see that the adaptive lasso is a special case of (2.2) with $\ell = 2$. Two important open questions for adaptive

Lasso are to obtain a good enough initial estimator in high dimensions and to select a suitable tuning parameter $\lambda$ which achieves the optimal statistical performance. Our solution to the first question is to use the approximate solution to Lasso with controlled computational complexity, which corresponds to $\ell = 1$ in (2.2). For the choice of $\lambda$, [7] suggested sequential tuning: in the first stage, they use cross validation to select the initial tuning parameter, denoted here by $\widehat{\lambda}_{\text{init,cv}}$ and the corresponding estimator $\widehat{\boldsymbol{\beta}}_{\text{init}}$; in the second stage, they again adopt cross validation to select the adaptive tuning parameter $\lambda$ in the adaptive Lasso. Despite the popularity of such tuning procedure, there are no theoretical guarantees to support it. As will be shown later in Theorem 4.2 and Corollary 4.3, our framework produces optimal solution by only tuning $\boldsymbol{\lambda}^{(0)} = \lambda \mathbf{1}$ in the contraction stage, indicating that sequential tuning may not be necessary for the adaptive Lasso if $\text{w}(\cdot)$ is chosen from the tightening function class $\mathcal{T}$.

It is worth noting that a classical weight $\text{w}(\beta_j) \equiv 1/|\beta_j|$ for the adaptive Lasso does not belong to the tightening function class $\mathcal{T}$. As pointed out by [10], zero is an absorbing state of the adaptive Lasso with this choice of weight function. Hence, when the Lasso estimator in the first stage misses any true positives, it will be missed forever in later stages as well. In contrast, the proposed tightening function class $\mathcal{T}$ overcomes such shortcomings by restricting the weight function $\text{w}(\cdot)$ to be bounded. This phenomenon is further elaborated via our numerical experiments in Section 6. The mean square error for the adaptive Lasso can be even worse than the Lasso estimator because the adaptive Lasso may miss true positives in the strongly correlated design case.

Our framework also reveals interesting connections between the adaptive Lasso and folded-concave regularization. Specifically, consider the following folded-concave penalized regression

$$(3.3) \quad \min_{\boldsymbol{\beta} \in \mathbb{R}^d} \left\{ \mathcal{L}(\boldsymbol{\beta}) + \mathcal{R}_\lambda(|\boldsymbol{\beta}|) \right\}, \text{ where } \mathcal{R}_\lambda(|\boldsymbol{\beta}|) \text{ is a folded concave penalty.}$$

We assume that $\mathcal{R}_\lambda(\cdot)$ is elementwisely decomposable, that is $\mathcal{R}_\lambda(|\boldsymbol{\beta}|) = \sum_{k=1}^d p_\lambda(|\beta_k|)$. Under this assumption, using the concave duality, we can rewrite $\mathcal{R}_\lambda(|\boldsymbol{\beta}|)$ as

$$(3.4) \qquad \mathcal{R}_\lambda(|\boldsymbol{\beta}|) = \inf_{\mathbf{v}} \left\{ |\boldsymbol{\beta}|^{\mathrm{T}} \mathbf{v} - \mathcal{R}_\lambda^\star(\mathbf{v}) \right\},$$

where $\mathcal{R}_\lambda^\star(\cdot)$ is the dual of $\mathcal{R}_\lambda(\cdot)$. By the duality theory, we know that the minimum of (3.4) is achieved at $\widehat{\mathbf{v}} = \nabla \mathcal{R}_\lambda(|\boldsymbol{\mu}|)|_{\boldsymbol{\mu}=\boldsymbol{\beta}}$. We can employ (3.4)

to reformulate (3.3) as

$$(\widehat{\boldsymbol{\beta}}, \widehat{\mathbf{v}}) = \underset{\boldsymbol{\beta}, \mathbf{v}}{\operatorname{argmin}} \left\{ \mathcal{L}(\boldsymbol{\beta}) + \mathbf{v}^{\mathrm{T}} |\boldsymbol{\beta}| - \mathcal{R}_\lambda^\star(\mathbf{v}) \right\}.$$

The optimization above can then be solved by exploiting the alternating minimization scheme. In particular, we repeatedly apply the following two steps:

(1) Optimize over $\boldsymbol{\beta}$ with $\mathbf{v}$ fixed: $\widehat{\boldsymbol{\beta}}^{(\ell)} = \operatorname{argmin}_{\boldsymbol{\beta}} \left\{ \mathcal{L}(\boldsymbol{\beta}) + (\widehat{\mathbf{v}}^{(\ell-1)})^{\mathrm{T}} |\boldsymbol{\beta}| \right\}$.
(2) Optimize over $\mathbf{v}$ with $\boldsymbol{\beta}$ fixed. We can obtain closed form solution: $\mathbf{v}^{(\ell)} = \nabla \mathcal{R}_\lambda(|\boldsymbol{\mu}|)|_{\boldsymbol{\mu} = \widehat{\boldsymbol{\beta}}^{(\ell)}}$.

This is a special case of (1.2) if we take $\mathrm{w}(\boldsymbol{\beta}) = \lambda^{-1} \nabla \mathcal{R}_\lambda(|\boldsymbol{\mu}|)|_{\boldsymbol{\mu} = \boldsymbol{\beta}}$ and let $\ell$ grow until convergence. Therefore, with a properly chosen weight function $\mathrm{w}(\cdot)$, our proposed algorithm bridges the adaptive Lasso and folded-concave penalized regression together under different choices of $\ell$. In Corollary 4.3, we will prove that, when $\ell$ is in the order of $\log(\lambda\sqrt{n})$, then the proposed estimator enjoys the optimal statistical rate $\|\widehat{\boldsymbol{\beta}}^{(\ell)} - \boldsymbol{\beta}^*\|_2 \propto \sqrt{s/n}$, under mild conditions.

**4. Theoretical Results.** We establish the optimal statistical rate of convergence and the computational complexity of the proposed algorithm. To establish these results in a general framework, we first introduce the localized versions of the sparse eigenvalue and restricted eigenvalue conditions.

4.1. *Localized Eigenvalues and Assumptions.* The sparse eigenvalue condition [30] is commonly used in the analysis of sparse learning problems. However, it is only valid for the least square loss. For a general loss function, the Hessian matrix depends on the parameter $\boldsymbol{\beta}$ and can become nearly singular in certain regions. For example, the Hessian matrix of the logistic loss is

$$\nabla^2 \mathcal{L}(\boldsymbol{\beta}) = \frac{1}{n} \sum_{i=1}^n \mathbf{x}_i \mathbf{x}_i^{\mathrm{T}} \cdot \frac{1}{1 + \exp(-\mathbf{x}_i^{\mathrm{T}} \boldsymbol{\beta})} \cdot \frac{1}{1 + \exp(\mathbf{x}_i^{\mathrm{T}} \boldsymbol{\beta})},$$

which tends to zero as $\|\boldsymbol{\beta}\|_2 \to \infty$, no matter what the data are. One of our key theoretical observations is that: what we really need are the localized conditions around the true parameters $\boldsymbol{\beta}^*$, which we now introduce.

4.1.1. *Localized Sparse Eigenvalue.*

**Definition 4.1** (Localized Sparse Eigenvalue, LSE)**.** The localized sparse eigenvalues are defined as

$$\rho_+(m,r) = \sup_{\mathbf{u},\boldsymbol{\beta}} \left\{ \mathbf{u}_J^T \nabla^2 \mathcal{L}(\boldsymbol{\beta}) \mathbf{u}_J : \|\mathbf{u}_J\|_2^2 = 1, |J| \leq m, \|\boldsymbol{\beta} - \boldsymbol{\beta}^*\|_2 \leq r \right\};$$

$$\rho_-(m,r) = \inf_{\mathbf{u},\boldsymbol{\beta}} \left\{ \mathbf{u}_J^T \nabla^2 \mathcal{L}(\boldsymbol{\beta}) \mathbf{u}_J : \|\mathbf{u}_J\|_2^2 = 1, |J| \leq m, \|\boldsymbol{\beta} - \boldsymbol{\beta}^*\|_2 \leq r \right\}.$$

Both $\rho_+(m,r)$ and $\rho_-(m,r)$ depend on the Hessian matrix $\nabla^2 \mathcal{L}(\boldsymbol{\beta})$, the true coefficient $\boldsymbol{\beta}^*$, the sparsity level $m$, and an extra locality parameter $r$. They reduce to the commonly-used sparse eigenvalues when $\nabla^2 \mathcal{L}(\boldsymbol{\beta})$ does not change with $\boldsymbol{\beta}$ as in the quadratic loss. The following assumption specifies the LSE condition in detail. Recall that $s = \|\boldsymbol{\beta}^*\|_0$.

**Assumption 4.1.** We say the LSE condition holds if there exist an integer $\widetilde{s} \geq cs$ for some constant $c$, $r$ and a constant $C$ such that

$$0 < \rho_* \leq \rho_-(2s + 2\widetilde{s}, r) < \rho_+(2s + 2\widetilde{s}, r) \leq \rho^* < +\infty \quad \text{and}$$
$$\rho_+(\widetilde{s}, r)/\rho_-(2s + 2\widetilde{s}, r) \leq 1 + C\widetilde{s}/s.$$

Assumption 4.1 is standard for linear regression problems and is commonly referred to as the sparse eigenvalue condition when $r = \infty$. Such conditions have been employed by [4, 24, 22, 20, 28]. The newly proposd LSE condition, to the best of our knowledge, is the weakest one in the literature.

4.1.2. *Localized Restricted Eigenvalue.* In this section, we introduce the localized version of the restricted eigenvalue condition [4]. This is an alternative condition to Assumption 4.1 that allows us to handle general Hessian matrices that depend on $\boldsymbol{\beta}$, under which the theoretical properties can be carried out parallelly.

**Definition 4.2** (Localized Restricted Eigenvalue, LRE)**.** The localized restricted eigenvalue is defined as

$$\kappa_+(m,\gamma,r) = \sup_{\mathbf{u},\boldsymbol{\beta}} \left\{ \mathbf{u}^T \nabla^2 \mathcal{L}(\boldsymbol{\beta}) \mathbf{u} : (\mathbf{u},\boldsymbol{\beta}) \in \mathcal{C}(m,\gamma,r) \right\};$$
$$\kappa_-(m,\gamma,r) = \inf_{\mathbf{u},\boldsymbol{\beta}} \left\{ \mathbf{u}^T \nabla^2 \mathcal{L}(\boldsymbol{\beta}) \mathbf{u} : (\mathbf{u},\boldsymbol{\beta}) \in \mathcal{C}(m,\gamma,r) \right\},$$

where $\mathcal{C}(m,\gamma,r) \equiv \left\{ \mathbf{u},\boldsymbol{\beta} : S \subseteq J, |J| \leq m, \|\mathbf{u}_{J^c}\|_1 \leq \gamma \|\mathbf{u}_J\|_1, \|\boldsymbol{\beta} - \boldsymbol{\beta}^*\|_2 \leq r \right\}$ is a local $\ell_1$ cone.

Similarly, the localized restricted eigenvalue reduces to the restricted eigenvalue when $\nabla^2 \mathcal{L}(\boldsymbol{\beta})$ does not depend on $\boldsymbol{\beta}$. We say the localized restricted eigenvalue condition holds if there exists $m, \gamma, r$ such that $0 < \kappa_-(m, \gamma, r) \leq \kappa_+(m, \gamma, r) < \infty$. In Appendix B, we give a geometric explanation of the local $\ell_1$ cone, $\mathcal{C}(m, \gamma, r)$, and the coresponding localized analysis.

4.2. *Statistical Theory.* In this section, we provide theoretical analysis of the proposed estimator under the LSE condition. For completeness, in Appendix B, we also establish similar results under localized restricted eigenvalue condition. We begin with the contraction stage. Recall that the initial value $\widetilde{\boldsymbol{\beta}}^{(0)}$ is taken as $\mathbf{0}$ for simplicity. We need the following assumption on the tightening function.

**Assumption 4.2.** Assume that $w(\cdot) \in \mathcal{T}$ and $w(u) \geq 1/2$ for $u = 18\rho_*^{-1}\delta^{-1}\lambda$. Here $\mathcal{T}$ is the tightening function class defined in (2.1).

Our first result characterizes the statistical convergence rate of the estimator in the contraction stage. The key ideas of the proofs are outlined in Section 5. Other technical lemmas and details can be found in the online supplement.

**Proposition 4.1** (Statistical Rate in the Contraction Stage). Suppose that Assumption 4.1 holds. If $\lambda$, $\varepsilon$ and $r$ satisfy

$$(4.1) \qquad 4(\|\nabla \mathcal{L}(\boldsymbol{\beta}^*)\|_\infty + \varepsilon) \leq \lambda \leq r\rho_*/(18\sqrt{s}),$$

then any $\varepsilon_c$-optimal solution $\widetilde{\boldsymbol{\beta}}^{(1)}$ satisfies

$$\|\widetilde{\boldsymbol{\beta}}^{(1)} - \boldsymbol{\beta}^*\|_2 \leq 18\rho_*^{-1}\lambda\sqrt{s} \lesssim \lambda\sqrt{s}.$$

The result above is a deterministic statement. Its proof is omitted as it directly follows from Lemma 5.1 with $\ell = 1$ and $\mathcal{E}_1$ there to be $S$, the support of the true parameter $\boldsymbol{\beta}^*$. The proof of Lemma 5.1 can be found in Appendix B. In Proposition 4.1, the approximation error $\varepsilon_c$, can be taken to be the order of $\lambda \asymp \sqrt{\log d/n}$ in the sub-Gaussian noise case. The contraction stage ensures that the $\ell_2$ estimation error is proportional to $\lambda\sqrt{s}$, which is identical to the optimal rate of convergence for the Lasso estimator [4, 31]. Our result can be regarded as a generalization of the usual Lasso analysis to more general losses which satisfy the localized sparse eigenvalue condition. We are ready to present the main theorem, which demonstrates the effects of optimization error, shrinkage bias and tightening steps on the statistical rate.

**Theorem 4.2** (Optimal Statistical Rate). Suppose Assumptions 4.1 and 4.2 hold. If $4(\|\nabla\mathcal{L}(\boldsymbol{\beta}^*)\|_\infty + (\varepsilon_t \vee \varepsilon_c)) \leq \lambda \lesssim r/\sqrt{s}$, then any $\varepsilon_t$-optimal solution $\widetilde{\boldsymbol{\beta}}^{(\ell)}$, $\ell \geq 2$, satisfies the following $\delta$-contraction property

$$\|\widetilde{\boldsymbol{\beta}}^{(\ell)} - \boldsymbol{\beta}^*\|_2 \leq C\big(\|\nabla\mathcal{L}(\boldsymbol{\beta}^*)_S\|_2 + \varepsilon_t\sqrt{s} + \lambda\|\mathrm{w}(|\boldsymbol{\beta}_S^*| - u)\|_2\big) + \delta\|\widetilde{\boldsymbol{\beta}}^{(\ell-1)} - \boldsymbol{\beta}^*\|_2,$$

where $C$ is a constant and $u = 18\rho_*^{-1}\delta^{-1}\lambda$. Consequently, there exists a constant $C'$ such that

$$\|\widetilde{\boldsymbol{\beta}}^{(\ell)} - \boldsymbol{\beta}^*\|_2 \leq C'\big(\underbrace{\|\nabla\mathcal{L}(\boldsymbol{\beta}^*)_S\|_2}_{\text{oracle rate}} + \overbrace{\varepsilon_t\sqrt{s}}^{\text{opt err}} + \underbrace{\lambda\|\mathrm{w}(|\boldsymbol{\beta}_S^*| - u)\|_2}_{\text{coefficient effect}}\big) + \overbrace{2C'\delta^{\ell-1}\lambda\sqrt{s}}^{\text{tightening effect}}.$$

The effect of the tightening stage can be clearly seen from the theorem above: each tightening step induces a $\delta$-contraction property, which reduces the influence of the estimation error from the previous step by a $\delta$-fraction. Therefore, in order to achieve the oracle rate $\sqrt{s/n}$, we shall carefully choose the optimization error such that $\varepsilon_t \lesssim \|\nabla\mathcal{L}(\boldsymbol{\beta}^*)\|_2/\sqrt{s}$ and make the tightening iterations $\ell$ large enough. As a corollary, we give the explicit statistical rate under the quadratic loss $\mathcal{L}(\boldsymbol{\beta}) = (2n)^{-1}\|\mathbf{y} - \mathbf{X}\boldsymbol{\beta}\|_2^2$. In this case, we take $\lambda \asymp \sqrt{n^{-1}\log d}$ so that the scaling condition (4.1) holds with high probability. We use sub-Gaussian$(0, \sigma^2)$ to denote a sub-Gaussian distribution random variable with mean 0 and variance proxy $\sigma^2$.

**Corollary 4.3.** Let $y_i = \mathbf{x}_i^T\boldsymbol{\beta}^* + \epsilon_i$, $1 \leq i \leq n$, be independently and identically distributed sub-Gaussian random variables with $\epsilon_i \sim$ sub-Gaussian$(0, \sigma^2)$. The columns of $\mathbf{X}$ are normalized such that $\max_j \|\mathbf{X}_{*j}\|_2 \leq \sqrt{n}$. Assume there exists an $\gamma > 0$ such that $\|\boldsymbol{\beta}_S^*\|_{\min} \geq u + \gamma\lambda$ and $\mathrm{w}(\gamma\lambda) = 0$. Under Assumptions 4.1 and 4.2, if $\lambda \asymp \sqrt{n^{-1}\log d}$, $\varepsilon_t \leq \sqrt{1/n}$ and $T \gtrsim \log\log d$, then with probability at least $1 - 2d^{-\eta_1} - 2\exp\{-\eta_2 s\}$, $\widetilde{\boldsymbol{\beta}}^{(T)}$ must satisfy

$$\|\widetilde{\boldsymbol{\beta}}^{(T)} - \boldsymbol{\beta}^*\|_2 \lesssim \sqrt{s/n},$$

where $\eta_1$ and $\eta_2$ are positive constants.

Corollary 4.3 indicates that I-LAMM can achieve the oracle statistical rate $\sqrt{s/n}$ as if the support for the true coefficients were known in advance. To achieve such rate, we require $\varepsilon_c \lesssim \sqrt{\log d/n}$ and $\varepsilon_t \lesssim \sqrt{1/n}$. In other words, we need only a more accurate estimator in the tightening stage rather than in both stages. This will help us to relax the computational burden, which will be discussed in detail in Theorem 4.7. Our last result concerns the oracle property of the obtained estimator $\widetilde{\boldsymbol{\beta}}^{(\ell)}$ for $\ell$ large enough, with the proof postponed to Appendix B in the online supplement. We first define the oracle estimator $\widehat{\boldsymbol{\beta}}^\circ$ as

$$\widehat{\boldsymbol{\beta}}^\circ = \operatorname*{argmin}_{\operatorname{supp}(\boldsymbol{\beta})=S} \mathcal{L}(\boldsymbol{\beta}).$$

**Theorem 4.4** (Strong Oracle Property). Suppose Assumptions 4.1 and 4.2 hold. Assume $\|\boldsymbol{\beta}_S^*\|_{\min} \geq u + \gamma\lambda$ and $\operatorname{w}(\gamma\lambda) = 0$ for some constant $\gamma$. Let $4(\|\nabla\mathcal{L}(\widehat{\boldsymbol{\beta}}^\circ)\|_\infty + \varepsilon_c \vee \varepsilon_t) \leq \lambda \lesssim r/\sqrt{s}$ and $\varepsilon_t \leq \lambda/\sqrt{s}$. If $\|\widehat{\boldsymbol{\beta}}^\circ - \boldsymbol{\beta}^*\|_{\max} \leq \eta_n \lesssim \lambda$, then for $\ell$ large enough such that $\ell \gtrsim \log\{(1+\varepsilon_c/\lambda)\sqrt{s}\}$, we have

$$\widetilde{\boldsymbol{\beta}}^{(\ell)} = \widehat{\boldsymbol{\beta}}^\circ.$$

The theorem above is again a deterministic result. Large probability bound can be obtained by bounding the probability of the event $\{4(\|\nabla\mathcal{L}(\widehat{\boldsymbol{\beta}}^\circ)\|_\infty + (\varepsilon_c \vee \varepsilon_t)) \leq \lambda\}$. The assumption that $\|\widehat{\boldsymbol{\beta}}^\circ - \boldsymbol{\beta}^*\|_{\max} \lesssim \lambda$ is very mild, because the oracle estimator only depends on the intrinsic dimension $s$ rather than $d$. For instance, under linear model with sub-Gaussian errors, it can be shown that $\|\widehat{\boldsymbol{\beta}}^\circ - \boldsymbol{\beta}\|_{\max} \leq \sqrt{\log s/n}$ with high probability.

Theorem 4.4 implies that the oracle estimator $\widehat{\boldsymbol{\beta}}^\circ$ is a fixed point of the I-LAMM algorithm, namely, once the initial estimator is $\widehat{\boldsymbol{\beta}}^\circ$, the next iteration produces the same estimator. This is in the same spirit as that proved in [12].

4.3. *Computational Theory.* In this section, we analyze the computational rate for all of our approximate solutions. We start with the following assumption.

**Assumption 4.3.** $\nabla\mathcal{L}(\boldsymbol{\beta})$ is locally $\rho_c$-Lipschitz continuous, i.e.

(4.2) $\quad \|\nabla\mathcal{L}(\boldsymbol{\beta}_1) - \nabla\mathcal{L}(\boldsymbol{\beta}_2)\|_2 \leq \rho_c \|\boldsymbol{\beta}_1 - \boldsymbol{\beta}_2\|_2, \text{ for } \boldsymbol{\beta}_1, \boldsymbol{\beta}_2 \in B_2(R/2, \boldsymbol{\beta}^*),$

where $\rho_c$ is the Lipschitz constant and $R \lesssim \|\boldsymbol{\beta}^*\|_2 + \lambda\sqrt{s}$.

We then give the explicit iteration complexity of the contraction stage in the following proposition. Recall the definition of $\phi_0$ and $\gamma_u$ in Algorithm 2.1, and $\rho_*$ in Assumption 4.1.

**Proposition 4.5** (Sublinear Rate in the Contraction Stage). Assume that Assumption 4.1 and 4.3 hold. Let $4(\|\nabla\mathcal{L}(\boldsymbol{\beta}^*)\|_\infty + \varepsilon_c) \leq \lambda \lesssim r/\sqrt{s}$. To achieve an approximate local solution $\widetilde{\boldsymbol{\beta}}^{(1)}$ such that $\omega_{\boldsymbol{\lambda}^{(0)}}(\widetilde{\boldsymbol{\beta}}^{(1)}) \leq \varepsilon_c$ in the contraction stage, we need no more than $((1 + \gamma_u)R\rho_c/\varepsilon_c)^2$ LAMM iterations, where $\rho_c$ is a constant defined in (4.2).

The sublinear rate is due to the lack of strong convexity of the loss function in the contraction stage, because we allow starting with arbitrary bad initial value, say **0**. Once it enters the contracting region (aka, the tightening stage), the problem becomes sparse strongly convex (see Proposition B.3 in Appendix B), which endows the algorithm a linear rate of convergence. This is empirically demonstrated in Figure 2. Our next proposition gives a formal statement on the geometric convergence rate for each subproblem in the tightening stage.

**Proposition 4.6** (Geometric Rate in the Tightening Stage). Suppose that the same conditions for Theorem 4.2 hold. To obtain an approximate solution $\widetilde{\boldsymbol{\beta}}^{(\ell)}$ satisfying $\omega_{\boldsymbol{\lambda}^{(\ell-1)}}(\widetilde{\boldsymbol{\beta}}^{(\ell)}) \leq \varepsilon$ in each step of the $\ell$-th tightening stage ($\ell \geq 2$), we need at most $C' \log(C'' \lambda \sqrt{s}/\varepsilon)$ LAMM iterations, where $C'$ and $C''$ are two positive constants.

Proposition 4.6 suggests that we only need to conduct a logarithmic number of LAMM iterations in each tightening step. Simply combining the computational rate in both the contraction and the tightening stages, we manage to obtain the global computational complexity.

**Theorem 4.7.** Assume that $\lambda \sqrt{s} = o(1)$. Suppose that the same conditions for Theorem 4.2 hold. To achieve an approximate solution $\widetilde{\boldsymbol{\beta}}^{(\ell)}$ such that $\omega_{\boldsymbol{\lambda}^{(0)}}(\widetilde{\boldsymbol{\beta}}^{(1)}) \leq \varepsilon_c \lesssim \lambda$ and $\omega_{\boldsymbol{\lambda}^{(k-1)}}(\widetilde{\boldsymbol{\beta}}^{(k)}) \leq \varepsilon_t \lesssim \sqrt{1/n}$ for $2 \leq k \leq T$, the total number of LAMM iterations we need is at most

$$C' \frac{1}{\varepsilon_c^2} + C''(T-1) \log\left(\frac{1}{\varepsilon_t}\right),$$

where $C'$ and $C''$ are two positive constants, and $T \asymp \log(\lambda\sqrt{n})$.

**Remark 4.8.** We complete this section with a remark on the sublinear rate in the contraction stage. Without further structures, the sublinear rate in the first stage is the best possible one for the proposed optimization procedure when $\lambda$ is held fixed. Linear rate can be achieved when we start from a sufficiently good initial value. Another strategy is to use the path-following algorithm which is developed in [28], where they gradually reduce the size of $\lambda$ to ensure the solution sequence to be sparse.

**5. Proof Strategy for Main Results.** In this section, we present the proof strategies for the main statistical and computational theorems, with technical lemmas and other details left in the supplementary material.

5.1. *Proof Strategy for Statistical Recovery Result in Section 4.2.* Proposition 4.1 indicates that the contraction estimator suffers from a suboptimal rate of convergence $\lambda\sqrt{s}$. The tightening stage helps refine the statistical rate adaptively. To suppress the noise in the $\ell$th subproblem, it is necessary to control $\min_j\{|\widetilde{\beta}_j^{(\ell-1)}| : j \in S^c\}$ in high dimensions. For this, we construct an entropy set $\mathcal{E}_\ell$ of $S$ in each tightening subproblem to bound the magnitude of $\|\boldsymbol{\lambda}_{\mathcal{E}_\ell^c}^{(\ell-1)}\|_{\min}$. The entropy set at the $\ell$th step is defined as

$$(5.1) \qquad \mathcal{E}_\ell = S \cup \left\{ j : \lambda_j^{(\ell-1)} < \lambda \mathrm{w}(u), u = 18\delta^{-1}\rho_*^{-1}\lambda \propto \lambda \right\}.$$

Under mild conditions, we will show that $|\mathcal{E}_\ell| \leq 2s$ and $\|\boldsymbol{\lambda}_{\mathcal{E}_\ell^c}^{(\ell)}\|_{\min} \geq \lambda \mathrm{w}(u) \geq \lambda/2$, which is more precisely stated in the following lemma.

**Lemma 5.1.** *Suppose that Assumption 4.1 and 4.2 hold. If $4(\|\nabla\mathcal{L}(\boldsymbol{\beta}^*)\|_\infty + \varepsilon_t \vee \varepsilon_c) \leq \lambda \lesssim r/\sqrt{s}$, we must have $|\mathcal{E}_\ell| \leq 2s$, and the $\varepsilon$-optimal solution $\widetilde{\boldsymbol{\beta}}^{(\ell)}$ satisfies*

$$\begin{aligned}\|\widetilde{\boldsymbol{\beta}}^{(\ell)} - \boldsymbol{\beta}^*\|_2 &\leq 12\rho_*^{-1}\left(\|\boldsymbol{\lambda}_S^{(\ell-1)}\|_2 + \|\nabla\mathcal{L}(\boldsymbol{\beta}^*)_{\mathcal{E}_\ell}\|_2 + \varepsilon\sqrt{|\mathcal{E}_\ell|}\right) \\ &\leq 18\rho_*^{-1}\lambda\sqrt{s} \lesssim \lambda\sqrt{s}.\end{aligned}$$

Lemma 5.1 bounds $\|\widetilde{\boldsymbol{\beta}}^{(\ell)} - \boldsymbol{\beta}^*\|_2$ in terms of $\|\boldsymbol{\lambda}_S^{(\ell-1)}\|_2$, which is further upper bounded by the order of $\lambda\sqrt{s}$. The rate $\lambda\sqrt{s}$ coincides with the convergence rate of the contraction estimator. Later, we will exploit this result in our localized analysis to secure that all the approximate solutions $\{\widetilde{\boldsymbol{\beta}}^{(\ell)}\}_{\ell=1,\ldots,T}$ fall in a local $\ell_2$-ball centered at $\boldsymbol{\beta}^*$ with radius $r \gtrsim \lambda\sqrt{s}$.

The next lemma further bounds $\|\boldsymbol{\lambda}_S^{(\ell-1)}\|_2$ using functionals of $\widetilde{\boldsymbol{\beta}}^{(\ell-1)}$, which connects the adaptive regularization parameter to the estimator from previous steps.

**Lemma 5.2.** *Assume $\mathrm{w} \in \mathcal{T}$. Let $\lambda_j^{(\ell-1)} = \lambda\mathrm{w}(|\widetilde{\beta}_j^{(\ell-1)}|)$ for $\widetilde{\boldsymbol{\beta}}^{(\ell-1)}$, then for any norm $\|\cdot\|_*$, we have*

$$\|\boldsymbol{\lambda}_S^{(\ell-1)}\|_* \leq \lambda\|\mathrm{w}(|\boldsymbol{\beta}_S^*| - u)\|_* + \lambda u^{-1}\|\boldsymbol{\beta}_S^* - \widetilde{\boldsymbol{\beta}}_S^{(\ell-1)}\|_*,$$

where $\mathrm{w}(|\boldsymbol{\beta}_S^*| - u) \equiv \left(\mathrm{w}(|\beta_j^*| - u)\right)_{j \in S}$.

Lemma 5.2 bounds the tightening weight $\boldsymbol{\lambda}^{(\ell-1)}$ in the $\ell$th subproblem by two terms. The first term describes the coefficient effects: when the coefficients are large enough (in absolute value) such that $\|\boldsymbol{\beta}^*\|_{\min} \geq u + \gamma\lambda$ and $\mathrm{w}(\gamma\lambda) = 0$, it becomes 0. The second term concerns the estimation error of the estimator from previous step. Combing the above two lemmas, we prove that $\widetilde{\boldsymbol{\beta}}^{(\ell)}$ benefits from the tightening stage and possesses a refined statistical rate of convergence. The proof of Corollary 4.3 is left in Appendix B in the online supplement.

PROOF OF THEOREM 4.2. Applying Lemma 5.1, we obtain the size of the entropy set $\mathcal{E}_\ell$ (see definition in (5.1) ) is bounded by $2s$ and

$$(5.2) \quad \|\widetilde{\boldsymbol{\beta}}^{(\ell)} - \boldsymbol{\beta}^*\|_2 \leq C_1 \Big( \|\boldsymbol{\lambda}_S^{(\ell-1)}\|_2 + \|\nabla \mathcal{L}(\boldsymbol{\beta}^*)_{\mathcal{E}_\ell}\|_2 + \varepsilon_t \sqrt{|\mathcal{E}_\ell|} \Big) \lesssim \lambda \sqrt{s},$$

where $C_1 = 12\rho_*^{-1}$. Using Lemma 5.2 yields that

$$\|\boldsymbol{\lambda}_S^{(\ell-1)}\|_2 \leq \lambda \|\mathrm{w}(|\boldsymbol{\beta}_S^*| - u)\|_2 + \lambda u^{-1} \|(\widetilde{\boldsymbol{\beta}}^{(\ell-1)} - \boldsymbol{\beta}^*)_S\|_2.$$

Plugging the inequality above into (5.2) obtains us that

$$(5.3) \quad \|\widetilde{\boldsymbol{\beta}}^{(\ell)} - \boldsymbol{\beta}^*\|_2 \leq C_1 \Big( \underbrace{\|\nabla \mathcal{L}(\boldsymbol{\beta}^*)_{\mathcal{E}_\ell}\|_2 + \varepsilon_t \sqrt{|\mathcal{E}_\ell|}}_{\mathrm{I}} + \lambda \|\mathrm{w}(|\boldsymbol{\beta}_S^*| - u)\|_2 \Big)$$
$$+ C_1 \lambda u^{-1} \|(\widetilde{\boldsymbol{\beta}}^{(\ell-1)} - \boldsymbol{\beta}^*)_S\|_2.$$

We now simplify the inequality above by providing an upper bound for term I. Decomposing the support set $\mathcal{E}_\ell$ into $S$ and $\mathcal{E}_\ell \setminus S$ and applying the triangle inequality along with the Hölder inequality, we have

$$(5.4) \quad \mathrm{I} \leq \|\nabla \mathcal{L}(\boldsymbol{\beta}^*)_S\|_2 + \varepsilon_t \sqrt{s} + \big( \|\nabla \mathcal{L}(\boldsymbol{\beta}^*)\|_\infty + \varepsilon_t \big) \sqrt{\mathcal{E}_\ell / S}.$$

Following the proof of Lemma 5.1 in Appendix B, $\sqrt{|\mathcal{E}_\ell \setminus S|}$ can be bounded by

$$\|\widetilde{\boldsymbol{\beta}}^{(\ell-1)}_{\mathcal{E}_\ell \setminus S}\|_2 / u \leq \|\widetilde{\boldsymbol{\beta}}^{(\ell-1)} - \boldsymbol{\beta}^*\|_2 / u, \text{ where } u = 18\rho^{*-1}\delta^{-1}\lambda \propto \lambda.$$

Therefore, (5.4) can be simplified to

$$\mathrm{I} \leq \|\nabla \mathcal{L}(\boldsymbol{\beta}^*)_S\|_2 + \varepsilon_t \sqrt{s} + \frac{\lambda}{4u} \|\widetilde{\boldsymbol{\beta}}^{(\ell-1)} - \boldsymbol{\beta}^*\|_2,$$

which, combining with (5.3), yields the contraction property with $\delta$. Consequently, we obtain

$$\|\widetilde{\boldsymbol{\beta}}^{(\ell)} - \boldsymbol{\beta}^*\|_2 \leq C \big( \|\nabla \mathcal{L}(\boldsymbol{\beta}^*)_{\mathcal{E}_\ell}\|_2 + \varepsilon_t \sqrt{s} + \lambda \|\mathbf{w}_S(|\boldsymbol{\beta}_S^*| - u)\|_2 \big) + \delta^{\ell-1} \|\widetilde{\boldsymbol{\beta}}^{(1)} - \boldsymbol{\beta}^*\|_2,$$
$$\leq C \big( \|\nabla \mathcal{L}(\boldsymbol{\beta}^*)_{\mathcal{E}_\ell}\|_2 + \varepsilon_t \sqrt{s} + \lambda \|\mathbf{w}_S(|\boldsymbol{\beta}_S^*| - u)\|_2 \big) + C\delta^{\ell-1}\lambda\sqrt{s},$$

where $C = C_1/(1-\delta)$ and the last inequality follows from Proposition 4.1. The proof is completed. $\square$

5.2. *Proof Strategy for Computational Result in Section 4.3.* In this section, we present the sketch for the proofs of the results in Section 4.3. We start with the contraction stage. The next lemma shows that the contraction stage enjoys a sublinear rate of convergence. The proof can be found in Appendix C.

**Lemma 5.3.** Recall that $F(\boldsymbol{\beta}, \boldsymbol{\lambda}) = \mathcal{L}(\boldsymbol{\beta}) + \sum_{j=1}^{d} \lambda_j |\beta_j|$. We have

$$F(\boldsymbol{\beta}^{(1,k)}, \boldsymbol{\lambda}^{(0)}) - F(\widehat{\boldsymbol{\beta}}^{(1)}, \boldsymbol{\lambda}^{(0)}) \leq \frac{\phi_c}{2k} \|\boldsymbol{\beta}^{(1,0)} - \widehat{\boldsymbol{\beta}}^{(1)}\|_2^2.$$

The result above suggests that the optimization error decreases to zero at the rate of $1/k$, while Proposition 4.1 indicates that the best statistical rate for the contraction stage is only in the order of $\lambda\sqrt{s}$. Therefore, one can early stop the LAMM iterations in the contraction stage as soon as it enters the contraction region $\{\boldsymbol{\beta} : \|\boldsymbol{\beta} - \boldsymbol{\beta}^*\|_2 \lesssim C\lambda\sqrt{s}, \boldsymbol{\beta} \text{ is sparse}\}$. It is this lemma that helps characterize the iteration complexity in terms of the total number of LAMM updates needed in the contraction stage, see Proposition 4.5.

To utilize the localized sparse eigenvalue condition in the tightening stage, we need the following proposition which characterizes the sparsity of all the approximate solutions produced by the contraction stage.

**Lemma 5.4.** Assume that Assumption 4.1 holds. If $4(\|\nabla \mathcal{L}(\boldsymbol{\beta}^*)\|_\infty + \varepsilon_c) \leq \lambda \lesssim r/\sqrt{s}$, then $\widetilde{\boldsymbol{\beta}}^{(1)}$ in the contraction stage is $s + \widetilde{s}$ sparse. In particular, we have $\|(\widetilde{\boldsymbol{\beta}}^{(1)})_{S^c}\|_0 \leq \widetilde{s}$.

Together with Proposition 4.1, it ensures that the contraction estimator $\widetilde{\boldsymbol{\beta}}^{(1)}$ falls in the contraction region $\{\boldsymbol{\beta} : \|\boldsymbol{\beta} - \boldsymbol{\beta}^*\|_2 \leq C\lambda\sqrt{s} \text{ and } \boldsymbol{\beta} \text{ is sparse}\}$. This makes the localized sparse eigenvalue condition useful and thus makes the geometric rate of convergence possible.

**Lemma 5.5** (Geometric Rate in the Tightening Stage)**.** Under the same conditions for Theorem 4.2, for any $\ell \geq 2$, $\{\boldsymbol{\beta}^{(\ell,k)}\}$ converges geometrically,

$$\begin{aligned}
F(\boldsymbol{\beta}^{(\ell,k)}, \boldsymbol{\lambda}^{(\ell-1)}) &- F(\widehat{\boldsymbol{\beta}}^{(\ell)}, \boldsymbol{\lambda}^{(\ell-1)}) \\
&\leq \left(1 - \frac{1}{4\gamma_u \kappa}\right)^k \left\{ F(\boldsymbol{\beta}^{(\ell,0)}, \boldsymbol{\lambda}^{(\ell-1)}) - F(\widehat{\boldsymbol{\beta}}^{(\ell)}, \boldsymbol{\lambda}^{(\ell-1)}) \right\}.
\end{aligned}$$

The above result suggests that each subproblem in the tightening stage enjoys a geometric rate of convergence, which is the fastest possible rate among all first-order optimization methods under the blackbox model. Lemma 5.5 can be used to obtain the computational complexity analysis of each single step of the tightening stage, i.e., Proposition 4.6.

**6. Numerical Examples.** In this section, we evaluate the statistical performance of the proposed framework through several numerical experiments. We consider the following three examples.

**Example 6.1** (Linear Regression). In the first example, continuous responses were generated according to the model

$$(6.1) \qquad y_i = \mathbf{x}_i^\mathrm{T} \boldsymbol{\beta}^* + \epsilon_i, \text{ where } \boldsymbol{\beta}^* = (5, 3, 0, 0, -2, \underbrace{0, \ldots, 0}_{d-5})^\mathrm{T},$$

and $n = 100$. Moreover, in model (6.1), $\{\mathbf{x}_i\}_{i \in [n]}$ are generated from $N(0, \boldsymbol{\Sigma})$ distribution with covariance matrix $\boldsymbol{\Sigma}$, which is independent of $\epsilon_i \sim N(0, 1)$. We take $\boldsymbol{\Sigma}$ as a correlation matrix $\boldsymbol{\Sigma} = (\rho_{ij})$ as follows.

- Case 1: independent correlation design with $(\rho_{ij}) = \mathrm{diag}(1, \cdots, 1)$;
- Case 2: constant correlation design with $\rho_{ij} = 0.75$ if $i \neq j$; $\rho_{ij} = 1$, otherwise;
- Case 3: autoregressive correlation design with $\rho_{ij} = 0.95^{|i-j|}$.

**Example 6.2** (Logistic Regression). In the second example, independent observations with binary responses are generated according to the model

$$\mathbb{P}(y_i = 1 | \mathbf{x}_i) = \frac{\exp\{\mathbf{x}_i^\mathrm{T} \boldsymbol{\beta}^*\}}{1 + \exp\{\mathbf{x}_i^\mathrm{T} \boldsymbol{\beta}^*\}} \quad i = 1, \ldots, n,$$

where $\boldsymbol{\beta}^*$ and $\{\mathbf{x}_i\}_{i \in [n]}$ are generated in the same manner as in the case 1 of Example 6.1.

**Example 6.3** (Varying Dimensions and Sample Sizes). In this example, we continue Example 6.1 with varying dimensions and sample sizes. Specifically, we consider linear regression under autoregressive correlation design with $\rho_{ij} = 0.90^{|i-j|}$ with $d$ varying from 1000 to 3500 and $n$ varying from 100 to 500.

In the first two cases, we fix the sample size $n$ at 100 and consider $d = 1000$. We investigate the sparsity recovery and estimation properties of the I-LAMM (or TAC) estimator via numerical simulations. We compared the I-LAMM estimator with the following methods: the oracle estimator which assumes the availability of the active set $S$; the refitted Lasso (Refit) which uses a post least square refit on the selected set from Lasso; the adaptive Lasso (ALasso) estimator with weight function $\mathrm{w}(\beta_j) = 1/|\beta_j|$ proposed by [33]; the smoothly clipped absolute deviation (SCAD) estimator [9] with $a = 3.7$; and the minimax concave penalty (MCP) estimator with $a = 3$ [29]. For I-LAMM, we used the 3-fold cross-validation to select the constant $c \in$

Table 1

*The median of MSE, TP, FP, Time in seconds under the Case 1, Case 2 and Case 3 for linear regression in Example 6.1 and logistic regression in Example 6.2.*

|  | MSE | TP | FP | Time | MSE | TP | FP | Time |
|---|---|---|---|---|---|---|---|---|
|  | Linear\Case 1 | | | | Linear\Case 2 | | | |
| I-LAMM | 0.0285 | 3.00 | 0.00 | 0.17 | 0.0659 | 3.00 | 0.00 | 0.19 |
| Lasso | 0.3114 | 3.00 | 17.00 | 0.02 | 1.3709 | 3.00 | 16.00 | 0.04 |
| Refit | 0.5585 | 3.00 | 17.00 | 0.02 | 2.1573 | 3.00 | 16.00 | 0.04 |
| ALasso | 0.4616 | 3.00 | 15.00 | 0.06 | 1.6077 | 3.00 | 13.00 | 0.08 |
| SCAD | 0.0397 | 3.00 | 0.00 | 0.21 | 0.0695 | 3.00 | 0.00 | 0.23 |
| MCP | 0.0344 | 3.00 | 0.00 | 0.17 | 0.0706 | 3.00 | 0.00 | 0.22 |
| Oralcle | 0.0258 | 3.00 | 0.00 | - | 0.0565 | 3.00 | 0.00 | - |
|  | Linear\Case 3 | | | | Logistic | | | |
| I-LAMM | 0.2819 | 3.00 | 3.00 | 0.22 | 8.94 | 3.00 | 0.00 | 0.20 |
| Lasso | 5.8061 | 2.00 | 20.00 | 0.03 | 26.92 | 3.00 | 20.00 | 0.03 |
| Refit | 2.6354 | 2.00 | 20.00 | 0.03 | 26.85 | 3.00 | 20.00 | 0.03 |
| ALasso | 4.4242 | 2.00 | 12.00 | 0.06 | 8.28 | 3.00 | 7.00 | 0.05 |
| SCAD | 14.8680 | 2.00 | 5.00 | 0.25 | 9.48 | 3.00 | 12.00 | 0.21 |
| MCP | 14.9381 | 1.00 | 1.00 | 0.18 | 11.84 | 3.00 | 3.00 | 0.22 |
| Oralcle | 0.1661 | 3.00 | 0.00 | - | 3.32 | 3.00 | 0.00 | - |

$0.5 \times \{1, 2, \ldots, 20\}$ in the tuning parameter $\lambda = c\sqrt{\log d/n}$ in the contraction stage, with regularization parameters updated automatically at later steps. We further took $\gamma_u = 2$, $\varepsilon_c = \sqrt{\log d/n}$ and $\varepsilon_t = \sqrt{1/n}$. For Lasso, we used the I-LAMM algorithm; for ALasso, sequential tuning in [7] was used: we employed 3-fold cross validation in each step with I-LAMM algorithm used; and the SCAD and MCP estimators were computed using the R package ncvreg and 3-fold cross-validation was used for tuning parameter selection.

For each simulation setting, we generated 100 simulated datasets and applied different estimators to each dataset. We report different statistics for each estimator in Table 1. To measure the sparsity recovery performance, we calculated the median of the number of zero coefficients incorrectly estimated to be nonzero (i.e. false positive, denoted as FP), the median of the number of nonzero coefficients correctly estimated to be nonzero (i.e. true positive, denoted by TP). To measure the estimation accuracy, we calculated the median of mean squared error (MSE). To evaluate the computational efficiency, we gave the median of time (in seconds) used to produce the final estimator for different methods. Note that the computational time provided here is merely for a reference. They depend on optimization errors and implementation.

We have several important observations. First, it is not surprising that

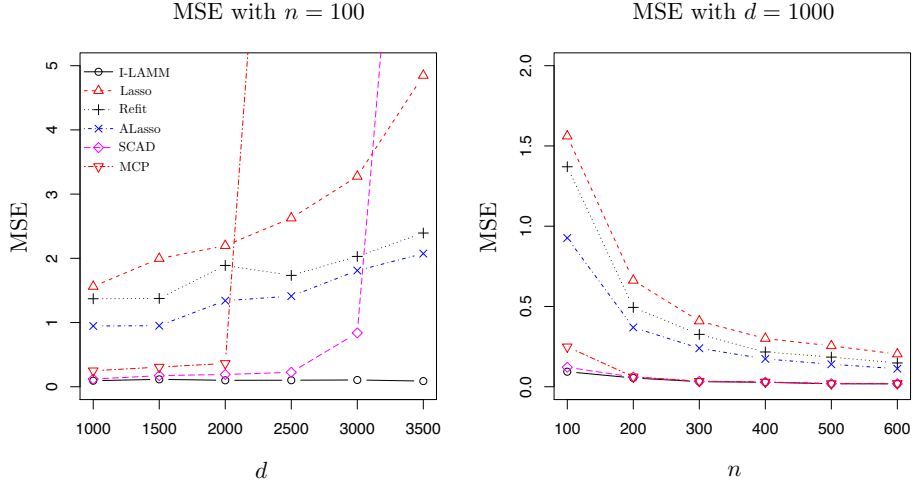

FIG 4. *The median of MSE with varying dimensions and sample sizes in Example 6.3.*

Lasso tends to overfit. Other procedures improve the performance of Lasso by reducing the estimation bias and the false positive rate. The best overall performance is achieved by the I-LAMM estimator with small MSE and FP in all cases. The MCP and SCAD estimators also have overall good performance in the logistic regression model, case 1 and case 2 of the linear regression model. However, all of MCP, SCAD and ALasso breaks down by missing true positives in case 3, where the design matrix exhibits a strong correlation between features, while I-LAMM remains the best followed by the Lasso estimator. This suggests the superiority of I-LAMM over other nonconvex penalized regression methods under strongly correlated designs. The MSE of the I-LAMM estimator keeps flat when the dimension $d$ varies, which justifies the oracle rate $\sqrt{s/n}$. SCAD and MCP have competitive performance when the dimension is relatively small, but they quickly break down when the dimension gets larger. This is possibly due to the numerical instability for directly solving nonconvex systems. This phenomenon is also observed in [28]. When the sample size is increasing, the performances of I-LAMM, SCAD and MCP are almost identical to each other while other convex methods suffer from slightly worse performance.

In addition, to demonstrate the phase transition phenomenon, in Figure 2, we plot the log estimation error verses the number of iterations for each tightening step for case 2 in Example 6.1. Indeed, the contraction stage suffers a sublinear rate of convergence before getting into the contracting region and enjoys a geometric rate afterwards, while the tightening stage has a geometric rate of convergence. These are in line with our asymptotic theory.

**7. Conclusions and Discussions.** We propose a computational framework, I-LAMM (or TAC), for simultaneous control of algorithmic complexity and statistical error when fitting high dimensional models. Even though I-LAMM only solves a sequence of convex programs approximately, the solution enjoys the same optimal statistical property of the unobtainable global optimum for the folded-concave penalized regression. Our theoretical treatment relies on a novel localized analysis which avoids the parameter bound contraint, such as $\|\boldsymbol{\beta}\|_1 \leq R$, used in all other recent works. Statistically, a $\delta$-contraction property is established: each convex program contracts the previous estimator by a $\delta$-fraction until the optimal statistical error is reached. Computationally, a phase transition in algorithmic convergence is established. The contraction stage enjoys only a sublinear rate of convergence while the tightening stage converges geometrically fast.

Recently, [22] proposed the restricted eigenvalue condition for unified M-estimators. [18] leveraged this condition, which is more related to our localized conditions. However, there are two major differences. First, their local parameter $r$ is fixed at a constant independent of $n, d, s$, while we allow it to go to 0 as long as $r \gtrsim \sqrt{s \log d/n}$. Second, their high dimensional regression problem relies on the $\ell_1$ ball constraint $\|\boldsymbol{\beta}\|_1 \leq R$, while our newly developed localized analysis, together with the localized conditions, removes such type of constraint. In [21], the authors only consider the solutions in a local cone, which makes their analysis much simpler than ours. In this paper, we provide a stronger result: with high probability, all local solutions must fall in a local sparse (or $\ell_1$) cone and thus makes the localized eigenvalue conditions applicable.

More recently, [27] proposed a two-step approach named calibrated CCCP which achieve strong oracle properties when using the Lasso estimator as initialization. Our work differs from theirs in two aspects. First, their work aims at analyzing the least square loss while our analysis handles much broader families of loss functions. Second, their procedure attains an oracle rate but requires the minimum signal strength to be in the order of $s\sqrt{\log d/n}$. Such a requirement is suboptimal. In contrast, our results requires only $\sqrt{\log d/n}$. This weakened assumption on minimum signal strength also distinguishes I-LAMM from other convex procedures, such as least squares refit after model selection [3]. In [27], the authors also proposed a high dimensional BIC criterion for variable selection and finding the oracle estimator along the solution path. We believe such a criterion can also be applied to our framework under general conditions. In further studies, [20], [19] and [28] study the theoretical properties of nonconvex penalized M-estimators. Specifically, [20] and [19] provide conditions under which all the local optima obtained by an $\ell_1$-ball constrained optimization enjoys desired statistical rates. [28] propose a path-following strategy to obtain optimal computational and statistical rates of convergence, which also relies an extra ball constraint.

Our work differs from the aforementioned literature at least in three aspects:

(1) Our theory exploits new notion of localized analysis, which is not available in [20], [19] and [28]. Such analysis allows us to eliminate the extra ball constraints in previous work, which introduce more tuning effort and are intuitively redundant given the penalty function.
(2) Our statistical results tolerate explicit computational precisions and are valid for all obtained approximate solutions, while the analysis in [20] only targets on the exact local solutions. Moreover, our computational result does not rely on the path-following type strategy as in [28] and is valid for any algorithm with desired statistical properties as basic building blocks within each of the tightening steps.
(3) We provide a refined oracle statistical rate $\sqrt{s/n}$ for the obtained approximation solution, while [20] and [28] do not provide such a result. [20] provide a statistical rate which is also achievable using the convex Lasso penalty. [28] only prove the oracle rate for exact local solutions.

Our work can be applied to many different topics: low rank matrix completion problems, high dimensional graphical models, quantile regression and many others. We conjecture that in all of the aforementioned topics, I-LAMM can give faster rate by approximately solving a sequence of convex programs, with controlled computing resources. It is also interesting to see how our algorithm works in large-scale distributed systems. Is there any fundamental tradeoffs between statistical efficiency, communication and time complexity? We leave these as future research projects.

**Supplementary Material.** The supplementary material contains proofs for Corollary 4.3, Theorem 4.4, Proposition 4.5, Proposition 4.6 and Theorem 4.7 in Section 4. It collects proofs of the lemmas presented in Section 5. An application to robust linear regression is given in Appendix D. Other technical lemmas are collected in Appendices E and F.

## References.

[1] AGARWAL, A., NEGAHBAN, S. and WAINWRIGHT, M. J. (2012). Fast global convergence rates of gradient methods for high-dimensional statistical recovery. *The Annals of Statistics* **40** 2452–2482.

[2] BECK, A. and TEBOULLE, M. (2009). A fast iterative shrinkage-thresholding algorithm for linear inverse problems. *SIAM journal on imaging sciences* **2** 183–202.

[3] BELLONI, A. and CHERNOZHUKOV, V. (2013). Least squares after model selection in high-dimensional sparse models. *Bernoulli* **19** 521–547.

[4] BICKEL, P. J., RITOV, Y. and TSYBAKOV, A. B. (2009). Simultaneous analysis of lasso and dantzig selector. *The Annals of Statistics* **37** 1705–1732.

[5] BOYD, S. and VANDENBERGHE, L. (2009). *Convex optimization.* Cambridge university press.

[6] BREHENY, P. and HUANG, J. (2011). Coordinate descent algorithms for nonconvex penalized regression, with applications to biological feature selection. *The Annals of Applied Statistics* **5** 232–253.

[7] BÜHLMANN, P. and VAN DE GEER, S. (2011). *Statistics for high-dimensional data: methods, theory and applications.* Springer Science & Business Media.

[8] BUNEA, F., TSYBAKOV, A., WEGKAMP, M. ET AL. (2007). Sparsity oracle inequalities for the lasso. *Electronic Journal of Statistics* **1** 169–194.

[9] FAN, J. and LI, R. (2001). Variable selection via nonconcave penalized likelihood and its oracle properties. *Journal of the American Statistical Association* **96** 1348–1360.

[10] FAN, J. and LV, J. (2008). Sure independence screening for ultrahigh dimensional feature space (with discussion). *Journal of Royal Statistical Society, Series B* **70** 849–911.

[11] FAN, J. and LV, J. (2011). Nonconcave Penalized Likelihood With NP-Dimensionality. *Information Theory, IEEE Transactions on* **57** 5467–5484.

[12] FAN, J., XUE, L. and ZOU, H. (2014). Strong oracle optimality of folded concave penalized estimation. *The Annals of Statistics* **42** 819–849.

[13] FRIEDMAN, J., HASTIE, T., HÖFLING, H. and TIBSHIRANI, R. (2007). Pathwise coordinate optimization. *Annals of Applied Statistics* **1** 302–332.

[14] HUNTER, D. R. and LANGE, K. (2004). A tutorial on mm algorithms. *The American Statistician* **58** 30–37.

[15] KIM, Y., CHOI, H. and OH, H.-S. (2008). Smoothly clipped absolute deviation on high dimensions. *Journal of the American Statistical Association* **103** 1665–1673.

[16] KIM, Y. and KWON, S. (2012). Global optimality of nonconvex penalized estimators. *Biometrika* **99** 315–325.

[17] LANGE, K., HUNTER, D. R. and YANG, I. (2000). Optimization transfer using surrogate objective functions. *Journal of computational and graphical statistics* **9** 1–20.

[18] LOH, P.-L. (2015). Statistical consistency and asymptotic normality for high-dimensional robust m-estimators. *arXiv preprint arXiv:1501.00312.*

[19] LOH, P.-L. and WAINWRIGHT, M. J. (2014). Support recovery without incoherence: A case for nonconvex regularization. *arXiv preprint arXiv:1412.5632.*

[20] LOH, P.-L. and WAINWRIGHT, M. J. (2015). Regularized $m$-estimators with nonconvexity: statistical and algorithmic theory for local optima. *Journal of Machine Learning Research* **16** 559–616.

[21] LOZANO, A. C. and MEINSHAUSEN, N. (2013). Minimum distance estimation for robust high-dimensional regression. *arXiv preprint arXiv:1307.3227.*

[22] NEGAHBAN, S. N., RAVIKUMAR, P., WAINWRIGHT, M. J. and YU, B. (2012). A unified framework for high-dimensional analysis of $m$-estimators with decomposable

LOCAL ADAPTIVE MAJORIZE-MINIMIZATION 29regularizers. *Statistical Science* **27** 538–557.
[23] NESTEROV, Y. (2013). Gradient methods for minimizing composite functions. *Mathematical Programming* **140** 125–161.
[24] RASKUTTI, G., WAINWRIGHT, M. J. and YU, B. (2010). Restricted eigenvalue properties for correlated gaussian designs. *The Journal of Machine Learning Research* **11** 2241–2259.
[25] TIBSHIRANI, R. (1996). Regression shrinkage and selection via the lasso. *Journal of the Royal Statistical Society, Series B* **58** 267–288.
[26] VAN DE GEER, S. A., BÜHLMANN, P. ET AL. (2009). On the conditions used to prove oracle results for the lasso. *Electronic Journal of Statistics* **3** 1360–1392.
[27] WANG, L., KIM, Y. and LI, R. (2013). Calibrating non-convex penalized regression in ultra-high dimension. *The Annals of Statistics* **41** 2505–2682.
[28] WANG, Z., LIU, H. and ZHANG, T. (2014). Optimal computational and statistical rates of convergence for sparse nonconvex learning problems. *The Annals of statistics* **42** 2164.
[29] ZHANG, C.-H. (2010). Nearly unbiased variable selection under minimax concave penalty. *The Annals of Statistics* **38** 894–942.
[30] ZHANG, C.-H. and ZHANG, T. (2012). A general theory of concave regularization for high-dimensional sparse estimation problems. *Statistical Science* **27** 576–593.
[31] ZHANG, T. (2009). Some sharp performance bounds for least squares regression with l1 regularization. *The Annals of Statistics* **37** 2109–2144.
[32] ZHANG, T. (2010). Analysis of multi-stage convex relaxation for sparse regularization. *The Journal of Machine Learning Research* **11** 1081–1107.
[33] ZOU, H. (2006). The adaptive lasso and its oracle properties. *Journal of the American Statistical Association* **101** 1418–1429.
[34] ZOU, H. and LI, R. (2008). One-step sparse estimates in nonconcave penalized likelihood models. *The Annals of statistics* **36** 1509.

# Supplementary Material to "I-LAMM: Simultaneous Control of Algorithmic Complexity and Statistical Error"

BY Jianqing Fan    Han Liu    Qiang Sun    Tong Zhang

The supplementary material contains proofs for Corollary 4.3, Theorem 4.4, Proposition 4.5, Proposition 4.6 and Theorem 4.7 in Section 4. It collects the proofs of the key lemmas presented in Section 5. An application to robust linear regression is given in Appendix D. Other technical lemmas are collected in Appendices E and F.

## APPENDIX A: GENERAL CONVEX LOSS FUNCTIONS

We request the convex loss $\mathcal{L}$ to have continuous first order derivative. In addition, we request it to be locally twice differentiable almost everywhere. Specifically, we consider the following family of loss functions

$$\mathcal{F}_\mathcal{L} = \Big\{ \mathcal{L} : \mathcal{L} \text{ is convex}, \nabla \mathcal{L} \text{ is continuous and differentiable in } B_2(r, \boldsymbol{\beta}^*) \Big\},$$

where $\boldsymbol{\beta}$ is any vector in $\mathbb{R}^d$ and $r \gtrsim \lambda\sqrt{s}$. This family includes many interesting loss functions. Some examples are given as below:

- **Logistic Loss**: Let $\{(\mathbf{x}_i, y_i)\}_{i \in [n]}$ be $n$ observed data points of $(\boldsymbol{X}, y)$, where $\boldsymbol{X}$ is a $d$-dimensional covariate and $y_i \in \{-1, +1\}$. The logistic loss is given in the form

$$\mathcal{L}(\boldsymbol{\beta}) = n^{-1} \sum_{i=1}^{n} \Big\{ \log\big(1 + \exp(-y_i \mathbf{x}_i^\mathrm{T} \boldsymbol{\beta})\big) \Big\},$$

  where $n$ is the sample size.

- **Huber Loss**: In the robust linear regression, the Huber loss takes the form

$$\mathcal{L}(\boldsymbol{\beta}) = n^{-1} \sum_{i=1}^{n} \Big\{ \ell_\alpha(y_i - \mathbf{x}_i^\mathrm{T} \boldsymbol{\beta}) \Big\},$$

  where $\ell_\alpha(x) = 2\alpha^{-1}|x| - \alpha^{-2}$ when $|x| > \alpha^{-1}$ and $\ell_\alpha(x) = x^2$ otherwise.

- **Gaussian Graphical Model Loss**: Let $\boldsymbol{\Theta} = \boldsymbol{\Sigma}^{-1}$ be the precision matrix and $\widehat{\boldsymbol{\Sigma}}$ be the sample covariance matrix. The negative log-likelihood loss of Gaussian graphical model is

$$\mathcal{L}(\boldsymbol{\Theta}) = \mathrm{tr}(\widehat{\boldsymbol{\Sigma}}\boldsymbol{\Theta}) - \log \det(\boldsymbol{\Theta}).$$

We note that the Huber loss is convex and has continuous first-order derivative. In addition, its second derivative exists in a neighborhood of $\boldsymbol{\beta}^*$. Therefore it belongs to the family of loss functions defined above. Other loss functions include the least square loss and locally twice differentiable convex composite likelihood loss.

We remark here that, the loss functions analyzed in our paper can be nonconvex. If $\mathcal{L}(\boldsymbol{\beta})$ is nonconvex, we can decompose it as $\mathcal{L}(\boldsymbol{\beta}) = \widetilde{\mathcal{L}}(\boldsymbol{\beta}) + \mathcal{H}(\boldsymbol{\beta})$ such that $\widetilde{\mathcal{L}}(\boldsymbol{\beta})$ is the convex part and $\mathcal{H}(\boldsymbol{\beta})$ is the concave part. We then write the objective function as $\mathcal{F}(\boldsymbol{\beta}) = \widetilde{\mathcal{L}}(\boldsymbol{\beta}) + \widetilde{\mathcal{R}}_\lambda(\boldsymbol{\beta})$ such that $\widetilde{\mathcal{R}}(\boldsymbol{\beta}) = \mathcal{H}(\boldsymbol{\beta}) + \mathcal{R}(\boldsymbol{\beta})$ and treat $\widetilde{\mathcal{R}}(\boldsymbol{\beta})$ as our new regularizer. If the corresponding weight function $\widetilde{\mathrm{w}}(\cdot)$ satisfies Assumption 4.1, our theory shall go through without any problems. A similar technique is exploited in [28].

## APPENDIX B: PROOFS OF STATISTICAL THEORY

**B.1. Statistical Theory under LSE Condition.** In this section, we collect the proofs for Corollary 4.3 and Theorem 4.4. We give proofs for the key technical lemmas in Section 5.1, which are used to prove theorem 4.2. Other technical lemmas are postponed to later sections. We then establish parallel results under the localized restricted eigenvalue condition.

B.1.1. *Proofs of Main Results.* In this section, we first prove Corollary 4.3 and then give the proof of Theorem 4.4.

PROOF OF COROLLARY 4.3. We start by bounding $\mathbb{P}(\|\nabla\mathcal{L}(\boldsymbol{\beta}^*)\|_\infty \geq \lambda/8)$, where $\nabla\mathcal{L}(\boldsymbol{\beta}^*) = n^{-1}\mathbf{X}^\mathrm{T}(\mathbf{y} - \mathbf{X}\boldsymbol{\beta}^*)$. For $\lambda \geq c\sqrt{\log d/n}$, using the union bound, we obtain

$$\mathbb{P}\Big(\|\nabla\mathcal{L}(\boldsymbol{\beta}^*)\|_\infty \geq \lambda/8\Big) \leq \mathbb{P}\Big(n^{-1}\|\mathbf{X}^\mathrm{T}(\mathbf{y} - \mathbf{X}\boldsymbol{\beta}^*)\|_\infty \geq 8^{-1}c\sqrt{\log d/n}\Big)$$

$$\text{(B.1)} \qquad \leq \sum_{j=1}^d \mathbb{P}\Big(1/n|\boldsymbol{X}_j^\mathrm{T}\boldsymbol{\epsilon}| \geq 8^{-1}c\sqrt{\log d/n}\Big).$$

Let $v_j = \boldsymbol{X}_j^\mathrm{T}\boldsymbol{\epsilon}$. Since $\epsilon_i$ is sub-Gaussian$(0, \sigma^2)$ for $i = 1, \ldots, n$, we obtain

$$\mathbb{E}\Big(\exp\{t_0 v_j\} + \exp\{-t_0 v_j\}\Big) \leq 2\exp\Big\{n^{-2}\|\mathbf{X}_{*j}\|^2 \sigma^2 t_0^2/2\Big\},$$

which implies $\mathbb{P}(|v_j| \geq t)\exp\{t_0 t\} \leq 2\exp\big\{n^{-2}\|\mathbf{X}_{*j}\|^2\sigma^2 t_0^2/2\big\}$. Taking $t_0 = t(n^{-2}\|\mathbf{X}_{*j}\|_2^2 \sigma^2)^{-1}$ yields that

$$\mathbb{P}(|v_j| \geq t) \leq 2\exp\bigg\{-\frac{t^2}{2\sigma^2\|\mathbf{X}_{*j}\|_2^2/n^2}\bigg\}.$$

Further taking $t = \lambda/8$ in the bound above and plugging it into (B.1) results

$$\mathbb{P}\Big(\|\nabla\mathcal{L}(\boldsymbol{\beta}^*)\|_\infty \geq \lambda/8\Big) \leq 2d\exp\bigg\{-\frac{(c/8)^2\log d}{2\sigma^2\max_j\{\|\mathbf{X}_{*j}\|_2^2/n\}}\bigg\} \leq 2d^{-c^2/(128\sigma^2)}.$$

Define the event $\mathcal{J}_1 = \Big\{\|\nabla\mathcal{L}(\boldsymbol{\beta}^*)\|_\infty \leq \lambda/8\Big\}$. Then with probability at least $1 - 2d^{-\eta_1}$ where $\eta_1 = c^2/(128\sigma^2)$, we have $\|\nabla^2\mathcal{L}(\boldsymbol{\beta}^*)\|_\infty \leq \lambda/8$.

It remains to show the oracle rate holds. Applying Theorem 4.2, we have

$$(\text{B.2}) \quad \|\widetilde{\boldsymbol{\beta}}^{(\ell)} - \boldsymbol{\beta}^*\|_2 \leq C\Big(\|\nabla\mathcal{L}(\boldsymbol{\beta}^*)_S\|_2 + \underbrace{\varepsilon_t\sqrt{s}}_{\text{I}} + \underbrace{\lambda\|\text{w}(|\boldsymbol{\beta}^*_S| - u)\|_2}_{\text{II}}\Big) + \underbrace{C\delta^{\ell-1}\lambda\sqrt{s}}_{\text{III}}.$$

For I, $\varepsilon_t\sqrt{s} \leq \sqrt{s/n}$ since $\varepsilon_t \leq \sqrt{1/n}$. Because $\|\boldsymbol{\beta}^*\|_{\min} \geq u + \gamma\lambda$, we have

$$\text{w}(|\boldsymbol{\beta}^*_S| - u) \leq \text{w}(\|\boldsymbol{\beta}^*_S\|_{\min}\mathbf{1}_S - u) \leq (\text{w}(\gamma\lambda), \ldots, \text{w}(\gamma\lambda))^{\text{T}} = \mathbf{0}.$$

This implies II=0. For III, because $\ell \geq \lfloor \log\lambda\sqrt{n}(\log 1/\delta)^{-1}\rfloor + 2 \gtrsim \log\lambda\sqrt{n}$, we obtain

$$\text{III} = C\delta^{\ell-1}\lambda\sqrt{s} \leq C\delta^{(\log 1/\delta)^{-1}\log\lambda\sqrt{n}}\lambda\sqrt{s} = C\sqrt{s/n}.$$

Plugging the bounds of I, II and III back into (B.2), we have

$$\|\widetilde{\boldsymbol{\beta}}^{(\ell)} - \boldsymbol{\beta}^*\|_2 \leq C\|\nabla\mathcal{L}(\boldsymbol{\beta}^*)_S\|_2 + C\sqrt{s/n}.$$

It remains to bound $\|\nabla\mathcal{L}(\boldsymbol{\beta}^*)_S\|_2$. For the quadratic loss,

$$\nabla\mathcal{L}(\boldsymbol{\beta}^*)_S = n^{-1}\mathbf{X}_{*S}^{\text{T}}(\mathbf{y} - \mathbf{X}\boldsymbol{\beta}^*) = n^{-1}\mathbf{X}_{*S}^{\text{T}}\boldsymbol{\epsilon}.$$

Taking $\mathbf{v} = \boldsymbol{\epsilon}$, $\mathbf{A} = n^{-1}\mathbf{X}_{*S}\mathbf{X}_{*S}^{\text{T}}$ and $t = \mathbb{E}\boldsymbol{\epsilon}^{\text{T}}\mathbf{A}\boldsymbol{\epsilon}$ in the Hanson-Wright inequality (Lemma F.3) yields that

$$\mathbb{P}\big(|\boldsymbol{\epsilon}^{\text{T}}\mathbf{A}\boldsymbol{\epsilon} - \mathbb{E}\boldsymbol{\epsilon}^{\text{T}}\mathbf{A}\boldsymbol{\epsilon}| > \mathbb{E}\boldsymbol{\epsilon}^{\text{T}}\mathbf{A}\boldsymbol{\epsilon}\big) \leq 2\exp\bigg[-C_h\min\bigg\{\frac{\mathbb{E}\boldsymbol{\epsilon}^{\text{T}}\mathbf{A}\boldsymbol{\epsilon}}{\sigma^2\|\mathbf{A}\|_2}, \frac{(\mathbb{E}\boldsymbol{\epsilon}^{\text{T}}\mathbf{A}\boldsymbol{\epsilon})^2}{\sigma^4\|\mathbf{A}\|_F^2}\bigg\}\bigg]$$

$$(\text{B.3}) \qquad\qquad \leq 2\exp\bigg[-C_h\min\bigg\{\frac{s\sigma^2}{\sigma^2\lambda_{\max}(\mathbf{A})}, \frac{s^2\sigma^4}{s\sigma^4\lambda_{\max}^2(\mathbf{A})}\bigg\}\bigg],$$

where $C_h$ is a universal constant that does not depend on $n, d, s$; and $\mathbb{E}\boldsymbol{\epsilon}^{\text{T}}\mathbf{A}\boldsymbol{\epsilon} = s\sigma^2$ using the expectation of a quadratic form. Note that the non-zero singular values of $\mathbf{X}_{*S}^{\text{T}}\mathbf{X}_{*S}$ and $\mathbf{X}_{*S}\mathbf{X}_{*S}^{\text{T}}$ are the same and $\rho_-(s, r)$ is bounded above by $\rho^*$, we have

$$\|\mathbf{A}\|_2 = \bigg\|\frac{1}{n}\mathbf{X}_{*S}\mathbf{X}_{*S}^{\text{T}}\bigg\|_2 = \bigg\|\frac{1}{n}\mathbf{X}_{*S}^{\text{T}}\mathbf{X}_{*S}\bigg\|_2 = \rho_+(s, r) \leq \rho^*,$$

which, together with (B.3), results that

$$\mathbb{P}\Big(|\boldsymbol{\epsilon}^{\mathrm{T}}\mathbf{A}\boldsymbol{\epsilon}-\mathbb{E}\boldsymbol{\epsilon}^{\mathrm{T}}\mathbf{A}\boldsymbol{\epsilon}|>\mathbb{E}\boldsymbol{\epsilon}^{\mathrm{T}}\mathbf{A}\boldsymbol{\epsilon}\Big)\leq 2\exp\{-C_h's\},$$

where $C_h' = C_h \min\{1/\rho^*, 1/(\rho^*)^2\}$.

Define $\mathbf{P}_S$ to be the projection matrix into the column space of $\mathbf{X}_{*S}$. Since $\mathbf{P}_S \mathbf{X}_{*S} = \mathbf{X}_{*S}$, we obtain that $\boldsymbol{\epsilon}^{\mathrm{T}}\mathbf{A}\boldsymbol{\epsilon} = (\mathbf{P}_S\boldsymbol{\epsilon})^{\mathrm{T}}\mathbf{A}(\mathbf{P}_S\boldsymbol{\epsilon})$. Thus we have $\mathbb{E}\boldsymbol{\epsilon}^{\mathrm{T}}\mathbf{A}\boldsymbol{\epsilon} \leq \lambda_{\max}(\mathbf{A})\|\mathbf{P}_S\boldsymbol{\epsilon}\|_2^2$. Further define the event set $\mathcal{J}_2 = \big\{|\boldsymbol{\epsilon}^{\mathrm{T}}\mathbf{A}\boldsymbol{\epsilon} - \mathbb{E}\boldsymbol{\epsilon}^{\mathrm{T}}\mathbf{A}\boldsymbol{\epsilon}| \leq \mathbb{E}\boldsymbol{\epsilon}^{\mathrm{T}}\mathbf{A}\boldsymbol{\epsilon}\big\}$. Then with probability at least $\mathbb{P}(\mathcal{J}_2) \geq 1 - 2\exp\{-C_h's\}$,

$$\|\nabla\mathcal{L}(\boldsymbol{\beta}^*)_S\|_2 = \sqrt{\frac{1}{n}\boldsymbol{\epsilon}^{\mathrm{T}}\mathbf{A}\boldsymbol{\epsilon}} \leq \sqrt{\frac{2}{n}\mathbb{E}\boldsymbol{\epsilon}^{\mathrm{T}}\mathbf{A}\boldsymbol{\epsilon}} \leq \sqrt{\frac{2}{n}\rho^*\mathbb{E}\big[\|\mathbf{P}_S\boldsymbol{\epsilon}\|_2^2\big]}$$
$$= \sqrt{2\rho^*}\sigma\sqrt{s/n} \leq \sqrt{2\rho^*}\sigma\sqrt{s/n}.$$

Define $\mathcal{J} = \mathcal{J}_1 \cap \mathcal{J}_2$. Then, in the event $\mathcal{J}$, we have

$$\|\widetilde{\boldsymbol{\beta}}^{(\ell)} - \boldsymbol{\beta}^*\|_2 \leq C\big(\sqrt{2\rho^*}\sigma + 1\big)\sqrt{s/n} \propto \sqrt{s/n},$$

where $\mathbb{P}(\mathcal{J}) \geq 1 - \mathbb{P}(\mathcal{J}_1) - \mathbb{P}(\mathcal{J}_2)$. In other words, the above bound holds with probability at least $1 - 2d^{-\eta_1} - 2\exp\{-\eta_2 s\}$, in which $\eta_1 = c^2/(128\sigma^2)$ and $\eta_2 = C_h'$. □

We then give the proof for the oracle property under the LSE condition. Similar result holds under the LRE condition.

PROOF OF THEOREM 4.4. Let us define $S^{(\ell)} = \big\{j : |\widetilde{\beta}^{(\ell)} - \beta_j^*| > u\big\}$, where $u$ is defined in Assumption 4.2. We have $S^{(0)} = \{(i,j) : |\beta_j^*| \geq u\} = S$. We need several lemmas. Our first lemma bounds the discrepancy between $\widetilde{\boldsymbol{\beta}}^{(\ell)}$ and $\widehat{\boldsymbol{\beta}}^{\circ}$. The proof is similar to that of Lemma B.7.

**Lemma B.1.** Suppose Assumption 4.1 and 4.2 hold. Let $C = 12/\rho_*$. If $4(\|\nabla\mathcal{L}(\widehat{\boldsymbol{\beta}}^{\circ})\|_{\infty} + \varepsilon_c \vee \varepsilon_t) \leq \lambda \lesssim r/\sqrt{s}$, we must have $|\mathcal{E}_\ell| \leq 2s$, and for $\ell \geq 2$, the $\varepsilon$-optimal solution $\widetilde{\boldsymbol{\beta}}^{(\ell)}$ must satisfy

$$\|\widetilde{\boldsymbol{\beta}}^{(\ell)} - \widehat{\boldsymbol{\beta}}^{\circ}\|_2 \leq C\Big(\|\boldsymbol{\lambda}_{\mathcal{E}_\ell}^{(\ell-1)}\|_2 + \varepsilon_t\sqrt{|\mathcal{E}_\ell|}\Big).$$

Our second lemma connects $\lambda^{(\ell-1)}$ to $\widetilde{\boldsymbol{\beta}}^{(\ell-1)}$. The proof follows a similar argument used in the proof of Lemma 5.2 and thus is omitted.

**Lemma B.2.** We have

$$\big\|\boldsymbol{\lambda}_{\mathcal{E}_\ell}^{(\ell-1)}\big\|_2 \leq \lambda\big\|\mathrm{w}(|\boldsymbol{\beta}_S^*|-u)\big\|_2 + \lambda\big|\big\{j \in S : |\widetilde{\beta}_j^{(\ell-1)} - \beta_j^*| \geq u\big\}\big|^{1/2} + \lambda\sqrt{|\mathcal{E}_\ell \setminus S|}.$$

Combining the two above lemmas together and using the definition of $S^{(\ell)}$, we obtain

$$\|\widetilde{\boldsymbol{\beta}}^{(\ell)}-\widehat{\boldsymbol{\beta}}^\circ\|_2 \leq C\Big\{\underbrace{\lambda\|\mathrm{w}(|\boldsymbol{\beta}^*_S|-u)\|_2}_{\text{I}}+\lambda\sqrt{|S^{(\ell-1)}\cap S|}+\underbrace{\lambda\sqrt{|\mathcal{E}_\ell\setminus S|}}_{\text{II}}+\varepsilon\sqrt{|\mathcal{E}_\ell|}\Big\}.$$

Since $\|\boldsymbol{\beta}_S\|_{\min} \geq u + \alpha\lambda$ and $\mathrm{w}(\alpha\lambda) = 0$, we have $\text{I} = \lambda\|\mathrm{w}(|\boldsymbol{\beta}^*_S| - u)\|_2 = 0$. For any $j \in \mathcal{E}_\ell \setminus S$, we must have $\lambda_j^{(\ell-1)} < \lambda\mathrm{w}(u)$, and thus $|\widetilde{\beta}_j^{(\ell-1)}| = |\widetilde{\beta}_j^{(\ell-1)} - \boldsymbol{\beta}^*_{ij}| \geq u$ since $\boldsymbol{\beta}^*_j = 0$ for $j \in S^c$. This implies $\mathcal{E}_\ell \setminus S \in S^{(\ell-1)} \setminus S$, or equivalently $\text{II} \leq \lambda\sqrt{|S^{(\ell-1)} \setminus S|}$. Therefore, for $\ell \geq 2$, we have

$$\begin{aligned}(B.4)\quad \|\widetilde{\boldsymbol{\beta}}^{(\ell)}-\widehat{\boldsymbol{\beta}}^\circ\|_2 &\leq C\Big\{\lambda\sqrt{|S^{(\ell-1)}\cap S|}+\lambda\sqrt{|S^{(\ell-1)}\setminus S|}+\varepsilon_t\sqrt{|\mathcal{E}_\ell|}\Big\}\\ &\leq C\Big\{\lambda\sqrt{2|S^{(\ell-1)}|}+\varepsilon_t\sqrt{|\mathcal{E}_\ell|}\Big\}\end{aligned}$$

On the other hand, since $\|\widehat{\boldsymbol{\beta}}^\circ - \boldsymbol{\beta}^*\|_{\max} \leq \eta_n \leq \delta^{-1}\rho_*^{-1}\lambda$, $j \in S^{(\ell)}$ implies that

$$|\widetilde{\beta}_j^{(\ell)}-\widehat{\beta}_j^\circ| \geq |\widetilde{\beta}_j^{(\ell)}-\beta_j^*|-|\widehat{\beta}_j^\circ-\beta_j^*| \geq u-\delta^{-1}\rho_*^{-1}\lambda \geq 12\sqrt{2}\delta^{-1}\rho_*^{-1}\lambda.$$

We then bound $\sqrt{|S^{(\ell)}|}$ in terms of $\|\widetilde{\boldsymbol{\beta}}^{(\ell)} - \widehat{\boldsymbol{\beta}}^\circ\|_2$:

$$\sqrt{|S^{(\ell)}|} \leq \frac{\|\widetilde{\boldsymbol{\beta}}^{(\ell)}-\widehat{\boldsymbol{\beta}}^\circ\|_2}{u-\eta_n} \leq \delta\sqrt{|S^{(\ell-1)}|}+\delta\frac{\varepsilon_t\sqrt{s}}{\lambda}.$$

Doing induction on $\ell$ and using the fact that $S^{(0)} = S$, we obtain

$$\sqrt{|S^{(\ell)}|} \leq \delta^\ell\sqrt{s}+\delta^\ell\frac{\varepsilon_c\sqrt{s}}{\lambda}+\frac{\delta}{1-\delta}\frac{\varepsilon_t\sqrt{s}}{\lambda}$$

Thus, for $\ell$ large enough such that $\ell \gtrsim \log\{(1+\varepsilon_c/\lambda)\sqrt{s}\}$ and $\varepsilon_t$ small enough such that $\varepsilon_t \lesssim \lambda/\sqrt{s}$, we must have the right hand side of the above inequality is small than 1, which implies that

$$S^{(\ell)} = \varnothing \text{ and thus } \widetilde{\boldsymbol{\beta}}^{(\ell)} = \widehat{\boldsymbol{\beta}}^\circ.$$

Therefore, the estimator enjoys the strong oracle property.

□

**B.2. Key Lemmas.** In this section, we collect proofs for Lemma 5.1 and Lemma 5.2. We start with a proposition that connects the LSE condition to the localized versions of the sparse strong convexity/sparse strong smoothness (SSC/SSM) in [1], which will be frequently used in our theoretical analysis. Let

$$D_{\mathcal{L}}(\boldsymbol{\beta}_1, \boldsymbol{\beta}_2) \equiv \mathcal{L}(\boldsymbol{\beta}_1) - \mathcal{L}(\boldsymbol{\beta}_2) - \langle \nabla \mathcal{L}(\boldsymbol{\beta}_2), \boldsymbol{\beta}_1 - \boldsymbol{\beta}_2 \rangle,$$

and $D_{\mathcal{L}}^s(\boldsymbol{\beta}_1, \boldsymbol{\beta}_2) \equiv D_{\mathcal{L}}(\boldsymbol{\beta}_1, \boldsymbol{\beta}_2) + D_{\mathcal{L}}(\boldsymbol{\beta}_2, \boldsymbol{\beta}_1)$.

**Proposition B.3.** *For any $\boldsymbol{\beta}_1, \boldsymbol{\beta}_2 \in B_2(r, \boldsymbol{\beta}^*) \equiv \{\boldsymbol{\beta} : \|\boldsymbol{\beta} - \boldsymbol{\beta}^*\|_2 \leq r\}$ such that $\|\boldsymbol{\beta}_1 - \boldsymbol{\beta}_2\|_0 \leq m$, we have*

$$\frac{1}{2}\rho_-(m,r)\|\boldsymbol{\beta}_1 - \boldsymbol{\beta}_2\|_2^2 \leq D_{\mathcal{L}}(\boldsymbol{\beta}_1, \boldsymbol{\beta}_2) \leq \frac{1}{2}\rho_+(m,r)\|\boldsymbol{\beta}_1 - \boldsymbol{\beta}_2\|_2^2,$$

$$\rho_-(m,r)\|\boldsymbol{\beta}_1 - \boldsymbol{\beta}_2\|_2^2 \leq D_{\mathcal{L}}^s(\boldsymbol{\beta}_1, \boldsymbol{\beta}_2) \leq \rho_+(m,r)\|\boldsymbol{\beta}_1 - \boldsymbol{\beta}_2\|_2^2.$$

PROOF OF PROPOSITION B.3. We prove the second inequality. By the mean value theorem, there exists a $\gamma \in [0,1]$ such that $\widetilde{\boldsymbol{\beta}} = \gamma \boldsymbol{\beta}_1 + (1-\gamma)\boldsymbol{\beta}_2 \in B_2(r, \boldsymbol{\beta}^*)$, $\|\widetilde{\boldsymbol{\beta}}\|_0 \leq m$ and

$$\langle \nabla \mathcal{L}(\boldsymbol{\beta}_1) - \nabla \mathcal{L}(\boldsymbol{\beta}_2), \boldsymbol{\beta}_1 - \boldsymbol{\beta}_2 \rangle = (\boldsymbol{\beta}_1 - \boldsymbol{\beta}_2)^T \{\nabla^2 \mathcal{L}(\widetilde{\boldsymbol{\beta}})\}(\boldsymbol{\beta}_1 - \boldsymbol{\beta}_2).$$

By the definition of the localized sparse eigenvalue, we obtain the desired result. The other inequality can be proved similarly. □

We then present the proof for Lemma 5.1 below.

PROOF OF LEMMA 5.1. If we assume that, for all $\ell \geq 1$, the following two inequalities hold

(B.5) $\quad |\mathcal{E}_\ell| = k \leq 2s$, where $\mathcal{E}_\ell$ is defined in (5.1), and

(B.6) $\quad \|\boldsymbol{\lambda}_{\mathcal{E}_\ell^c}^{(\ell-1)}\|_{\min} \geq \lambda/2 \geq \|\nabla \mathcal{L}(\boldsymbol{\beta}^*)\|_\infty + \varepsilon$.

Applying Lemma E.2 in the online supplement, we obtain the desired bound:

$$\|\widetilde{\boldsymbol{\beta}}^{(\ell)} - \boldsymbol{\beta}^*\|_2 \leq \frac{12}{\rho_*}(\|\boldsymbol{\lambda}_S^{(\ell-1)}\|_2 + \|\nabla \mathcal{L}(\boldsymbol{\beta}^*)_{\mathcal{E}_\ell}\|_2 + \varepsilon\sqrt{|\mathcal{E}_\ell|}) < \frac{18}{\rho_*}\lambda\sqrt{s} \leq r.$$

Therefore, it remains to show (B.5) and (B.6) hold for all $\ell \geq 1$. We prove these by induction. For $\ell = 1$, $\lambda \geq \lambda w(u)$ and thus $\mathcal{E}_1 = S$, which implies (B.5) and (B.6). Assume these two statements hold at $\ell - 1$. Since $j \in \mathcal{E}_\ell \backslash S$ implies $j \notin S$ and $\lambda w(|\widetilde{\beta}_j^{(\ell-1)}|) = \lambda_j^{(\ell)} < \lambda w(u)$ by definition, and since

w($x$) is non-increasing, we obtain that $|\widetilde{\beta}_j^{\ell-1}| \geq u$. Therefore by induction hypothesis, we have

$$\sqrt{|\mathcal{E}_\ell \setminus S|} \leq \frac{\|\widetilde{\boldsymbol{\beta}}_{\mathcal{E}_\ell \setminus S}^{(\ell-1)}\|}{u} \leq \frac{\|\widetilde{\boldsymbol{\beta}}^{(\ell-1)} - \boldsymbol{\beta}^*\|_2}{u} \leq \frac{18\lambda}{\rho_* u}\sqrt{s} \leq \sqrt{s},$$

where the last inequality uses the definition of $u$ in (5.1). The inequality above implies that $|\mathcal{E}_\ell| \leq 2s$. For such $\mathcal{E}_\ell$, we have $\|\boldsymbol{\lambda}_{\mathcal{E}_\ell^c}\|_{\min} \geq \lambda w(u) \geq \lambda/2 \geq \|\nabla \mathcal{L}(\boldsymbol{\beta}^*)\|_\infty + \varepsilon$, which completes the induction step. This completes the proof. □

PROOF OF LEMMA 5.2. If $|\beta_j^* - \widetilde{\beta}_j| \geq u$, then $w(|\widetilde{\beta}_j|) \leq 1 \leq u^{-1}|\widetilde{\beta}_j - \beta_j^*|$; otherwise, $w(|\widetilde{\beta}_j|) \leq w(|\beta_j^*| - u)$. Therefore, the following inequality always hold

$$w(|\widetilde{\beta}_j|) \leq w(|\beta_j^*| - u) + u^{-1}|\beta_j^* - \widetilde{\beta}|.$$

Applying the triangle inequality completes the proof. □

**B.3. Statistical Theory under LRE Condition.** In this section, we present the main theorem and its proof, with some technical lemmas postponed to later sections. We formally introduce the LRE condition below.

**Assumption B.1.** There exist $k \leq 2s$, $\gamma$ and $r \gtrsim \lambda\sqrt{s}$ such that $0 < \kappa_* \leq \kappa_-(k, \gamma, r) \leq \kappa_+(k, \gamma, r) \leq \kappa^* < \infty$.

B.3.1. *Proofs of Main Theorems.* We begin with a proposition, which establishes the relationship between localized restricted eigenvalue and the localized version of the restricted strong convexity/smoothness. The proof is similar to that of Proposition B.3, and thus is omitted.

**Proposition B.4.** For any $\boldsymbol{\beta}_1, \boldsymbol{\beta}_2 \in \mathcal{C}(k, \gamma, r) \cap B_2(r/2, \boldsymbol{\beta}^*)$, we have

$$\frac{1}{2}\kappa_-(k, \gamma, r)\|\boldsymbol{\beta}_1 - \boldsymbol{\beta}_2\|_2^2 \leq D_\mathcal{L}(\boldsymbol{\beta}_1, \boldsymbol{\beta}_2) \leq \frac{1}{2}\kappa_+(k, \gamma, r)\|\boldsymbol{\beta}_1 - \boldsymbol{\beta}_2\|_2^2, \text{ and}$$
$$\kappa_-(k, \gamma, r)\|\boldsymbol{\beta}_1 - \boldsymbol{\beta}_2\|_2^2 \leq D_\mathcal{L}^s(\boldsymbol{\beta}_1, \boldsymbol{\beta}_2) \leq \kappa_+(k, \gamma, r)\|\boldsymbol{\beta}_1 - \boldsymbol{\beta}_2\|_2^2.$$

Next, we bound the $\ell_2$ error using the regularization parameter. The proof is similar to that of Lemma 5.1 and depends on Lemma B.7, where we introduce the localized analysis such that the localized restricted eigenvalue condition can be applied.

**Lemma B.5.** Suppose that Assumption 4.2 holds with $u = 18\lambda/\kappa^*$, and Assumption B.1 holds with a $k \leq 2s$, $\gamma = 5$. If $\lambda$, $\varepsilon$ and $r$ satisfy

$$4(\|\nabla \mathcal{L}(\boldsymbol{\beta}^*)\|_\infty + \varepsilon) \leq \lambda \leq r\kappa_*/(2\sqrt{s}),$$

then $|\mathcal{E}_\ell| \leq 2s$ and any $\varepsilon$-optimal solution $\widetilde{\boldsymbol{\beta}}^{(\ell)}$ ($\ell \geq 1$) must satisfy

$$\|\widetilde{\boldsymbol{\beta}}^{(\ell)} - \boldsymbol{\beta}^*\|_2 \leq \kappa_*^{-1}\big(\|\boldsymbol{\lambda}_S^{(\ell-1)}\|_2 + \|\nabla \mathcal{L}(\boldsymbol{\beta}^*)_{\mathcal{E}_\ell}\|_2 + \varepsilon\sqrt{|\mathcal{E}_\ell|}\big) \lesssim \lambda\sqrt{s}.$$

Recall Lemma 5.2, which bounds the regularization parameter using the functional of the estimator from previous step. Combining these two lemmas together, we obtain the following main theorem.

**Theorem B.6** (Optimal Statistical Rate under Localized Restricted Eigenvalue Condition)**.** Suppose the same conditions of Lemma B.5 hold, but with $\varepsilon$ replaced by $\varepsilon_c \vee \varepsilon_t$. Then, for $\ell \geq 2$ and some constant $C$, any $\varepsilon_t$-optimal solution must satisfy

$$\|\widetilde{\boldsymbol{\beta}}^{(\ell)} - \boldsymbol{\beta}^*\|_2 \lesssim \underbrace{\|\nabla \mathcal{L}(\boldsymbol{\beta}^*)_S\|_2}_{\text{oracel rate}} + \overbrace{\varepsilon_t\sqrt{s}}^{\text{opt err}} + \underbrace{\lambda\|\mathrm{w}_S(|\boldsymbol{\beta}_S^*| - u)\|_2}_{\text{coefficient effect}} + \overbrace{\delta^{\ell-1}\lambda\sqrt{s}}^{\text{tightening effect}}.$$

PROOF OF THEOREM B.6. Under the conditions of the theorem, Lemma B.5 directly implies $|\mathcal{E}_\ell| \leq 2s$ and

$$\|\boldsymbol{\lambda}_{\mathcal{E}_\ell^c}^{(\ell-1)}\|_{\min} \geq \|\nabla\mathcal{L}(\boldsymbol{\beta}^*)\|_\infty + \varepsilon, \text{ for all } \ell \geq 1,$$

where $\varepsilon = \varepsilon_c \vee \varepsilon_t$. Using Lemma B.7 then obtains us that

$$(\text{B.7}) \qquad \|\widetilde{\boldsymbol{\beta}}^{(\ell)} - \boldsymbol{\beta}^*\|_2 \leq \kappa_*^{-1}\big(\|\boldsymbol{\lambda}_S^{(\ell-1)}\|_2 + \|\nabla \mathcal{L}(\boldsymbol{\beta}^*)_{\mathcal{E}_\ell}\|_2 + \varepsilon\sqrt{|\mathcal{E}_\ell|}\big).$$

To bound the first term in the inequality above, we apply Lemma 5.2 and obtain

$$(\text{B.8}) \qquad \|\widetilde{\boldsymbol{\beta}}^{(\ell)} - \boldsymbol{\beta}^*\|_2 \leq \frac{1}{\kappa_*}\Big(\underbrace{\|\nabla\mathcal{L}(\boldsymbol{\beta}^*)_{\mathcal{E}_\ell}\|_2 + \varepsilon_t\sqrt{|\mathcal{E}_\ell|}}_{\text{I}} + \lambda\|\mathrm{w}_S(|\boldsymbol{\beta}_S^*| - u)\|_2\Big)$$
$$+ \frac{\lambda}{u\kappa_*}\|\widetilde{\boldsymbol{\beta}}^{(\ell-1)} - \boldsymbol{\beta}^*\|_2.$$

Following a similar argument in the proof of Theorem 4.2, the term I can be bounded by

$$\|\nabla\mathcal{L}(\boldsymbol{\beta}^*)_S\|_2 + \varepsilon_t\sqrt{s} + \frac{\lambda}{4u}\|\widetilde{\boldsymbol{\beta}}^{(\ell-1)} - \boldsymbol{\beta}^*\|_2,$$

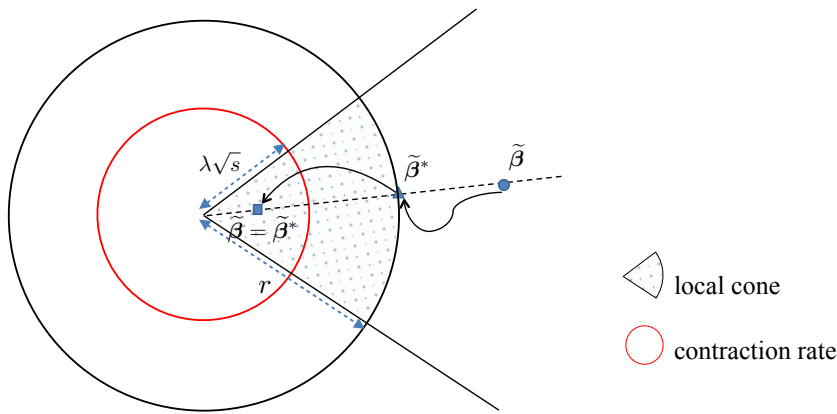

FIG 5. *The localized analysis under localized restricted eigenvalue condition: an intermediate solution, $\widetilde{\boldsymbol{\beta}}^\star$, is constructed such that $\|\widetilde{\boldsymbol{\beta}}^\star - \boldsymbol{\beta}^*\|_2 = r$ if $\|\widetilde{\boldsymbol{\beta}} - \boldsymbol{\beta}^*\|_2 > r$, $\widetilde{\boldsymbol{\beta}}^\star = \widetilde{\boldsymbol{\beta}}$, else. Employing the localized eigenvalue condition, this intermediate solution converges at the rate of $\lambda\sqrt{s}$ in the contraction stage and $\sqrt{s/n}$ in the tightening stage. Under mild conditions, the approximate solution $\widetilde{\boldsymbol{\beta}} = \widetilde{\boldsymbol{\beta}}^\star$ and thus has same convergence rate as the intermediate solution. The shaded area is a local $\ell_1$-cone where the restricted eigenvalue condition holds.*

which, combining with (B.8), yields that

$$\|\widetilde{\boldsymbol{\beta}}^{(\ell)} - \boldsymbol{\beta}^*\|_2 \lesssim \|\nabla \mathcal{L}(\boldsymbol{\beta}^*)_S\|_2 + \varepsilon_t \sqrt{s} + \lambda \|\mathrm{w}_S(|\boldsymbol{\beta}_S^*| - u)\|_2 \\ + \delta \|\widetilde{\boldsymbol{\beta}}^{(\ell-1)} - \boldsymbol{\beta}^*\|_2,$$

where $\delta = 2\lambda/(u\kappa_*) \leq 1/2$. The proof is completed by applying the above inequality recursively. $\square$

B.3.2. *Localized Analysis.* In this section, we carry out the localized analysis. A geometric explanation of the localized analysis is given in Figure 5. The following lemma shows that the approximate solution always falls in the neighborhood of $\boldsymbol{\beta}^*$ by a novel localized analysis and thus the localized restricted eigenvalue condition can be exploited. Recall the definitions of $D_\mathcal{L}(\cdot, \cdot)$, $D_\mathcal{L}^s(\cdot, \cdot)$ and $\mathcal{C}(k, \gamma, r)$ in Section B.2.

**Lemma B.7.** Suppose Assumption B.1 holds. Take $\mathcal{E}$ such that $S \cap \mathcal{E}^c = \emptyset$, $|\mathcal{E}| \leq k \leq 2s$. Further assume that $\|\boldsymbol{\lambda}_{\mathcal{E}^c}\|_{\min} \geq \|\nabla \mathcal{L}(\boldsymbol{\beta}^*)\|_\infty + \varepsilon$ and $2\kappa_*^{-1} \lambda \sqrt{s} \leq r$. Then any $\varepsilon$-optimal solution $\widetilde{\boldsymbol{\beta}}$ must satisfy

$$\|\widetilde{\boldsymbol{\beta}} - \boldsymbol{\beta}^*\|_2 \leq \kappa_*^{-1} \big(\|\boldsymbol{\lambda}_S\|_2 + \|\nabla \mathcal{L}(\boldsymbol{\beta}^*)_\mathcal{E}\|_2 + \varepsilon \sqrt{|\mathcal{E}|}\big) \lesssim \lambda \sqrt{s}.$$

PROOF OF LEMMA B.7. Let $\widetilde{\boldsymbol{\beta}}^{\star} = \boldsymbol{\beta}^* + t(\widetilde{\boldsymbol{\beta}} - \boldsymbol{\beta}^*)$, where $t = 1$, if $\|\widetilde{\boldsymbol{\beta}} - \boldsymbol{\beta}^*\|_2 \leq r$; $t \in (0, 1)$ such that $\|\widetilde{\boldsymbol{\beta}}^{\star} - \boldsymbol{\beta}^*\|_2 = r$, otherwise. By the definition of $\widetilde{\boldsymbol{\beta}}^*$, we know that $\|\widetilde{\boldsymbol{\beta}}^{\star} - \boldsymbol{\beta}^*\|_2 \leq r$. Using the $\ell_1$-cone lemma, we know that the approximate solution falls in the $\ell_1$-cone, i.e.

$$\text{(B.9)} \qquad \|(\widetilde{\boldsymbol{\beta}} - \boldsymbol{\beta}^*)_{\mathcal{E}}\|_1 \leq 5\|(\widetilde{\boldsymbol{\beta}} - \boldsymbol{\beta}^*)_{\mathcal{E}^c}\|_1,$$

From the construction of $\widetilde{\boldsymbol{\beta}}^*$, we know $\widetilde{\boldsymbol{\beta}}^* - \boldsymbol{\beta}^* = t(\widetilde{\boldsymbol{\beta}} - \boldsymbol{\beta}^*)$. Thus, we have

$$\|(\widetilde{\boldsymbol{\beta}}^{\star} - \boldsymbol{\beta}^*)_{\mathcal{E}}\|_1 \leq 5\|(\widetilde{\boldsymbol{\beta}}^{\star} - \boldsymbol{\beta}^*)_{\mathcal{E}^c}\|_1.$$

Combining the inequality above with the assumption $|\mathcal{E}| \leq k$ results that $\widetilde{\boldsymbol{\beta}}^*$ falls in the local $\ell_1$-cone, i.e. $\widetilde{\boldsymbol{\beta}}^* \in \mathcal{C}(k, c_0, r)$. Then Proposition B.4 implies the localized restricted strong convexity, i.e.

$$\text{(B.10)} \qquad \kappa_-(k, 5, r)\|\widetilde{\boldsymbol{\beta}}^{\star} - \boldsymbol{\beta}^*\|_2^2 \leq D_{\mathcal{L}}^s(\widetilde{\boldsymbol{\beta}}^*, \boldsymbol{\beta}^*).$$

To bound the right hand side of the above inequality, we use Lemma F.2 in Appendix F:

$$\text{(B.11)} \qquad D_{\mathcal{L}}^s(\widetilde{\boldsymbol{\beta}}^*, \boldsymbol{\beta}^\circ) \leq t D_{\mathcal{L}}^s(\widetilde{\boldsymbol{\beta}}, \boldsymbol{\beta}^*) = t\langle \nabla \mathcal{L}(\widetilde{\boldsymbol{\beta}}) - \nabla \mathcal{L}(\boldsymbol{\beta}^*), \widetilde{\boldsymbol{\beta}} - \boldsymbol{\beta}^*\rangle.$$

It suffices to bound the right hand side of (B.11). Plugging (B.11) back into (B.10) and adding $\langle \boldsymbol{\lambda} \odot \boldsymbol{\xi}, \widetilde{\boldsymbol{\beta}} - \boldsymbol{\beta}^*\rangle$ to both sides, we obtain

$$\text{(B.12)} \qquad \kappa_-(k, 5, r)\|\widetilde{\boldsymbol{\beta}}^* - \boldsymbol{\beta}^*\|^2 + t\underbrace{\langle \nabla \mathcal{L}(\boldsymbol{\beta}^*), \widetilde{\boldsymbol{\beta}} - \boldsymbol{\beta}^*\rangle}_{\text{I}} + t\underbrace{\langle \boldsymbol{\lambda} \odot \boldsymbol{\xi}, \widetilde{\boldsymbol{\beta}} - \boldsymbol{\beta}^*\rangle}_{\text{II}}$$
$$\leq t\underbrace{\langle \nabla \mathcal{L}(\widetilde{\boldsymbol{\beta}}) + \boldsymbol{\lambda} \odot \boldsymbol{\xi}, \widetilde{\boldsymbol{\beta}} - \boldsymbol{\beta}^*\rangle}_{\text{III}}.$$

It remains to bound terms I, II and III respectively. For I, separating the support of $\nabla \mathcal{L}(\boldsymbol{\beta}^*)$ and $\widetilde{\boldsymbol{\beta}} - \boldsymbol{\beta}^*$ to $\mathcal{E}$ and $\mathcal{E}^c$ and using the Hölder inequality, we obtain

$$\text{I} = \langle (\nabla \mathcal{L}(\boldsymbol{\beta}^*))_{\mathcal{E}}, (\widetilde{\boldsymbol{\beta}} - \boldsymbol{\beta}^*)_{\mathcal{E}}\rangle + \langle (\nabla \mathcal{L}(\boldsymbol{\beta}^*))_{\mathcal{E}^c}, (\widetilde{\boldsymbol{\beta}} - \boldsymbol{\beta}^*)_{\mathcal{E}^c}\rangle$$
$$\text{(B.13)} \quad \geq -\|(\nabla \mathcal{L}(\boldsymbol{\beta}^*))_{\mathcal{E}}\|_2\|(\widetilde{\boldsymbol{\beta}} - \boldsymbol{\beta}^*)_{\mathcal{E}}\|_2 - \|(\nabla \mathcal{L}(\boldsymbol{\beta}^*))_{\mathcal{E}^c}\|_\infty\|(\widetilde{\boldsymbol{\beta}} - \boldsymbol{\beta}^*)_{\mathcal{E}^c}\|_1.$$

For II, separating the support of $\boldsymbol{\lambda} \odot \boldsymbol{\xi}$ and $\widetilde{\boldsymbol{\beta}} - \boldsymbol{\beta}^*$ to $S, \mathcal{E} \setminus S$ and $\mathcal{E}^c$ results

$$\text{II} = \langle (\boldsymbol{\lambda} \odot \boldsymbol{\xi})_S, (\widetilde{\boldsymbol{\beta}} - \boldsymbol{\beta}^*)_S\rangle + \langle (\boldsymbol{\lambda} \odot \boldsymbol{\xi})_{\mathcal{E} \setminus S}, (\widetilde{\boldsymbol{\beta}} - \boldsymbol{\beta}^*)_{\mathcal{E} \setminus S}\rangle$$
$$+ \langle (\boldsymbol{\lambda} \odot \boldsymbol{\xi})_{\mathcal{E}^c}, (\widetilde{\boldsymbol{\beta}} - \boldsymbol{\beta}^*)_{\mathcal{E}^c}\rangle.$$

To bound the last term in the above equality, note that $\mathcal{E}^c \cap S = \emptyset$ and thus

$$\langle (\boldsymbol{\lambda} \odot \boldsymbol{\xi})_{\mathcal{E}^c}, (\widetilde{\boldsymbol{\beta}} - \boldsymbol{\beta}^*)_{\mathcal{E}^c} \rangle = \langle \boldsymbol{\lambda}_{\mathcal{E}^c}, |\widetilde{\boldsymbol{\beta}}_{\mathcal{E}^c}| \rangle = \langle \boldsymbol{\lambda}_{\mathcal{E}^c}, |(\widetilde{\boldsymbol{\beta}} - \boldsymbol{\beta}^*)_{\mathcal{E}^c}| \rangle,$$

which yields that

$$\begin{aligned}
\text{II} &= \langle (\boldsymbol{\lambda} \odot \boldsymbol{\xi})_S, (\widetilde{\boldsymbol{\beta}} - \boldsymbol{\beta}^*)_S \rangle + \langle (\boldsymbol{\lambda} \odot \boldsymbol{\xi})_{\mathcal{E} \setminus S}, (\widetilde{\boldsymbol{\beta}} - \boldsymbol{\beta}^*)_{\mathcal{E} \setminus S} \rangle \\
&\quad + \langle \boldsymbol{\lambda}_{\mathcal{E}^c}, |(\widetilde{\boldsymbol{\beta}} - \boldsymbol{\beta}^*)_{\mathcal{E}^c}| \rangle \\
&\geq -\|\boldsymbol{\lambda}_S\|_2 \|(\widetilde{\boldsymbol{\beta}} - \boldsymbol{\beta}^*)_S\|_2 + \langle \boldsymbol{\lambda}_{\mathcal{E}^c}, |(\widetilde{\boldsymbol{\beta}} - \boldsymbol{\beta}^*)_{\mathcal{E}^c}| \rangle \\
(\text{B.14}) \quad &\geq -\|\boldsymbol{\lambda}_S\|_2 \|(\widetilde{\boldsymbol{\beta}} - \boldsymbol{\beta}^*)_S\|_2 + \|\boldsymbol{\lambda}_{\mathcal{E}^c}\|_{\min} \|(\widetilde{\boldsymbol{\beta}} - \boldsymbol{\beta}^*)_{\mathcal{E}^c}\|_1,
\end{aligned}$$

where the first inequality is due to the fact that $\langle (\boldsymbol{\lambda} \odot \boldsymbol{\xi})_{\mathcal{E} \setminus S}, (\widetilde{\boldsymbol{\beta}} - \boldsymbol{\beta}^*)_{\mathcal{E} \setminus S} \rangle = 0$ and we use Hölder inequality in the last inequality. For III, we first write $\mathbf{u} = \nabla \mathcal{L}(\widetilde{\boldsymbol{\beta}}) + \boldsymbol{\lambda} \odot \boldsymbol{\xi}$. Using similar arguments, we obtain

$$\begin{aligned}
\langle \nabla \mathcal{L}(\widetilde{\boldsymbol{\beta}}) + \boldsymbol{\lambda} \odot \boldsymbol{\xi}, \widetilde{\boldsymbol{\beta}} - \boldsymbol{\beta}^* \rangle &= \langle \mathbf{u}_{\mathcal{E}}, (\widetilde{\boldsymbol{\beta}} - \boldsymbol{\beta}^*)_{\mathcal{E}} \rangle + \langle \mathbf{u}_{\mathcal{E}^c}, (\widetilde{\boldsymbol{\beta}} - \boldsymbol{\beta}^*)_{\mathcal{E}^c} \rangle \\
(\text{B.15}) \quad &\leq \|\mathbf{u}_{\mathcal{E}}\|_2 \|(\widetilde{\boldsymbol{\beta}} - \boldsymbol{\beta}^*)_{\mathcal{E}}\|_2 + \|\mathbf{u}_{\mathcal{E}^c}\|_{\infty} \|(\widetilde{\boldsymbol{\beta}} - \boldsymbol{\beta}^*)_{\mathcal{E}^c}\|_1.
\end{aligned}$$

Plugging (B.13), (B.14), (B.15) into (B.12) and taking inf over $\boldsymbol{\xi} \in \partial \|\widetilde{\boldsymbol{\beta}}\|_1$, we obtain

$$\begin{aligned}
\kappa_-(k, 5, r) &\|\widetilde{\boldsymbol{\beta}}^* - \boldsymbol{\beta}^*\|^2 + t(\|\boldsymbol{\lambda}_{\mathcal{E}^c}\|_{\min} - \|\nabla \mathcal{L}(\boldsymbol{\beta}^*)\|_{\infty}) \|(\widetilde{\boldsymbol{\beta}} - \boldsymbol{\beta}^*)_{\mathcal{E}^c}\|_1 \\
&- t(\|\nabla \mathcal{L}(\boldsymbol{\beta}^*)\|_2 + \|\boldsymbol{\lambda}_{\mathcal{E}}\|_2) \|(\widetilde{\boldsymbol{\beta}} - \boldsymbol{\beta}^*)_{\mathcal{E}}\|_2 \\
&\leq t \inf_{\boldsymbol{\xi} \in \partial \|\widetilde{\boldsymbol{\beta}}\|_1} \|\mathbf{u}_{\mathcal{E}}\|_2 \|(\widetilde{\boldsymbol{\beta}} - \boldsymbol{\beta}^*)_{\mathcal{E}}\|_2 + t \inf_{\boldsymbol{\xi} \in \partial \|\widetilde{\boldsymbol{\beta}}\|_1} \|\mathbf{u}_{\mathcal{E}^c}\|_{\infty} \|(\widetilde{\boldsymbol{\beta}} - \boldsymbol{\beta}^*)_{\mathcal{E}^c}\|_1 \\
&\leq \varepsilon \sqrt{|\mathcal{E}|} \times t \|(\widetilde{\boldsymbol{\beta}} - \boldsymbol{\beta}^*)_{\mathcal{E}}\|_2 + \varepsilon \times t \|(\widetilde{\boldsymbol{\beta}} - \boldsymbol{\beta}^*)_{\mathcal{E}^c}\|_1,
\end{aligned}$$

where we use the fact $\inf_{\boldsymbol{\xi} \in \partial \|\widetilde{\boldsymbol{\beta}}\|_1} \|\mathbf{u}_{\mathcal{E}}\|_2 \leq \inf_{\boldsymbol{\xi} \in \partial \|\widetilde{\boldsymbol{\beta}}\|_1} \sqrt{|\mathcal{E}|} \|\mathbf{u}_{\mathcal{E}}\|_{\infty} \leq \varepsilon \sqrt{|\mathcal{E}|}$ in the last inequality. After some algebra, we obtain

$$\begin{aligned}
\kappa_-(k, 5, r) &\|\widetilde{\boldsymbol{\beta}}^* - \boldsymbol{\beta}^*\|_2^2 + t(\|\boldsymbol{\lambda}_{\mathcal{E}^c}\|_{\min} - (\|\nabla \mathcal{L}(\boldsymbol{\beta}^*)\|_{\infty} + \varepsilon)) \|(\widetilde{\boldsymbol{\beta}} - \boldsymbol{\beta}^*)_{\mathcal{E}^c}\|_1 \\
&\leq (\|\boldsymbol{\lambda}_S\|_2 + \|\nabla \mathcal{L}(\boldsymbol{\beta}^*)_{\mathcal{E}}\|_2 + \varepsilon \sqrt{|\mathcal{E}|}) \times t \|(\widetilde{\boldsymbol{\beta}} - \boldsymbol{\beta}^*)_{\mathcal{E}}\|_2.
\end{aligned}$$

Using the assumption that $\|\boldsymbol{\lambda}_{\mathcal{E}^c}\|_{\min} \geq \|\nabla \mathcal{L}(\boldsymbol{\beta}^*)\|_{\infty} + \varepsilon$, the inequality above can be simplified to

$$(\text{B.16}) \quad \kappa_-(k, 5, r) \|\widetilde{\boldsymbol{\beta}}^* - \boldsymbol{\beta}^*\|_2^2 \leq \underbrace{(\|\boldsymbol{\lambda}_S\|_2 + \|\nabla \mathcal{L}(\boldsymbol{\beta}^*)_{\mathcal{E}}\|_2 + \varepsilon \sqrt{|\mathcal{E}|})}_{\text{(i)}}$$
$$\times \underbrace{t \|(\widetilde{\boldsymbol{\beta}} - \boldsymbol{\beta}^*)_{\mathcal{E}}\|_2}_{\text{(ii)}}.$$

For (i), using the fact that $\|\boldsymbol{\lambda}_S\|_2 \leq \|\boldsymbol{\lambda}_S\|_\infty \sqrt{|S|} \leq \lambda\sqrt{s}$, we have

$$(\text{i}) \leq \lambda\sqrt{|S|} + (\|\nabla\mathcal{L}(\boldsymbol{\beta}^*)_{\mathcal{E}}\|_\infty + \varepsilon)\sqrt{|\mathcal{E}|} \leq \lambda\sqrt{s} + \frac{1}{4}\lambda\sqrt{k},$$

where, in the first inequality, we have used

$$\|\nabla\mathcal{L}(\boldsymbol{\beta}^*)_{\mathcal{E}}\|_2 + \varepsilon\sqrt{|\mathcal{E}|} \leq (\|\nabla\mathcal{L}(\boldsymbol{\beta}^*)_{\mathcal{E}}\|_\infty + \varepsilon)\sqrt{|\mathcal{E}|} \leq \frac{1}{4}\lambda\sqrt{|\mathcal{E}|}.$$

For (ii), we have $t\|(\widetilde{\boldsymbol{\beta}} - \boldsymbol{\beta}^*)_{\mathcal{E}}\|_2 = \|(\widetilde{\boldsymbol{\beta}}^* - \boldsymbol{\beta}^*)_{\mathcal{E}}\|_2$. Plugging the bounds for (i) and (ii) into (B.16) and using the assumption $2\kappa_*^{-1}\lambda\sqrt{s} \leq r$, we obtain

$$\|\widetilde{\boldsymbol{\beta}}^* - \boldsymbol{\beta}^*\|_2 \leq \frac{1 + \sqrt{2}/2}{\kappa_*}\lambda\sqrt{s} < r,$$

which is a contraction with the construction of $\widetilde{\boldsymbol{\beta}}^*$. This indicates that $\widetilde{\boldsymbol{\beta}}^* = \widetilde{\boldsymbol{\beta}}$. Therefore, the desired bound hold for $\widetilde{\boldsymbol{\beta}}$. □

## APPENDIX C: COMPUTATIONAL THEORY

In this section, we collect proofs for Proposition 4.5, Proposition 4.6 and Theorem 4.7. We then give the proofs for Lemma 5.3, Lemma 5.4 and Lemma 5.5. The proofs of technical lemmas are postponed to later sections. We denote the quadratic coefficient $\phi$ by $\phi_c$ in the contraction stage, and by $\phi_t$ in the tightening stage.

**C.1. Proofs of Main Results.** We start with the contraction stage and give the proof of Proposition 4.5, followed by that of Proposition 4.6.

PROOF OF PROPOSITION 4.5. We omit the super script $\ell$ in $\boldsymbol{\beta}^{(\ell,k)}$, 1 in $\widehat{\boldsymbol{\beta}}^{(1)}$, and 0 in $\boldsymbol{\lambda}^{(0)}$ for simplicity. Applying Lemma E.5 results that

$$(\text{C.1}) \qquad \omega_{\boldsymbol{\lambda}}(\boldsymbol{\beta}^{(k+1)}) \leq (\phi_c + \rho_c)\|\boldsymbol{\beta}^{(k+1)} - \boldsymbol{\beta}^{(k)}\|_2.$$

On the other hand, taking $\boldsymbol{\beta} = \boldsymbol{\beta}^{(k)}$ in Lemma E.4, we obtain

$$F(\boldsymbol{\beta}^{(k)}, \boldsymbol{\lambda}^{(0)}) - F(\boldsymbol{\beta}^{(k+1)}, \boldsymbol{\lambda}^{(0)}) \geq \frac{\phi_c}{2}\|\boldsymbol{\beta}^{(k+1)} - \boldsymbol{\beta}^{(k)}\|_2^2.$$

Plugging the inequality back into (C.1), we obtain a bound for the suboptimality measure

$$(\text{C.2}) \quad \omega_{\boldsymbol{\lambda}}(\boldsymbol{\beta}^{(k+1)}) \leq (\phi_c + \rho_c)\left\{\frac{2}{\phi_c}\left[F(\boldsymbol{\beta}^{(k)}, \boldsymbol{\lambda}) - F(\boldsymbol{\beta}^{(k+1)}, \boldsymbol{\lambda})\right]\right\}^{1/2}.$$

Since $\{F(\boldsymbol{\beta}^{(k)}, \boldsymbol{\lambda}^{(0)})\}_{k=0}^{\infty}$ is a monotone decreasing sequence, we have

(C.3) $\qquad F(\widehat{\boldsymbol{\beta}}^{(1)}, \boldsymbol{\lambda}) \leq \ldots \leq F(\boldsymbol{\beta}^{(k)}, \boldsymbol{\lambda}) \leq \ldots \leq F(\boldsymbol{\beta}^{(0)}, \boldsymbol{\lambda}).$

Plugging (C.3) back into (C.2) and using Lemma 5.3, we obtain

(C.4) $\quad \omega_{\boldsymbol{\lambda}}(\boldsymbol{\beta}^{(k+1)}) \leq (\phi_c + \rho_c) \Big\{ \frac{2}{\phi_c} [F(\boldsymbol{\beta}^{(k)}, \boldsymbol{\lambda}) - F(\widehat{\boldsymbol{\beta}}, \boldsymbol{\lambda})] \Big\}^{1/2} \leq \frac{\phi_c + \rho_c}{\sqrt{k}} \|\widehat{\boldsymbol{\beta}}\|_2,$

where we have the used the fact that $\widetilde{\boldsymbol{\beta}}^{(0)} = \mathbf{0}$. To further simplify the above bound, we observe that $\phi_c \leq \gamma_u \rho_c$. Using triangle inequality, we have

$$\omega_{\boldsymbol{\lambda}}(\boldsymbol{\beta}^{(k+1)}) \leq \frac{(1+\gamma_u)\rho_c}{\sqrt{k}} \|\widehat{\boldsymbol{\beta}}\|_2 \leq \frac{(1+\gamma_u)\rho_c}{\sqrt{k}} (\|\boldsymbol{\beta}^*\|_2 + \|\widehat{\boldsymbol{\beta}} - \boldsymbol{\beta}^*\|_2).$$

Taking $\ell = 1$ and $\varepsilon = 0$ in Lemma 5.1, we have $\|\widehat{\boldsymbol{\beta}} - \boldsymbol{\beta}^*\|_2 \leq 18\lambda\sqrt{s}$. Plugging this back into (C.4) yields that

$$\omega_{\boldsymbol{\lambda}}(\boldsymbol{\beta}^{(k+1)}) \leq \frac{(1+\gamma_u)\rho_c}{\sqrt{k}} \{\|\boldsymbol{\beta}^*\|_2 + 18\lambda\sqrt{s}\} \leq \frac{(1+\gamma_u)R\rho_c}{2\sqrt{k}},$$

where $R = 2(\|\boldsymbol{\beta}^*\|_2 + 18\lambda\sqrt{s}) \lesssim \|\boldsymbol{\beta}^*\|_2 + \lambda\sqrt{s}$. Therefore, in the contraction stage, to ensure that $\omega_{\boldsymbol{\lambda}^{(0)}}(\boldsymbol{\beta}^{(k+1)}) \leq \varepsilon$, it suffices to make $k$ satisfies that

$$\frac{(1+\gamma_u)R\rho_c}{2\sqrt{k}} \leq \varepsilon_c, \quad \text{which implies} \quad k \geq \Big( \frac{(1+\gamma_u)R\rho_c}{\varepsilon_c} \Big)^2.$$

$\square$

PROOF OF PROPOSITION 4.6. Write $\varepsilon_t$ as $\varepsilon$ and assume $\ell \geq 2$. Apply Lemma E.8 in the supplement, we obtain

$$\omega_{\boldsymbol{\lambda}^{(\ell-1)}}(\boldsymbol{\beta}^{(\ell,k+1)}) \leq \big(\rho_+(2s + 2\widetilde{s}, r) + \phi_t\big) \|\boldsymbol{\beta}^{(\ell,k+1)} - \boldsymbol{\beta}^{(\ell,k)}\|_2,$$

which, combining with Lemma E.9, yields

$$\omega_{\boldsymbol{\lambda}^{(\ell-1)}}(\boldsymbol{\beta}^{(\ell,k+1)}) \leq (\phi_t + \rho_+) \sqrt{\frac{2}{\phi_t} (F(\boldsymbol{\beta}^{(\ell,k)}, \boldsymbol{\lambda}^{(\ell-1)}) - F(\boldsymbol{\beta}^{(\ell,k+1)}, \boldsymbol{\lambda}^{(\ell-1)}))}$$

$$\leq (1+\kappa) \sqrt{2\gamma_u \rho_+ (F(\boldsymbol{\beta}^{(\ell,k)}, \boldsymbol{\lambda}^{(\ell-1)}) - F(\boldsymbol{\beta}^{(\ell,k+1)}, \boldsymbol{\lambda}^{(\ell-1)}))},$$

where we use $\rho_- \leq \phi_t \leq \gamma_u \phi_t$ in the last inequality. Since the sequence $\{F(\boldsymbol{\beta}^{(\ell,k)}; \boldsymbol{\lambda}^{(\ell-1)})\}_{k=0}^{\infty}$ decrease monotonically, we obtain

$$\omega_{\boldsymbol{\lambda}^{(\ell-1)}}(\boldsymbol{\beta}^{(\ell,k+1)}) \leq (1+\kappa) \sqrt{2\gamma_u \rho_+ (F(\boldsymbol{\beta}^{(\ell,k)}, \boldsymbol{\lambda}^{(\ell-1)}) - F(\widehat{\boldsymbol{\beta}}^{(\ell)}, \boldsymbol{\lambda}^{(\ell-1)}))}$$

$$\leq (1+\kappa) \sqrt{2\gamma_u \rho_+ \Big(1 - \frac{1}{4\gamma_u \kappa}\Big)^k (F(\boldsymbol{\beta}^{(\ell,k)}, \boldsymbol{\lambda}^{(\ell-1)}) - F(\widehat{\boldsymbol{\beta}}^{(\ell)}, \boldsymbol{\lambda}^{(\ell-1)}))}$$

$$\leq (1+\kappa) \sqrt{C\gamma_u \rho_+ \Big(1 - \frac{1}{4\gamma_u \kappa}\Big)^k \lambda^2 s},$$

where the second inequality is due to Lemma 5.5, and the last one due to Lemma E.14. Here $C$ is some positive constant. Therefore, for $\ell \geq 2$, to ensure that $\boldsymbol{\beta}^{(\ell,k+1)}$ satisfies $\omega_{\boldsymbol{\lambda}^{(\ell-1)}}(\boldsymbol{\beta}^{(\ell,k+1)}) \leq \varepsilon$, it suffices to choose $k$ such that

$$(1+\kappa)\sqrt{C\gamma_u\rho_+\left(1-\frac{1}{4\gamma_u\kappa}\right)^k \lambda^2 s} \leq \varepsilon.$$

Equivalently, we obtain

$$k \geq C' \log\left(C''\frac{\lambda\sqrt{s}}{\varepsilon}\right),$$

where $C' = 2/\log(4\gamma_u\kappa/\{4\gamma_u\kappa - 1\})$, $C'' = 2(1+\kappa)\sqrt{C\gamma_u\rho_+}$. □

**C.2. Key Lemmas.** In this section, we give the proofs of the key lemmas in Section 5.2. We start with the proof of Lemma 5.3, followed by the proofs of Lemma 5.4 and Lemma 5.5.

PROOF OF LEMMA 5.3. For simplicity, we omit the super-script $\ell$ in $(\ell, k)$ and denote $\boldsymbol{\beta}^{(\ell,k)}$, $\boldsymbol{\lambda}^{(0)}$ as $\boldsymbol{\beta}^{(k)}$ and $\boldsymbol{\lambda}$ respectively. Taking $\boldsymbol{\beta} = \widehat{\boldsymbol{\beta}}$ in Lemma E.4 and simplifying the inequality, we have

$$F(\widehat{\boldsymbol{\beta}}, \boldsymbol{\lambda}) - F(\boldsymbol{\beta}^{(j)}, \boldsymbol{\lambda}) \geq \frac{\phi_c}{2}\left\{\|\boldsymbol{\beta}^{(j)} - \boldsymbol{\beta}^{(j-1)}\|_2^2 - 2\langle \boldsymbol{\beta} - \boldsymbol{\beta}^{(j-1)}, \boldsymbol{\beta}^{(j)} - \boldsymbol{\beta}^{(j-1)}\rangle\right\}$$

$$\text{(C.5)} \qquad = \frac{\phi_c}{2}\left\{\|\widehat{\boldsymbol{\beta}} - \boldsymbol{\beta}^{(j)}\|_2^2 - \|\widehat{\boldsymbol{\beta}} - \boldsymbol{\beta}^{(j-1)}\|_2^2\right\}.$$

Multiplying both sides of (C.5) by $2/\phi_c$ and summing over $j$ results

$$\frac{2}{\phi_c}\sum_{j=1}^k \left\{F(\widehat{\boldsymbol{\beta}}, \boldsymbol{\lambda}) - F(\boldsymbol{\beta}^{(j)}, \boldsymbol{\lambda})\right\} \geq \sum_{j=1}^k \left\{\|\boldsymbol{\beta}^{(j)} - \widehat{\boldsymbol{\beta}}\|_2^2 - \|\boldsymbol{\beta}^{(j-1)} - \widehat{\boldsymbol{\beta}}\|_2^2\right\},$$

or equivalently

$$\text{(C.6)} \qquad \frac{2}{\phi_c}\left\{kF(\widehat{\boldsymbol{\beta}}, \boldsymbol{\lambda}) - \sum_{j=1}^k F(\boldsymbol{\beta}^{(j)}, \boldsymbol{\lambda})\right\} \geq \|\boldsymbol{\beta}^{(k)} - \widehat{\boldsymbol{\beta}}\|_2^2 - \|\boldsymbol{\beta}^{(0)} - \widehat{\boldsymbol{\beta}}\|_2^2.$$

On the other hand, taking $\boldsymbol{\beta} = \boldsymbol{\beta}^{(k-1)}$ in Lemma E.4 and replacing $k$ with $j$ yields

$$F(\boldsymbol{\beta}^{(j-1)}, \boldsymbol{\lambda}) - F(\boldsymbol{\beta}^{(j)}, \boldsymbol{\lambda}) \geq \frac{\phi_c}{2}\|\boldsymbol{\beta}^{(j)} - \boldsymbol{\beta}^{(j-1)}\|_2^2.$$

Multiplying both sides of the inequality above by $j-1$ and summing over $j$, we obtain

$$\frac{2}{\phi_c}\sum_{j=1}^{k}\Big\{(j-1)F(\boldsymbol{\beta}^{(j-1)},\boldsymbol{\lambda})-jF(\boldsymbol{\beta}^{(j)},\boldsymbol{\lambda})+F(\boldsymbol{\beta}^{(j)},\boldsymbol{\lambda})\Big\}$$
$$\geq \sum_{j=1}^{k}(j-1)\|\boldsymbol{\beta}^{(j)}-\boldsymbol{\beta}^{(j-1)}\|_2^2,$$

or equivalently,

$$(\text{C.7}) \quad \frac{2}{\phi_c}\Big\{-kF(\boldsymbol{\beta}^{(k)},\boldsymbol{\lambda})+\sum_{j=1}^{k}F(\boldsymbol{\beta}^{(j)},\boldsymbol{\lambda})\Big\} \geq \sum_{j=1}^{k}(j-1)\|\boldsymbol{\beta}^{(j)}-\boldsymbol{\beta}^{(j-1)}\|_2^2.$$

Adding (C.6) and (C.7) together and canceling the term $2/\phi_c \sum_{j=1}^{k} F(\boldsymbol{\beta}^{(j)},\boldsymbol{\lambda})$, we obtain that

$$\frac{2k}{\phi_c}\Big\{F(\widehat{\boldsymbol{\beta}},\boldsymbol{\lambda})-F(\boldsymbol{\beta}^{(k)},\boldsymbol{\lambda})\Big\} \geq \|\boldsymbol{\beta}^{(k)}-\widehat{\boldsymbol{\beta}}\|_2^2+\sum_{j=1}^{k}(j-1)\|\boldsymbol{\beta}^{(j)}-\boldsymbol{\beta}^{(j-1)}\|_2^2-\|\boldsymbol{\beta}^{(0)}-\widehat{\boldsymbol{\beta}}\|_2^2,$$

which, multiply both sides by $-1$, yields that

$$\frac{2k}{\phi_c}\{F(\boldsymbol{\beta}^{(k)},\boldsymbol{\lambda})-F(\widehat{\boldsymbol{\beta}},\boldsymbol{\lambda})\} \leq \|\boldsymbol{\beta}^{(0)}-\widehat{\boldsymbol{\beta}}\|_2^2.$$

Therefore, the proof is completed. □

PROOF OF LEMMA 5.4. For simplicity, we omit the super-script in $\widetilde{\boldsymbol{\beta}}^{(1)}$. Define the active set $S_n$ as $\{j:|\nabla\mathcal{L}(\widetilde{\boldsymbol{\beta}})_j|=\lambda\}$. Then we must have $\{j:\widetilde{\beta}_j\neq 0\}\subseteq S_n$. It suffices to show $|S_n|\leq s+\widetilde{s}$. To achieve this goal, we decompose $S_n$ into two parts and bound the size of them separately. Specifically, let

$$S_n \subseteq S \cup \underbrace{\{j\notin S^c:|(\nabla\mathcal{L}(\widetilde{\boldsymbol{\beta}})-\nabla\mathcal{L}(\boldsymbol{\beta}^*))_j|\geq\lambda/2\}}_{S_n^1} \cup \underbrace{\{j\notin S^c:|\nabla\mathcal{L}(\boldsymbol{\beta}^*)_j|>\lambda/2\}}_{S_n^2}.$$

For $S_n^2$, the assumption that $\|\nabla\mathcal{L}(\boldsymbol{\beta}^*)\|_\infty+\varepsilon\leq\lambda/4$, implies $S_n^2=\emptyset$ and thus $|S_n^2|=0$. For $S_n^1$, consider $S'$ with maximum size $s'=|S'|\leq\widetilde{s}$ such that $S'\subseteq S_n^1$. Then there exists a $d$-dimensional sign vector $\mathbf{u}$ satisfying $\|\mathbf{u}\|_\infty=1$ and $\|\mathbf{u}\|_0=s'$, such that

$$\lambda s'/2 \leq \mathbf{u}^{\mathrm{T}}\big(\nabla\mathcal{L}(\widetilde{\boldsymbol{\beta}})-\nabla\mathcal{L}(\boldsymbol{\beta}^*)\big).$$

Then, by the Mean Value theorem, there exist some $\gamma \in [0,1]$ such that

$$\nabla \mathcal{L}(\widetilde{\boldsymbol{\beta}}) - \nabla \mathcal{L}(\boldsymbol{\beta}^*) = \nabla \mathcal{L}^2(\gamma \widetilde{\boldsymbol{\beta}} + (1-\gamma)\boldsymbol{\beta}^*)(\widetilde{\boldsymbol{\beta}} - \boldsymbol{\beta}^*) \equiv \mathbf{H}(\widetilde{\boldsymbol{\beta}} - \boldsymbol{\beta}^*).$$

Here $\mathbf{H} = \nabla \mathcal{L}^2(\gamma \widetilde{\boldsymbol{\beta}} + (1-\gamma)\boldsymbol{\beta}^*)$. Writing $\mathbf{u}^{\mathrm{T}}\mathbf{H}(\widetilde{\boldsymbol{\beta}} - \boldsymbol{\beta}^*)$ as $\langle \mathbf{H}^{1/2}\mathbf{u}, \mathbf{H}^{1/2}(\widetilde{\boldsymbol{\beta}} - \boldsymbol{\beta}^*) \rangle$ and applying the Cauchy-Schwartz inequality, we have

$$(C.8) \qquad \lambda s'/2 \leq \langle \mathbf{H}^{1/2}\mathbf{u}, \mathbf{H}^{1/2}(\widetilde{\boldsymbol{\beta}} - \boldsymbol{\beta}^*) \rangle \leq \underbrace{\|\mathbf{H}^{1/2}\mathbf{u}\|_2}_{\mathrm{I}} \underbrace{\|\mathbf{H}^{1/2}(\widetilde{\boldsymbol{\beta}} - \boldsymbol{\beta}^*)\|_2}_{\mathrm{II}}.$$

Now we bound terms I and II respectively. Since $\widetilde{\boldsymbol{\beta}}, \boldsymbol{\beta}^* \in B_2(r, \boldsymbol{\beta}^*)$, any convex combination of $\widetilde{\boldsymbol{\beta}}, \boldsymbol{\beta}^*$ also falls in $B_2(r, \boldsymbol{\beta}^*)$. The localized sparse eigenvalue condition be used on $\mathbf{H}$. For I, it follows from Definition 4.1 that

$$\|\mathbf{H}^{1/2}\mathbf{u}\|_2 \leq \sqrt{\rho_+(s',r)}\|\mathbf{u}\|_2 \leq \sqrt{\rho_+(s',r)}\{\|\mathbf{u}\|_1\|\mathbf{u}\|_\infty\}^{1/2} \leq \sqrt{\rho_+(s',r)}\sqrt{s'}.$$

For II, write $C = \rho_*^{-1/2}$. It follows from Lemma E.6 in the supplement that the following inequality holds

$$\|\mathbf{H}^{1/2}(\widetilde{\boldsymbol{\beta}} - \boldsymbol{\beta}^*)\|_2^2 = \langle \nabla \mathcal{L}(\widetilde{\boldsymbol{\beta}}) - \nabla \mathcal{L}(\boldsymbol{\beta}^*), \widetilde{\boldsymbol{\beta}} - \boldsymbol{\beta}^* \rangle \leq C\lambda^2 s.$$

Thus by plugging the bounds for I and II back into (C.8), we obtain

$$\lambda s'/2 \leq \sqrt{\rho_+(s',r)}\sqrt{s'} \times C\lambda\sqrt{s}.$$

Multiplying both sides of the above inequality by $(\lambda/2)^{1/2}$ and taking squares results

$$(C.9) \qquad s' \leq 4C\rho_+(s',r)s \leq 4C\rho_+(\widetilde{s},r)s < \widetilde{s}.$$

where the last inequality is due to the assumption. Because $s' = |S'|$ achieves the maximum possible value such that $s' \leq \widetilde{s}$ for any subset $S'$ of $S_n^1$ and (C.9) shows that $s' < \widetilde{s}$, we must have $S' = S_n^1$, and thus

$$|S_n^1| = s' \leq \lfloor 4C\rho_+(\widetilde{s},r)s \rfloor < \widetilde{s}.$$

This proves the desired result. $\square$

PROOF OF LEMMA 5.5. For notational simplicity, we omit the tightening step index $\ell$ in $\boldsymbol{\beta}^{(\ell,k)}, \boldsymbol{\lambda}^{(\ell)}, \mathcal{E}_\ell^c$; and write $\boldsymbol{\beta}^{(\ell,k)}, \boldsymbol{\lambda}^{(\ell)}, \mathcal{E}_\ell^c$ as $\boldsymbol{\beta}^{(k)}, \boldsymbol{\lambda}$ and $\mathcal{E}_\ell^c$ respectively. Define $\boldsymbol{\beta}(\alpha) = \alpha\widehat{\boldsymbol{\beta}} + (1-\alpha)\boldsymbol{\beta}^{(k-1)}$. Since $F(\boldsymbol{\beta}^{(k)}, \boldsymbol{\lambda})$ is majorized at $\Psi(\boldsymbol{\beta}^{(k)}, \boldsymbol{\beta}^{(k-1)})$, we have

$$F(\boldsymbol{\beta}^{(k)}, \boldsymbol{\lambda}) \leq \min_{\boldsymbol{\beta}(\alpha)} \left\{ \mathcal{L}(\boldsymbol{\beta}^{(k-1)}) + \langle \nabla \mathcal{L}, \boldsymbol{\beta} - \boldsymbol{\beta}^{(k-1)} \rangle + \frac{\phi_t}{2}\|\boldsymbol{\beta} - \boldsymbol{\beta}^{(k-1)}\|_2^2 + \|\boldsymbol{\lambda} \odot \boldsymbol{\beta}\|_1 \right\}$$

$$\leq \min_{\boldsymbol{\beta}(\alpha)} \left\{ F(\boldsymbol{\beta}, \boldsymbol{\lambda}) + \frac{\phi_t}{2}\|\boldsymbol{\beta} - \boldsymbol{\beta}^{(k-1)}\|_2^2 \right\},$$

where we restrict $\boldsymbol{\beta}$ on the line segment $\alpha\widehat{\boldsymbol{\beta}} + (1-\alpha)\boldsymbol{\beta}^{(k-1)}$ in the first inequality and the last inequality follows from the convexity of $\mathcal{L}(\boldsymbol{\beta})$. Let $\widehat{\boldsymbol{\beta}}^{(\ell)}$ be a solution to $\operatorname{argmin}_{\boldsymbol{\beta}\in\mathbb{R}^d}\{\mathcal{L}(\boldsymbol{\beta})+\|\boldsymbol{\lambda}^{(\ell-1)}\odot\boldsymbol{\beta}\|_1\}$. Using the convexity of $F(\boldsymbol{\beta},\boldsymbol{\lambda})$, we obtain that

$$F(\boldsymbol{\beta}^{(k)},\boldsymbol{\lambda}) \leq \min_{\boldsymbol{\beta}(\alpha)}\left\{F(\boldsymbol{\beta},\boldsymbol{\lambda}) + \frac{\phi_t}{2}\|\boldsymbol{\beta}-\boldsymbol{\beta}^{(k-1)}\|_2^2\right\}$$

$$\leq \min_{\alpha}\left\{\alpha F(\widehat{\boldsymbol{\beta}},\boldsymbol{\lambda}) + (1-\alpha)F(\boldsymbol{\beta}^{(k-1)},\boldsymbol{\lambda}) + \frac{\alpha^2\phi_t}{2}\|\boldsymbol{\beta}^{(k-1)}-\widehat{\boldsymbol{\beta}}\|_2^2\right\}$$

$$\leq \min_{\alpha}\left\{F(\boldsymbol{\beta}^{(k-1)},\boldsymbol{\lambda}) - \alpha\bigl[F(\boldsymbol{\beta}^{(k-1)},\boldsymbol{\lambda}) - F(\widehat{\boldsymbol{\beta}},\boldsymbol{\lambda})\bigr] + \frac{\alpha^2\phi_t}{2}\|\boldsymbol{\beta}^{(k-1)}-\widehat{\boldsymbol{\beta}}\|_2^2\right\}.$$

Next, we bound the last term in the inequality above. Applying Lemma E.14 in the supplementary material, we obtain

$$\|(\boldsymbol{\beta}^{(k-1)})_{\mathcal{E}^c}\|_0 \leq \widetilde{s},\ \|\boldsymbol{\beta}^{(k-1)}-\boldsymbol{\beta}^*\|_2 \leq C'\lambda\sqrt{s} \leq r,\ \|\widehat{\boldsymbol{\beta}}-\boldsymbol{\beta}^*\|_2 \leq r,\text{ and }\|\widehat{\boldsymbol{\beta}}_{\mathcal{E}}\|_0 \leq \widetilde{s}.$$

Recall $\widehat{\boldsymbol{\xi}}$ is some subgradient of $\|\widehat{\boldsymbol{\beta}}\|_1$. Using the convexity of $\mathcal{L}(\cdot)$ and the $\ell_1$-norm, $F(\boldsymbol{\beta}^{(k-1)},\boldsymbol{\lambda})-F(\widehat{\boldsymbol{\beta}},\boldsymbol{\lambda})$ can be bounded in the following way

$$F(\boldsymbol{\beta}^{(k-1)},\boldsymbol{\lambda})-F(\widehat{\boldsymbol{\beta}},\boldsymbol{\lambda}) \geq \langle\nabla\mathcal{L}(\widehat{\boldsymbol{\beta}})+\boldsymbol{\lambda}\odot\widehat{\boldsymbol{\xi}},\boldsymbol{\beta}^{(k-1)}-\widehat{\boldsymbol{\beta}}\rangle + D_{\mathcal{L}}(\boldsymbol{\beta}^{(k-1)},\widehat{\boldsymbol{\beta}})$$

$$\geq \frac{\rho_-}{2}\|\boldsymbol{\beta}^{(k-1)}-\widehat{\boldsymbol{\beta}}\|_2^2,$$

where the last inequality is due to the first order optimality condition and Proposition B.3. Thus we conclude that

$$F(\boldsymbol{\beta}^{(k)},\boldsymbol{\lambda}) \leq \min\left\{F(\boldsymbol{\beta}^{(k-1)},\boldsymbol{\lambda})-\alpha\bigl[F(\boldsymbol{\beta}^{(k-1)},\boldsymbol{\lambda})-F(\widehat{\boldsymbol{\beta}},\boldsymbol{\lambda})\bigr]\right.$$

$$\left.+\frac{\alpha^2\phi_t}{\rho_-}\bigl[F(\boldsymbol{\beta}^{(k-1)},\boldsymbol{\lambda})-F(\widehat{\boldsymbol{\beta}},\boldsymbol{\lambda})\bigr]\right\}$$

$$\leq F(\boldsymbol{\beta}^{(k-1)},\boldsymbol{\lambda}) - \frac{\rho_-}{4\phi_t}\bigl[F(\boldsymbol{\beta}^{(k-1)},\boldsymbol{\lambda})-F(\widehat{\boldsymbol{\beta}},\boldsymbol{\lambda})\bigr].$$

which, combining with the fact $\phi_t \leq \gamma_u\rho_+$, yields

$$F(\boldsymbol{\beta}^{(k)},\boldsymbol{\lambda}) - F(\widehat{\boldsymbol{\beta}},\boldsymbol{\lambda}) \leq \left(1-\frac{1}{4\gamma_u\kappa}\right)^k\bigl\{F(\widetilde{\boldsymbol{\beta}}^{(0)},\boldsymbol{\lambda}) - F(\widehat{\boldsymbol{\beta}},\boldsymbol{\lambda})\bigr\},$$

in which $\kappa = \rho_+/\rho_-$.

$\square$

## APPENDIX D: AN APPLICATION TO ROBUST LINEAR REGRESSION

In this section, we give an application of Theorem 4.2 to robust linear regression. The Huber loss, defined in Section A, is used to robustify the heavy tailed error. We allow the cutoff parameter $\alpha$ to scale with $(n, d, s)$ for bias-robustness tradeoff. Let $y_i = \mathbf{x}_i^\top \boldsymbol{\beta} + \epsilon_i$, $1 \leq i \leq n$, be independently and identically distributed random variables, with mean $\mu = \mathbf{x}_i^\top \boldsymbol{\beta}$ and finite second moment $M$. Then, the following corollary suggests that, under only finite second moments, the sparse Huber estimator with an adaptive $\alpha$ can perform as good as the sparse ordinary least square estimator as if sub-Gaussian errors were assumed.

**Corollary D.1.** Suppose the same conditions in Theorem 4.2 hold. Assume the columns of $\mathbf{X}$ are normalized such that $\max_j \|\mathbf{X}_{*j}\|_2 \leq \sqrt{n}$. Assume there exists an $\alpha > 0$ such that $\|\boldsymbol{\beta}_S^*\|_{\min} \geq u + \gamma\lambda$ and $\mathrm{w}(\gamma\lambda) = 0$. If $\alpha \propto \lambda \propto \sqrt{n^{-1}\log d}$, $\varepsilon_t \leq \sqrt{1/n}$ and $T \gtrsim \log\log d$, then with probability at least $1 - 2d^{-\eta_1} - 2\exp(-\eta_2 s)$, $\widetilde{\boldsymbol{\beta}}^{(T)}$ must satisfy

$$\|\widetilde{\boldsymbol{\beta}}^{(T)} - \boldsymbol{\beta}^*\|_2 \lesssim \sqrt{s/n},$$

where $\eta_1$ and $\eta_2$ are positive constants.

PROOF OF COROLLARY D.1. The proof follows from that of Corollary 4.3 by bounding $\|\nabla \mathcal{L}(\boldsymbol{\beta}^*)_S\|_\infty$ and the probability of the event $\{\|\nabla \mathcal{L}(\boldsymbol{\beta}^*)\|_\infty \gtrsim \lambda\}$. The derivative $\nabla \mathcal{L}(\boldsymbol{\beta}^*)$ can be written as

$$\nabla \mathcal{L}(\boldsymbol{\beta}^*) = \frac{1}{n} \sum_{i=1}^n \nabla \ell_\alpha(\epsilon_i) \mathbf{x}_i = \frac{1}{n} \mathbf{X}^\mathrm{T} \boldsymbol{\epsilon}_\alpha,$$

where $\boldsymbol{\epsilon}_\alpha = (\epsilon_{1,\alpha}, \ldots, \epsilon_{n,\alpha})^T$. Therefore, it suffices to show that $\epsilon_\alpha$ has sub-Gaussian tail. Let $\psi(\alpha x) = 2^{-1}\alpha \nabla \ell_\alpha(x)$. Then $\psi(x)$ satisfies that

$$-\log(1 - x + x^2) \leq \psi(x) \leq \log(1 + x + x^2),$$

which yields that

$$\mathbb{E}\big[\exp\{\psi(\alpha\epsilon)\}\big] \leq 1 + \alpha^2 M \quad \text{and} \quad \mathbb{E}\big[\exp\{-\psi(\alpha\epsilon)\}\big] \leq 1 + \alpha^2 M.$$

Using Markov inequality, we obtain,

$$\mathbb{P}\big(\psi(\alpha\epsilon) \geq Mt^2\big) \leq \frac{\mathbb{E}\big[\exp\{\psi(\alpha\epsilon)\}\big]}{\exp(Mt^2)} \leq (1 + \alpha^2 M)\exp\{-Mt^2\},$$

or equivalently

$$\mathbb{P}\big(\nabla\ell_\alpha(\epsilon) \geq 2Mt^2/\alpha\big) \leq \frac{\mathbb{E}\big[\exp\{\psi(\alpha\epsilon)\}\big]}{\exp(Mt^2)} \leq (1+\alpha^2 M)\exp\{-Mt^2\}.$$

Taking $\alpha = t/2$, we obtain

$$\mathbb{P}\big(\nabla\ell_\alpha(\epsilon) \geq Mt\big) \leq \frac{\mathbb{E}\big[\exp\{\psi(t\epsilon/2)\}\big]}{\exp(Mt^2)} \leq (1+Mt^2/4)\exp\{-Mt^2\}$$
$$\leq \exp\big\{-Mt^2/2\big\},$$

where the last inequality follows from the fact that $1+Mt^2/4 \leq \exp\{Mt^2/2\}$. The rest of the proof follows from that of Corollary 4.3. $\square$

## APPENDIX E: TECHNICAL LEMMAS

**E.1. Statistical Theory.** We collect the technical lemmas that are used to prove Theorem 4.2. We start by defining the following localized sparse relative covariance.

**Definition E.1** (Localized Sparse Relative Covariance). The localized sparse relative covariance with parameter $r$ is defined as

$$\pi(i,j;\boldsymbol{\beta}^*,r) = \sup_{\mathbf{v},\mathbf{u},\|\boldsymbol{\beta}-\boldsymbol{\beta}^*\|_2 \leq r} \left\{ \frac{\mathbf{v}_I^T \nabla^2 \mathcal{L}(\boldsymbol{\beta})\mathbf{u}_J/\|\mathbf{u}_J\|_2}{\mathbf{v}_I^T \nabla^2 \mathcal{L}(\boldsymbol{\beta})\mathbf{v}_I/\|\mathbf{v}_I\|_2} : I \cap J = \varnothing, |I| \leq i, |J| \leq j \right\}.$$

This is different from restricted correlation defined in [4]. We measure the relative covariance between set $I$ and set $J$ with respect to that of set $I$. In the sequel, we omit the arguments $\boldsymbol{\beta}^*, r$ in $\pi(i,j;\boldsymbol{\beta}^*,r)$ for simplicity. Our next result bounds the localized sparse relative covariance in terms of sparse eigenvalues.

**Lemma E.1.** It holds that

$$\pi(i,j;\boldsymbol{\beta}^*,r) \leq \frac{1}{2}\sqrt{\frac{\rho_+(j,r)}{\rho_-(i+j,r)} - 1}.$$

Then we are ready to bound the estimation error by functionals of the regularization parameter under localized sparse eigenvalue condition, which is proved in the following lemma.

**Lemma E.2.** Take $\mathcal{E}$ such that $S \subseteq \mathcal{E}$ and $|\mathcal{E}| = k \leq 2s$. Let $J$ be the index set of the largest $m$ coefficients (in absolute value) in $\mathcal{E}^c$. Assume Assumption 4.1 holds, $\|\nabla \mathcal{L}(\boldsymbol{\beta}^*)\|_\infty + \varepsilon \leq \lambda/4$, $\|\boldsymbol{\lambda}_{\mathcal{E}^c}\|_{\min} \geq \lambda/2$ and $r\rho_-(k+m,r) > 2(1+5\sqrt{k/m})(\sqrt{k/s}/4+1)\lambda\sqrt{s}$. Then any $\varepsilon$-optimal solution $\widetilde{\boldsymbol{\beta}}$ must satisfy

$$\|\widetilde{\boldsymbol{\beta}} - \boldsymbol{\beta}^*\|_2 \leq \frac{2(1+5\sqrt{k/m})}{\rho_-(k+m,r)}(\|\boldsymbol{\lambda}_S\|_2 + \|\nabla \mathcal{L}(\boldsymbol{\beta}^*)_\mathcal{E}\|_2 + \sqrt{k}\varepsilon) \lesssim \lambda\sqrt{s}.$$

We now present the proofs for the two lemmas above by starting with proof of Lemma E.1.

PROOF OF LEMMA E.1. For simplicity, we omit the arguments $\boldsymbol{\beta}^*, r$ in $\pi(i,j;\boldsymbol{\beta}^*,r)$, $\rho_+(i,,r)$, and $\rho_-(i,r)$. Let $I = \mathcal{E} \cup J$ and $L = I \cup J$. For any $\alpha \in \mathbb{R}$, let $w = \mathbf{v}_I + \alpha \mathbf{u}_J$. Without loss of generality, we assume that $\|\mathbf{u}_J\|_2 = 1$, $\|\mathbf{v}_I\|_2 = 1$ and $\boldsymbol{\beta} \in B_2(r, \boldsymbol{\beta}^*)$. Using the definition of $\pi(i,j)$ and w, we have

$$\rho_-(i+j)\|w\|_2^2 \leq \underbrace{\mathbf{v}_I^T \nabla^2 \mathcal{L}(\boldsymbol{\beta})\mathbf{v}_I}_{c_1} + 2\alpha \underbrace{\mathbf{v}_I \nabla^2 \mathcal{L}(\boldsymbol{\beta})\mathbf{u}_J}_{b} + \alpha^2 \underbrace{\mathbf{u}_J \nabla^2 \mathcal{L}(\boldsymbol{\beta})\mathbf{u}_J}_{c_2},$$

which simplifies to

(E.1) $$(c_2 - \rho_-(i+j))\alpha^2 + 2b\alpha + (c_1 - \rho_-(i+j)) \geq 0.$$

Since the left hand side of (E.1) is positive semidefinite for all $\alpha$, we must have

$$(c_2 - \rho_-(i+j))(c_1 - \rho_-(i+j)) \geq b^2.$$

Multiplying by $4/c_1^2$ on both sides of the inequality above, we obtain

$$\frac{4b^2}{c_1^2} \leq 4c_1^{-1}(1 - \rho_-(i+j)/c_1)(c_2 - \rho_-(i+j))$$

$$\leq 4c_1^{-1}\rho_-(i+j)(1 - \rho_-(i+j)/c_1) \times \frac{c_2 - \rho_-(i+j)}{\rho_-(i+j)}$$

$$\leq \frac{c_2 - \rho_-(i+j)}{\rho_-(i+j)} \leq \frac{\rho_+(j)}{\rho_-(i+j)} - 1,$$

where ,in the last second inequality, we use $4c_1^{-1}\rho_-(i+j)(1-\rho_-(i+j)/c_1) \leq 1$; and the last inequality is due to $c_2 = \mathbf{u}_J \nabla^2 \mathcal{L}(\boldsymbol{\beta})\mathbf{u}_J \leq \rho_+(j)$. This yields

$$\frac{\mathbf{v}_I^T \nabla^2 \mathcal{L}(\boldsymbol{\beta})\mathbf{u}_J/\|\mathbf{u}_J\|_2}{\mathbf{v}_I^T \nabla \mathcal{L}(\boldsymbol{\beta})\mathbf{v}_I/\|\mathbf{v}_I\|_2} \leq \frac{|\mathbf{v}_I^T \nabla^2 \mathcal{L}(\boldsymbol{\beta})\mathbf{u}_J|}{\mathbf{v}_I^T \nabla \mathcal{L}(\boldsymbol{\beta})\mathbf{v}_I} \leq \frac{1}{2}\sqrt{\frac{\rho_+(j)}{\rho_-(i+j)} - 1}.$$

The proof is completed by taking sup of the left hand side with respect to $\boldsymbol{\beta}, \mathbf{u}, \mathbf{v}$. □

PROOF OF LEMMA E.2. For simplicity, we write $\nabla^2 \mathcal{L}(\boldsymbol{\beta})$ as $\nabla^2 \mathcal{L}$, whenever we have $\boldsymbol{\beta} \in B_2(r, \boldsymbol{\beta}^*) \equiv \{\boldsymbol{\beta} : \|\boldsymbol{\beta} - \boldsymbol{\beta}^*\|_2 \leq r\}$. Since we do not know whether $\widetilde{\boldsymbol{\beta}}$ belongs to $B_r^2(\boldsymbol{\beta}^*)$ in advance. we need to construct an intermediate estimator $\widetilde{\boldsymbol{\beta}}^*$ such that $\|\widetilde{\boldsymbol{\beta}}^* - \boldsymbol{\beta}^*\|_2 \leq r$. Let $\widetilde{\boldsymbol{\beta}}^\star = \boldsymbol{\beta}^* + t(\widetilde{\boldsymbol{\beta}} - \boldsymbol{\beta}^*)$ where $t = 1$ if $\|\widetilde{\boldsymbol{\beta}} - \boldsymbol{\beta}^*\| \leq r$; $t \in (0, 1)$ such that $\|\widetilde{\boldsymbol{\beta}}^* - \boldsymbol{\beta}^*\|_2 = r$ otherwise. Using the Cauchy Schwartz inequality, we obtain

$$\text{(E.2)} \quad \|(\widetilde{\boldsymbol{\beta}} - \boldsymbol{\beta}^*)_{I^c}\|_2 \leq (\underbrace{\|(\widetilde{\boldsymbol{\beta}} - \boldsymbol{\beta}^*)_{I^c}\|_1}_{\text{I}} \underbrace{\|(\widetilde{\boldsymbol{\beta}} - \boldsymbol{\beta}^*)_{I^c}\|_\infty}_{\text{II}})^{1/2}.$$

We bound I and II respectively. For I, since $I^c \subseteq \mathcal{E}^c$, we apply Lemma F.1 and obtain

$$\text{(E.3)} \quad \|(\widetilde{\boldsymbol{\beta}} - \boldsymbol{\beta}^*)_{I^c}\|_1 \leq \|(\widetilde{\boldsymbol{\beta}} - \boldsymbol{\beta}^*)_{\mathcal{E}^c}\|_1 \leq 5\|(\widetilde{\boldsymbol{\beta}} - \boldsymbol{\beta}^*)_\mathcal{E}\|_1.$$

For II, note that $I = \mathcal{E} \cup J$ and using the definition of $J$, we obtain

$$\text{(E.4)} \quad \|(\widetilde{\boldsymbol{\beta}} - \boldsymbol{\beta}^*)_{I^c}\|_\infty \leq \|(\widetilde{\boldsymbol{\beta}} - \boldsymbol{\beta}^*)_{\mathcal{E}^c}\|_1/m \leq \frac{5}{m}\|(\widetilde{\boldsymbol{\beta}} - \boldsymbol{\beta}^*)_\mathcal{E}\|_1.$$

Plugging (E.3) and (E.4) into (E.2) results

$$\|(\widetilde{\boldsymbol{\beta}} - \boldsymbol{\beta}^*)_{I^c}\|_2 \leq 5\sqrt{\frac{1}{m}}\|(\widetilde{\boldsymbol{\beta}} - \boldsymbol{\beta}^*)_\mathcal{E}\|_1 \leq 5\sqrt{\frac{k}{m}}\|(\widetilde{\boldsymbol{\beta}} - \boldsymbol{\beta}^*)_I\|_2.$$

Using triangle inequality along with the result above yields

$$\|\widetilde{\boldsymbol{\beta}} - \boldsymbol{\beta}^*\|_2 \leq \|(\widetilde{\boldsymbol{\beta}} - \boldsymbol{\beta}^*)_I\|_2 + \|(\widetilde{\boldsymbol{\beta}} - \boldsymbol{\beta}^*)_{I^c}\|_2 \leq (1 + c_0\sqrt{k/m})\|(\widetilde{\boldsymbol{\beta}} - \boldsymbol{\beta}^*)_I\|_2.$$

Since $\widetilde{\boldsymbol{\beta}}^* - \boldsymbol{\beta}^* = t(\widetilde{\boldsymbol{\beta}} - \boldsymbol{\beta}^*)$, we have

$$\text{(E.5)} \quad \|\widetilde{\boldsymbol{\beta}}^* - \boldsymbol{\beta}^*\|_2 = t\|\widetilde{\boldsymbol{\beta}} - \boldsymbol{\beta}^*\|_2 \leq (1 + 5\sqrt{k/m})\|(\widetilde{\boldsymbol{\beta}}^* - \boldsymbol{\beta}^*)_I\|_2.$$

Thus to bound $\|\widetilde{\boldsymbol{\beta}}^* - \boldsymbol{\beta}^*\|_2$, it suffices to bound $\|(\widetilde{\boldsymbol{\beta}}^* - \boldsymbol{\beta}^*)_I\|_2$.

**Bounding $\|(\widetilde{\boldsymbol{\beta}}^* - \boldsymbol{\beta}^*)_I\|_2$ by $D_\mathcal{L}^s(\widetilde{\boldsymbol{\beta}}, \boldsymbol{\beta}^*)$:**

For notational simplicity, we write $\mathbf{u} = \widetilde{\boldsymbol{\beta}}^* - \boldsymbol{\beta}^*$ sometimes. Let $\mathbf{u} = (u_1, u_2, \ldots, u_p)^T$. Without loss of generality, we assume the first $k$ elements of $\boldsymbol{\beta}^*$ contains the true support $S$. When $j > k$, $u_j$ is arranged such that $|u_{k+1}| \geq |u_{k+2}| \cdots \geq |u_p|$. Let $J_0 = \mathcal{E} = \{1, \ldots, k\}$, and $J_i = \{k+(i-1)m+1, \ldots, k+im\}$, for $i = 1, 2, \ldots$, with the size of last block smaller or equal than $m$. In this manner, we have $J_1 = J$ and $I = J_0 \cup J_1$. Moreover, we have $\|\mathbf{u}_{J_i}\|_2 \leq \|\mathbf{u}_{J_i}\|_\infty \sqrt{m} \leq \|\mathbf{u}_{J_{i-1}}\|_1/\sqrt{m}$ when $i > 1$, which implies that $\sum_{i>1} \|\mathbf{u}_{J_i}\|_2 \leq \|\mathbf{u}_{\mathcal{E}^c}\|_1/\sqrt{m}$. Now if

$$\text{(E.6)} \quad 1 - 2\pi(|I|, m) m^{-1/2} \frac{\|\mathbf{u}_{\mathcal{E}^c}\|_1}{\|\mathbf{u}_I\|_2} \geq 0,$$

separating the support of $\mathbf{u}$ into $I$, $I^c$ and using $\mathbf{u}_{I^c}^T \nabla^2 \mathcal{L}(\boldsymbol{\beta}) \mathbf{u}_{I^c} \geq 0$, we obtain

$$\mathbf{u}^T \nabla^2 \mathcal{L}(\boldsymbol{\beta}) \mathbf{u} \geq \mathbf{u}_I^T \nabla^2 \mathcal{L}(\boldsymbol{\beta}) \mathbf{u}_I + 2 \sum_{i>1} \mathbf{u}_I^T \nabla^2 \mathcal{L}(\boldsymbol{\beta}) \mathbf{u}_{J_i}$$
$$\geq \mathbf{u}_I^T \nabla^2 \mathcal{L}(\boldsymbol{\beta}) \mathbf{u}_I \left(1 - 2\pi(|I|, m) \sum_{i>1} \frac{\|\mathbf{u}_{J_i}\|_2}{\|\mathbf{u}_I\|_2}\right)$$
$$\geq \rho_-(|I|) \left(1 - 2\pi(|I|, m) \sqrt{\frac{1}{m}} \frac{\|\mathbf{u}_{\mathcal{E}^c}\|_1}{\|\mathbf{u}_I\|_2}\right) \|\mathbf{u}_I\|_2^2,$$

where we use the definition of $\pi(|I|, m)$ and $\rho_-(|I|)$ in the last two inequalities. Notice that $\|(\widetilde{\boldsymbol{\beta}} - \boldsymbol{\beta}^*)_{\mathcal{E}}\|_1 \leq \sqrt{k} \|(\widetilde{\boldsymbol{\beta}} - \boldsymbol{\beta}^*)_{\mathcal{E}}\|_2$ and applying Lemma F.1, we obtain

$$\|(\widetilde{\boldsymbol{\beta}}^* - \boldsymbol{\beta}^*)_{\mathcal{E}^c}\|_1 = t \|(\widetilde{\boldsymbol{\beta}} - \boldsymbol{\beta}^*)_{\mathcal{E}^c}\|_1 \leq 5 \times t \|(\widetilde{\boldsymbol{\beta}} - \boldsymbol{\beta}^*)_{\mathcal{E}}\|_1 \leq 5\sqrt{k} \times \|(\widetilde{\boldsymbol{\beta}}^* - \boldsymbol{\beta}^*)_{\mathcal{E}}\|_2.$$

Further note that $\mathcal{E} \subseteq I$, we obtain

$$1 - 2\pi(k+m, m) \sqrt{\frac{1}{m}} \frac{\|(\widetilde{\boldsymbol{\beta}}^* - \boldsymbol{\beta}^*)_{\mathcal{E}^c}\|_1}{\|(\widetilde{\boldsymbol{\beta}}^* - \boldsymbol{\beta}^*)_I\|_2} \geq 1 - 10\pi(k+m, m) \times \sqrt{\frac{k}{m}}.$$

Using $\pi(k+m, m) \leq 2^{-1} \left(\rho_+(m)/\rho_-(k+2m) - 1\right)^{1/2}$ and Assumption 4.1 with $c = 100$ results

$$1 - 2\pi(k+m, m) m^{-1} \frac{\|(\widetilde{\boldsymbol{\beta}}^* - \boldsymbol{\beta}^*)_{\mathcal{E}^c}\|_1}{\|(\widetilde{\boldsymbol{\beta}}^* - \boldsymbol{\beta}^*)_I\|_2} \geq 1 - 5\sqrt{\frac{k}{m}} \sqrt{\frac{\rho_+(m)}{\rho_-(k+2m)} - 1} \geq 1/2$$

Therefore for any $\boldsymbol{\beta} \in B_2(r, \boldsymbol{\beta}^*)$, we have that $(\widetilde{\boldsymbol{\beta}}^* - \boldsymbol{\beta}^*)^T [\nabla^2 \mathcal{L}(\boldsymbol{\beta})] (\widetilde{\boldsymbol{\beta}}^* - \boldsymbol{\beta}^*) \geq 1/2 \rho_-(k+m) \|(\widetilde{\boldsymbol{\beta}}^* - \boldsymbol{\beta}^*)_I\|_2^2$. By the Mean Value theorem, there exists a $\gamma \in [0, 1]$ such that

$$\langle \nabla \mathcal{L}(\widetilde{\boldsymbol{\beta}}^*) - \nabla \mathcal{L}(\boldsymbol{\beta}^*), \widetilde{\boldsymbol{\beta}}^* - \boldsymbol{\beta}^* \rangle = (\widetilde{\boldsymbol{\beta}}^* - \boldsymbol{\beta}^*)^T [\nabla^2 \mathcal{L}(\gamma \boldsymbol{\beta}^* + (1-\gamma) \widetilde{\boldsymbol{\beta}}^*)] (\widetilde{\boldsymbol{\beta}}^* - \boldsymbol{\beta}^*)$$
$$\geq 2^{-1} \rho_-(k+m) \|(\widetilde{\boldsymbol{\beta}}^* - \boldsymbol{\beta}^*)_I\|_2^2.$$

We further bound the left hand side $D_{\mathcal{L}}^s(\widetilde{\boldsymbol{\beta}}^*, \boldsymbol{\beta}^*)$ in the following.

**Bounding $D_{\mathcal{L}}^s(\widetilde{\boldsymbol{\beta}}^*, \boldsymbol{\beta}^*)$:**
Define $\mathbf{u} = \nabla \mathcal{L}(\widetilde{\boldsymbol{\beta}}) + \boldsymbol{\lambda} \odot \boldsymbol{\xi}$, where $\boldsymbol{\xi} \in \partial \|\widetilde{\boldsymbol{\beta}}\|_1$. Then by Lemma F.2 and the definition of $D_{\mathcal{L}}^s(\widetilde{\boldsymbol{\beta}}, \boldsymbol{\beta}^*)$, we obtain

$$(\text{E.7}) \qquad D_{\mathcal{L}}^s(\widetilde{\boldsymbol{\beta}}^*, \boldsymbol{\beta}^*) \leq t D_{\mathcal{L}}^s(\widetilde{\boldsymbol{\beta}}, \boldsymbol{\beta}^*) = t \langle \nabla \mathcal{L}(\widetilde{\boldsymbol{\beta}}) - \nabla \mathcal{L}(\boldsymbol{\beta}^*), \widetilde{\boldsymbol{\beta}} - \boldsymbol{\beta}^* \rangle$$

Adding and subtracting the term $t \langle \boldsymbol{\lambda} \odot \boldsymbol{\xi}, \widetilde{\boldsymbol{\beta}} - \boldsymbol{\beta}^* \rangle$, we have

$$t \langle \nabla \mathcal{L}(\widetilde{\boldsymbol{\beta}}) - \nabla \mathcal{L}(\boldsymbol{\beta}^*), \widetilde{\boldsymbol{\beta}} - \boldsymbol{\beta}^* \rangle = t \langle \nabla \mathcal{L}(\widetilde{\boldsymbol{\beta}}) + \boldsymbol{\lambda} \odot \boldsymbol{\xi}, \widetilde{\boldsymbol{\beta}} - \boldsymbol{\beta}^* \rangle \\ - t \langle \nabla \mathcal{L}(\boldsymbol{\beta}^*), \widetilde{\boldsymbol{\beta}} - \boldsymbol{\beta}^* \rangle - t \langle \boldsymbol{\lambda} \odot \boldsymbol{\xi}, \widetilde{\boldsymbol{\beta}} - \boldsymbol{\beta}^* \rangle.$$

Using a similar argument in the proof of Lemma B.7, we obtain

$$D_{\mathcal{L}}^s(\widetilde{\boldsymbol{\beta}}^*, \boldsymbol{\beta}^*) \leq \big( \|\boldsymbol{\lambda}_S\|_2 + \|\nabla \mathcal{L}(\boldsymbol{\beta}^*)_{\mathcal{E}}\|_2 + \varepsilon \sqrt{|\mathcal{E}|} \big) \|(\widetilde{\boldsymbol{\beta}}^* - \boldsymbol{\beta}^*)_I\|_2.$$

**Bounding $\|\widetilde{\boldsymbol{\beta}}^* - \boldsymbol{\beta}^*\|_2$ and $\|\widetilde{\boldsymbol{\beta}} - \boldsymbol{\beta}^*\|_2$:**
Combing the upper and lower bound for $D_L^s(\widetilde{\boldsymbol{\beta}}^*, \boldsymbol{\beta}^*)$, we have

$$\|(\widetilde{\boldsymbol{\beta}}^* - \boldsymbol{\beta}^*)_I\|_2 \leq \frac{2}{\rho_-(k+m, r)} \big( \|\boldsymbol{\lambda}_S\|_2 + \|\nabla \mathcal{L}(\boldsymbol{\beta}^*)_{\mathcal{E}}\|_2 + \varepsilon \sqrt{|\mathcal{E}|} \big).$$

Plugging the above bound into (E.5) yields

$$\|\widetilde{\boldsymbol{\beta}}^* - \boldsymbol{\beta}^*\|_2 \leq \frac{2\big(1 + 5\sqrt{k/m}\big)}{\rho_-(k+m, r)} \big( \|\boldsymbol{\lambda}_S\|_2 + \|\nabla \mathcal{L}(\boldsymbol{\beta}^*)_{\mathcal{E}}\|_2 + \sqrt{k}\varepsilon \big) < r.$$

If $t \neq 1$, by the construction of $\widetilde{\boldsymbol{\beta}}^*$, we must have $\|\widetilde{\boldsymbol{\beta}}^* - \boldsymbol{\beta}^*\|_2 = r$, which contradicts our the above inequality. Thus $t$ must be 1, which implies $\widetilde{\boldsymbol{\beta}}^* = \widetilde{\boldsymbol{\beta}}$. This completes the proof. $\square$

**E.2. Computational Theory.** In this section, we prove technical lemmas used in Appendix C.

E.2.1. *Contraction Stage.* We start with a lemma that characterizes the locality of the solution sequence. It also provides the lower and upper bounds of $\phi_c$, which will be exploited in our final localized iteration complexity analysis.

**Lemma E.3.** Under Assumption 4.3 and the same conditions of Theorem 4.2, we have

$$\|\boldsymbol{\beta}^{(1,k)} - \boldsymbol{\beta}^*\|_2 \leq R/2 \quad \text{and} \quad \phi_0 \leq \phi_c \leq \gamma_u \rho_c.$$

The next two lemmas are critical for the analysis of computational complexity in the contraction stage.

**Lemma E.4.** Recall that $F(\boldsymbol{\beta}, \boldsymbol{\lambda}) = \mathcal{L}(\boldsymbol{\beta}) + \|\boldsymbol{\lambda} \odot \boldsymbol{\beta}\|_1$. We have

$$F(\boldsymbol{\beta}; \boldsymbol{\lambda}^{(0)}) - F(\boldsymbol{\beta}^{(1,k)}; \boldsymbol{\lambda}^{(0)}) \geq 2^{-1}\phi_c \|\boldsymbol{\beta}^{(1,k)} - \boldsymbol{\beta}^{(1,k-1)}\|_2^2 \\ - \phi_c \langle \boldsymbol{\beta} - \boldsymbol{\beta}^{(1,k-1)}, \boldsymbol{\beta}^{(1,k)} - \boldsymbol{\beta}^{(1,k-1)} \rangle.$$

Our next lemma describes the relationship between suboptimality measure $\omega_{\boldsymbol{\lambda}}(\boldsymbol{\beta}^{(1,k)})$ and $\|\boldsymbol{\beta}^{(1,k)} - \boldsymbol{\beta}^{(1,k-1)}\|$, which is critical to establish the iteration complexity of the contraction stage.

**Lemma E.5.** $\omega_{\boldsymbol{\lambda}}(\boldsymbol{\beta}^{(1,k)}) \leq (\phi_c + \rho_c)\|\boldsymbol{\beta}^{(1,k)} - \boldsymbol{\beta}^{(1,k-1)}\|_2$.

PROOF OF LEMMA E.3. We first prove the second statement. If we assume for $\forall\, k \geq 1$, it holds that

(E.8) $$\|\boldsymbol{\beta}^{(k-1)} - \boldsymbol{\beta}^*\|_2 \leq R/2.$$

Then for any $\boldsymbol{\beta}$ such that $\|\boldsymbol{\beta} - \boldsymbol{\beta}^*\|_2 \leq R/2$, we have $\|\boldsymbol{\beta}^{(k-1)} - \boldsymbol{\beta}\|_2 \leq \|\boldsymbol{\beta}^{(k-1)} - \boldsymbol{\beta}^*\|_2 + \|\boldsymbol{\beta}^* - \boldsymbol{\beta}\|_2 \leq R$, by triangle inequality. Let $\mathbf{v} = \boldsymbol{\beta} - \boldsymbol{\beta}^{(k-1)}$. Using taylor expansion, we have

$$\nabla \mathcal{L}(\boldsymbol{\beta}) = \nabla \mathcal{L}(\boldsymbol{\beta}^{(k-1)}) + \langle \nabla \mathcal{L}(\boldsymbol{\beta}^{(k-1)}), \mathbf{v} \rangle \\ + \int_0^1 \langle \nabla \mathcal{L}(\boldsymbol{\beta}^{(k-1)} + t\mathbf{v}) - \nabla \mathcal{L}(\boldsymbol{\beta}^{(k-1)}), \mathbf{v} \rangle dt.$$

Applying Cauchy-Schwartz inequality and using Assumption 4.3, we obtain

$$\nabla \mathcal{L}(\boldsymbol{\beta}) \leq \nabla \mathcal{L}(\boldsymbol{\beta}^{(k-1)}) + \langle \nabla \mathcal{L}(\boldsymbol{\beta}^{(k-1)}), \mathbf{v} \rangle + \int_0^1 \rho_c t \|\mathbf{v}\|_2^2 dt \\ \leq \nabla \mathcal{L}(\boldsymbol{\beta}^{(k-1)}) + \langle \nabla \mathcal{L}(\boldsymbol{\beta}^{(k-1)}), \boldsymbol{\beta} - \boldsymbol{\beta}^{(k-1)} \rangle + \frac{\rho_c}{2}\|\boldsymbol{\beta} - \boldsymbol{\beta}^{(k-1)}\|_2^2.$$

The iterative LAMM algorithm implies that $\phi_0 \leq \phi_c \leq (1+\gamma_u)\rho_c$. Otherwise, if $\phi_c > (1+\gamma_u)\rho_c > \rho_c$, then $\phi_c' \equiv \phi_c/\gamma_u = \gamma_u^{-1}(1+\gamma_u)\rho_c$ is the quadratic parameter in the previous LAMM iteration. Let $\Psi'(\boldsymbol{\beta}; \boldsymbol{\beta}^{(1,k-1)})$ t be the corresponding local quadratic approximation. Then for any $\boldsymbol{\beta} \in B_2(R/2, \boldsymbol{\beta}^*)$, it holds that

$$\Psi'(\boldsymbol{\beta}; \boldsymbol{\beta}^{(1,k-1)}) = \mathcal{L}(\boldsymbol{\beta}^{(1,k-1)}) + \langle \nabla \mathcal{L}(\boldsymbol{\beta}^{(1,k-1)}), \boldsymbol{\beta} - \boldsymbol{\beta}^{(1,k-1)} \rangle \\ + \frac{\phi_c'}{2}\|\boldsymbol{\beta} - \boldsymbol{\beta}^{(1,k-1)}\| + \|\boldsymbol{\lambda}^{(0)} \odot \boldsymbol{\beta}\|_1 \geq F(\boldsymbol{\beta}, \boldsymbol{\lambda}^{(0)}).$$

However by the stopping rule of the I-LAMM algorithm, we must have $\Psi'(\boldsymbol{\beta}; \boldsymbol{\beta}^{(1,k-1)}) < F(\boldsymbol{\beta}, \boldsymbol{\lambda}^{(0)})$. This contradiction shows that $\phi_c \leq (1+\gamma_u)\rho_c$, and $\phi_0 < \phi_c$ can be ensured by taking $\phi_0$ small enough.

Therefore, it remains to show (E.8) holds by induction. For $k=1$, it obviously holds. Now suppose that $\|\boldsymbol{\beta}^{(k-1)} - \boldsymbol{\beta}^*\|_2 \leq R/2$. Taking $\boldsymbol{\beta} = \widehat{\boldsymbol{\beta}}$ in Lemma E.4, we obtain

$$0 \geq F(\widehat{\boldsymbol{\beta}}, \boldsymbol{\lambda}) - F(\boldsymbol{\beta}^{(j-1)}, \boldsymbol{\lambda}) \geq \frac{\phi_c}{2}\Big\{\|\widehat{\boldsymbol{\beta}} - \boldsymbol{\beta}^{(j-1)}\|_2^2 - \|\widehat{\boldsymbol{\beta}} - \boldsymbol{\beta}^{(0)}\|_2^2\Big\},$$

which implies

(E.9) $$\|\boldsymbol{\beta}^{(j)} - \widehat{\boldsymbol{\beta}}\|_2 \leq \|\boldsymbol{\beta}^{(j-1)} - \widehat{\boldsymbol{\beta}}\|_2.$$

Taking $j = 1, \ldots, k$, repeating (E.9) yields that

$$\|\boldsymbol{\beta}^{(k)} - \widehat{\boldsymbol{\beta}}\|_2 \leq \|\boldsymbol{\beta}^{(k-1)} - \widehat{\boldsymbol{\beta}}\|_2 \leq \ldots \leq \|\boldsymbol{\beta}^{(0)} - \widehat{\boldsymbol{\beta}}\|_2 \leq \|\boldsymbol{\beta}^{(0)} - \widehat{\boldsymbol{\beta}}\|_2.$$

Therefore, applying Lemma 5.1, we obtain

$$\|\boldsymbol{\beta}^{(k)} - \boldsymbol{\beta}^*\|_2 \leq \|\boldsymbol{\beta}^{(0)} - \boldsymbol{\beta}^*\|_2 + 2\|\widehat{\boldsymbol{\beta}} - \boldsymbol{\beta}^*\|_2 \leq \|\boldsymbol{\beta}^{(0)} - \boldsymbol{\beta}^*\|_2 + 18\rho_*\lambda\sqrt{s} \leq R/2.$$

This completes the induction step and thus finishes the proof. $\square$

We now give the proofs of Lemma E.4 and E.5.

PROOF OF LEMMA E.4. We omit the subscript in $\Psi_{\boldsymbol{\lambda}^{(\ell-1)},\phi_c}(\boldsymbol{\beta}, \boldsymbol{\beta}^{(\ell,k)})$, the super-script $\ell$ in $(\ell, k)$ and denote $\boldsymbol{\beta}^{(\ell,k)}$ as $\boldsymbol{\beta}^{(k)}$, where $\ell = 1$. Lemma E.3 makes us able to use the localized Lipschitz condition. First, we have

(E.10) $$F(\boldsymbol{\beta}, \boldsymbol{\lambda}^{(0)}) - F(\boldsymbol{\beta}^{(k)}, \boldsymbol{\lambda}^{(0)}) \geq F(\boldsymbol{\beta}, \boldsymbol{\lambda}^{(0)}) - \Psi(\boldsymbol{\beta}^{(k)}, \boldsymbol{\beta}^{(k-1)}).$$

The convexity of both $\mathcal{L}(\boldsymbol{\beta})$ and $\|\boldsymbol{\lambda} \odot \boldsymbol{\beta}\|_1$ implies

(E.11) $$L(\boldsymbol{\beta}) \geq \mathcal{L}(\boldsymbol{\beta}^{(k-1)}) + \langle \nabla \mathcal{L}(\boldsymbol{\beta}^{(k-1)}), \boldsymbol{\beta} - \boldsymbol{\beta}^{(k-1)} \rangle;$$

(E.12) $$\|\boldsymbol{\lambda}^{(0)} \odot \boldsymbol{\beta}\|_1 \geq \|\boldsymbol{\lambda}^{(0)} \odot \boldsymbol{\beta}^{(k)}\|_1 + \langle \boldsymbol{\lambda}^{(0)} \odot \boldsymbol{\xi}^{(k)}, \boldsymbol{\beta} - \boldsymbol{\beta}^{(k)} \rangle.$$

Adding (E.11) and (E.12) together, we obtain

(E.13) $$F(\boldsymbol{\beta}, \boldsymbol{\lambda}^{(0)}) \geq \mathcal{L}(\boldsymbol{\beta}^{(k-1)}) + \langle \nabla \mathcal{L}(\boldsymbol{\beta}^{(k-1)}), \boldsymbol{\beta} - \boldsymbol{\beta}^{(k-1)} \rangle \\ + \|\boldsymbol{\lambda}^{(0)} \odot \boldsymbol{\beta}^{(k)}\|_1 + \langle \boldsymbol{\lambda}^{(0)} \odot \boldsymbol{\xi}^{(k)}, \boldsymbol{\beta} - \boldsymbol{\beta}^{(k)} \rangle.$$

On the other side, $\Psi(\boldsymbol{\beta}^{(k)}, \boldsymbol{\beta}^{(k-1)})$ can be written as

$$
\begin{aligned}
\text{(E.14)} \quad & \mathcal{L}(\boldsymbol{\beta}^{(k-1)}) + \langle \nabla \mathcal{L}(\boldsymbol{\beta}^{(k-1)}), \boldsymbol{\beta}^{(k)} - \boldsymbol{\beta}^{(k-1)} \rangle \\
& + \frac{\phi_c}{2} \|\boldsymbol{\beta}^{(k)} - \boldsymbol{\beta}^{(k-1)}\|_2^2 + \|\boldsymbol{\lambda}^{(0)} \odot \boldsymbol{\beta}^{(k)}\|_1.
\end{aligned}
$$

Plugging (E.13) and (E.14) back into (E.10), we obtain

$$
\text{(E.15)} \quad F(\boldsymbol{\beta}, \boldsymbol{\lambda}^{(0)}) - F(\boldsymbol{\beta}^{(k)}, \boldsymbol{\lambda}^{(0)}) = -\frac{\phi_c}{2} \|\boldsymbol{\beta}^{(k)} - \boldsymbol{\beta}^{(k-1)}\|_2^2 \\
+ \langle \nabla \mathcal{L}(\boldsymbol{\beta}^{(k-1)}) + \boldsymbol{\lambda}^{(0)} \odot \boldsymbol{\xi}^{(k)}, \boldsymbol{\beta} - \boldsymbol{\beta}^{(k)} \rangle.
$$

By the first order optimality condition, there exists some $\boldsymbol{\xi}^{(k)}$ such that

$$
\nabla \mathcal{L}(\boldsymbol{\beta}^{(k-1)}) + \phi_c(\boldsymbol{\beta}^{(k)} - \boldsymbol{\beta}^{(k-1)}) + \boldsymbol{\lambda}^{(0)} \odot \boldsymbol{\xi}^{(k)} = 0.
$$

Plugging the equality above to (E.15), we complete the proof. $\square$

Lemma E.5 bounds the suboptimality measure $\omega_{\boldsymbol{\lambda}}(\boldsymbol{\beta}^{(1,k)})$ by $\|\boldsymbol{\beta}^{(1,k)} - \boldsymbol{\beta}^{(1,k-1)}\|_2$, which is critical to establish the iteration complexity of the contraction stage.

PROOF OF LEMMA E.5. We omit the super script 1 in $\boldsymbol{\beta}^{(1,k)}$ for simplicity. Since $\boldsymbol{\beta}^{(k)}$ is the exact solution to the $k$th iteration at $\ell = 1$, the first order optimality condition holds: there exists a $\boldsymbol{\xi}^{(k)} \in \partial \|\boldsymbol{\beta}^{(k)}\|_1$ such that

$$
\nabla \mathcal{L}(\boldsymbol{\beta}^{(k-1)}) + \phi_c(\boldsymbol{\beta}^{(k)} - \boldsymbol{\beta}^{(k-1)}) + \boldsymbol{\lambda} \odot \boldsymbol{\xi}^{(k)} = 0.
$$

Then for any $\mathbf{u}$ such that $\|\mathbf{u}\|_1 = 1$, we have

$$
\begin{aligned}
\langle \nabla \mathcal{L}(\boldsymbol{\beta}^{(k)}) + \boldsymbol{\lambda} \odot \boldsymbol{\xi}^{(k)}, \mathbf{u} \rangle &= \langle \nabla \mathcal{L}(\boldsymbol{\beta}^{(k)}), \mathbf{u} \rangle - \langle \nabla \mathcal{L}(\boldsymbol{\beta}^{(k-1)}) + \phi_c(\boldsymbol{\beta}^{(k)} - \boldsymbol{\beta}^{(k-1)}), \mathbf{u} \rangle \\
&= \langle \nabla \mathcal{L}(\boldsymbol{\beta}^{(k)}) - \nabla \mathcal{L}(\boldsymbol{\beta}^{(k-1)}), \mathbf{u} \rangle - \langle \phi_c(\boldsymbol{\beta}^{(k)} - \boldsymbol{\beta}^{(k-1)}), \mathbf{u} \rangle \\
&\leq \|\nabla \mathcal{L}(\boldsymbol{\beta}^{(k)}) - \nabla \mathcal{L}(\boldsymbol{\beta}^{(k-1)})\|_\infty + \phi_c \|\boldsymbol{\beta}^{(k)} - \boldsymbol{\beta}^{(k-1)}\|_\infty \\
&\leq (\phi_c + \rho_c) \|\boldsymbol{\beta}^{(k)} - \boldsymbol{\beta}^{(k-1)}\|_2,
\end{aligned}
$$

where the last inequality is due the the localized Lipchitz continuity, since $\|\boldsymbol{\beta}^{(k)} - \boldsymbol{\beta}^*\|_2 \leq R/2$, $\forall k \geq 1$ by Lemma E.3 in the supplement. The proof is completed by taking sup over $\|\mathbf{u}\|_1 \leq 1$ in the inequality above. $\square$

We then prove a technical lemma that is critical for the proof of Lemma 5.4.

**Lemma E.6** (Basic Inequality). Suppose the same conditions of Theorem 4.2 hold. Let $C = 225/(2\rho_*)$ and $\widetilde{\boldsymbol{\beta}}^{(1)}$ be the $\varepsilon$-optimal solution. Then we have the following basic inequality

$$\langle \nabla \mathcal{L}(\widetilde{\boldsymbol{\beta}}^{(1)}) - \nabla \mathcal{L}(\boldsymbol{\beta}^*), \widetilde{\boldsymbol{\beta}}^{(1)} - \boldsymbol{\beta}^* \rangle \leq C\lambda^2 s.$$

PROOF. We omit the superscript in $\widetilde{\boldsymbol{\beta}}^{(1)}$, and write $\widetilde{\boldsymbol{\beta}}^{(1)}$ as $\widetilde{\boldsymbol{\beta}}$ for simplicity. Proposition 4.1 implies that

$$(\text{E.16}) \quad \|\widetilde{\boldsymbol{\beta}} - \boldsymbol{\beta}^*\|_2^2 \leq \frac{12}{\rho_*}(\|\boldsymbol{\lambda}_S^{(0)}\|_2 + \|\nabla\mathcal{L}(\boldsymbol{\beta}^*)_S\|_2 + \varepsilon\sqrt{|S|}) \leq 15\lambda\sqrt{s}/\rho_*.$$

On the other side, applying Lemma F.1 yields that

$$\|\widetilde{\boldsymbol{\beta}} - \boldsymbol{\beta}^*\|_1 \leq \|(\widetilde{\boldsymbol{\beta}} - \boldsymbol{\beta}^*)_{\mathcal{E}_1^c}\|_1 + \|(\widetilde{\boldsymbol{\beta}} - \boldsymbol{\beta}^*)_{\mathcal{E}_1}\|_1 \leq 6\|(\widetilde{\boldsymbol{\beta}} - \boldsymbol{\beta}^*)_{\mathcal{E}_1}\|_1,$$

where $\mathcal{E}_1$ can be taken as $S$. This, combined with (E.16), results

$$(\text{E.17}) \quad \|(\widetilde{\boldsymbol{\beta}} - \boldsymbol{\beta})_S\|_1 \leq \sqrt{s}\|\widetilde{\boldsymbol{\beta}} - \boldsymbol{\beta}\|_2 \leq 15\lambda s/\rho_*.$$

Therefore, we obtain $\|\widetilde{\boldsymbol{\beta}} - \boldsymbol{\beta}^*\|_1 \leq 6\|(\widetilde{\boldsymbol{\beta}} - \boldsymbol{\beta}^*)_S\|_1 \leq 90\lambda s/\rho_*$. Because $\widetilde{\boldsymbol{\beta}}$ is a $\varepsilon$-optimal solution, we have

$$\langle \nabla\mathcal{L}(\widetilde{\boldsymbol{\beta}}) - \nabla\mathcal{L}(\boldsymbol{\beta}^*), \widetilde{\boldsymbol{\beta}} - \boldsymbol{\beta}^* \rangle \leq \|\nabla\mathcal{L}(\widetilde{\boldsymbol{\beta}}) + \boldsymbol{\lambda} \odot \widetilde{\boldsymbol{\xi}} - \boldsymbol{\lambda} \odot \widetilde{\boldsymbol{\xi}} - \nabla\mathcal{L}(\boldsymbol{\beta}^*)\|_\infty \|\widetilde{\boldsymbol{\beta}} - \boldsymbol{\beta}^*\|_1$$
$$\leq (1 + 1/4)\lambda\|\widetilde{\boldsymbol{\beta}} - \boldsymbol{\beta}\|_1 \leq 225\lambda s^2/(2\rho_*).$$

Therefore, the proof is completed. □

E.2.2. *Tightening Stage.* We collect technical lemmas that are needed to prove Lemma 5.5 and Proposition 4.6 in this section. We start by giving a lemma that ensures the sparsity along the approximate solution sequence $\{\boldsymbol{\beta}^{(\ell,k)}\}_{k=0}^\infty$ for the tightening stage ($\ell \geq 2$). We first need several technical lemmas. We remind the reader that the quadratic isotropic parameter in the tightening stage is denoted by $\phi_t$.

**Lemma E.7.** Suppose the same conditions in Theorem 4.2 hold. Assume $\boldsymbol{\beta}^{(k+1)}, \boldsymbol{\beta}^{(k)} \in B_2(r, \boldsymbol{\beta}^*)$ such that $\max\{\|\boldsymbol{\beta}_{S^c}^{(k+1)}\|_0, \|\boldsymbol{\beta}_{S^c}^{(k)}\|_0\} \leq \widetilde{s}$. For the LAMM algorithm, we have

$$\rho_-(2s + 2\widetilde{s}, r) \leq \phi_t \leq \gamma_u \rho_+(2s + 2\widetilde{s}, r).$$

PROOF. The proof follows a similar argument as that of Lemma E.3 and thus is omitted here for simplicity. □

The next two lemmas connects the suboptimality measure $\omega_{\boldsymbol{\lambda}}(\boldsymbol{\beta}^{(\ell,k)})$ to $\ell_2$ parameter bound and the objective functions. They are similar to the ones proved in the contraction stage but with different constants and we omit the proofs here.

**Lemma E.8.** If $\boldsymbol{\beta}^{(k-1)}, \boldsymbol{\beta}^{(k)} \in B_2(r/2, \boldsymbol{\beta}^*)$, $\|(\boldsymbol{\beta}^{(\ell,k)})_{S^c}\|_0 \leq \widetilde{s}$ and $\|(\boldsymbol{\beta}^{(\ell,k-1)})_{S^c}\|_0 \leq \widetilde{s}$, then for any $\ell \geq 2$ and $k \geq 1$, we have

$$\omega_{\boldsymbol{\lambda}}(\boldsymbol{\beta}^{(\ell,k)}) \leq (1+\gamma_u)\rho_+(2s+2\widetilde{s},r)\|\boldsymbol{\beta}^{(\ell,k)} - \boldsymbol{\beta}^{(\ell,k-1)}\|_2.$$

**Lemma E.9.** We have

$$\mathcal{F}(\boldsymbol{\beta}^{(\ell,k)}, \boldsymbol{\lambda}^{(\ell-1)}) - F(\boldsymbol{\beta}^{(\ell,k-1)}, \boldsymbol{\lambda}^{(\ell-1)}) \leq -\frac{\phi_t}{2}\|\boldsymbol{\beta}^{(\ell,k)} - \boldsymbol{\beta}^{(\ell,k-1)}\|_2.$$

Next we give a lemma that characterizes the parameter estimation and objective function bound for sparse approximate solutions.

**Lemma E.10.** Assume Assumption 4.1 holds. Let $\|\boldsymbol{\lambda}_{\mathcal{E}^c}\|_{\min} \geq \lambda/2$, $S \subset \mathcal{E}$ and $|\mathcal{E}| \leq 2s$. If $\|(\boldsymbol{\beta} - \boldsymbol{\beta}^*)_{S^c}\|_0 \leq \widetilde{s}$, $\omega_{\boldsymbol{\lambda}}(\boldsymbol{\beta}) \leq \varepsilon$ and $\boldsymbol{\beta} \in B_2(r, \boldsymbol{\beta}^*)$, then we must have

$$\|\boldsymbol{\beta} - \boldsymbol{\beta}^*\|_2 \leq 3\rho_*^{-1}\lambda\sqrt{s}/2,$$
$$F(\boldsymbol{\beta}, \boldsymbol{\lambda}) - F_{\boldsymbol{\lambda}}(\boldsymbol{\beta}^*, \boldsymbol{\lambda}) \leq 15\varepsilon\rho_*^{-1}\lambda s.$$

PROOF OF LEMMA E.10. For simplicity, we omit arguments $k, r$ in $\rho_-(k,r)$ and $\rho_+(k,r)$ when $k$ and $r$ are clear from the context. Since the sparse localized condition implies the localized sparse strong convexity, the following inequality follows from Proposition B.3:

$$\langle \nabla \mathcal{L}(\boldsymbol{\beta}) - \nabla \mathcal{L}(\boldsymbol{\beta}^*), \boldsymbol{\beta} - \boldsymbol{\beta}^* \rangle \geq \rho_- \|\boldsymbol{\beta} - \boldsymbol{\beta}^*\|_2^2.$$

Following a similar argument in the proof of Lemma B.7, we have

$$\text{(E.18)} \qquad \|\boldsymbol{\beta} - \boldsymbol{\beta}^*\|_2 \leq \frac{3\lambda\sqrt{s}}{2\rho_*}.$$

Next, we prove the desired bound for $F(\boldsymbol{\beta}, \boldsymbol{\lambda}) - F(\boldsymbol{\beta}^*, \boldsymbol{\lambda})$. Using the convexity of $F(\cdot, \boldsymbol{\lambda})$, we obtain

$$F(\boldsymbol{\beta}^*, \boldsymbol{\lambda}) \geq F(\boldsymbol{\beta}, \boldsymbol{\lambda}) + \langle \nabla \mathcal{L}(\boldsymbol{\beta}) + \boldsymbol{\lambda} \odot \boldsymbol{\xi}, \boldsymbol{\beta}^* - \boldsymbol{\beta}\rangle,$$

which yields that

$$\text{(E.19)} \quad F(\boldsymbol{\beta}, \boldsymbol{\lambda}) - F(\boldsymbol{\beta}^*, \boldsymbol{\lambda}) \leq -\langle \nabla \mathcal{L}(\boldsymbol{\beta}) + \boldsymbol{\lambda} \odot \boldsymbol{\xi}, \boldsymbol{\beta}^* - \boldsymbol{\beta}\rangle \leq \varepsilon\|\boldsymbol{\beta}^* - \boldsymbol{\beta}\|_1.$$

On the other hand, we know from Lemma F.1 that the approximate solution $\boldsymbol{\beta}$ falls in the $\ell_1$ cone:

$$\|(\boldsymbol{\beta}-\boldsymbol{\beta}^*)_{\mathcal{E}^c}\|_1 \leq 5\|(\boldsymbol{\beta}-\boldsymbol{\beta}^*)_{\mathcal{E}}\|_1,$$

which, together with (E.18), implies

$$(\text{E.20}) \quad \|\boldsymbol{\beta}-\boldsymbol{\beta}^*\|_1 \leq 6\|(\boldsymbol{\beta}-\boldsymbol{\beta}^*)_{\mathcal{E}}\|_1 \leq 6\sqrt{2s}\|(\boldsymbol{\beta}-\boldsymbol{\beta}^*)_{\mathcal{E}}\|_2 \leq 15\rho_*^{-1}\lambda s.$$

Plugging (E.20) into (E.19) completes the proof. □

**Lemma E.11** (Basic Inequality II). Assume Assumption 4.1 holds. Take $\mathcal{E}$ such that $S \subseteq \mathcal{E}$ and $|\mathcal{E}| \leq 2s$. Let $\lambda \geq 4(\|\nabla \mathcal{L}(\boldsymbol{\beta}^*)\|_\infty + \varepsilon)$ and $\|\boldsymbol{\lambda}_{\mathcal{E}^c}\|_{\min} \geq \lambda/2$. If $\|\boldsymbol{\beta}_{\mathcal{E}^c}\|_0 \leq \widetilde{s}$, $\boldsymbol{\beta} \in B_2(r, \boldsymbol{\beta}^*)$ and $F(\boldsymbol{\beta}, \boldsymbol{\lambda}) - F(\boldsymbol{\beta}^*, \boldsymbol{\lambda}) \leq C\lambda^2 s$, then

$$\frac{\rho_-(s+\widetilde{s}, r)}{2}\|\boldsymbol{\beta}-\boldsymbol{\beta}^*\|_2^2 + \frac{\lambda}{4}\|(\boldsymbol{\beta}-\boldsymbol{\beta}^*)_{\mathcal{E}^c}\|_1 \leq \frac{5\lambda}{4}\|(\boldsymbol{\beta}-\boldsymbol{\beta}^*)_{\mathcal{E}}\|_1 + C\lambda^2 s.$$

PROOF. Since $\|\boldsymbol{\beta}_{S^c}\|_0 \leq \widetilde{s}$ and $\|\boldsymbol{\beta}^*_{S^c}\|_0 = 0$, we have $\|(\boldsymbol{\beta}-\boldsymbol{\beta}^*)_{S^c}\|_0 \leq \widetilde{s}$. Proposition B.3 implies the localized sparse strong convexity:

$$(\text{E.21}) \quad L(\boldsymbol{\beta}^*) + \langle \nabla \mathcal{L}(\boldsymbol{\beta}^*), \boldsymbol{\beta}-\boldsymbol{\beta}^* \rangle + \frac{\rho_-}{2}\|\boldsymbol{\beta}-\boldsymbol{\beta}^*\|_2^2 \leq \mathcal{L}(\boldsymbol{\beta}).$$

Recall that $F(\boldsymbol{\beta}) = \mathcal{L}(\boldsymbol{\beta}) + \|\boldsymbol{\lambda} \odot \boldsymbol{\beta}\|_1$. We have $F(\boldsymbol{\beta}) - F(\boldsymbol{\beta}^*) \leq C\lambda^2 s$, or equivalently,

$$(\text{E.22}) \quad \mathcal{L}(\boldsymbol{\beta}) - \mathcal{L}(\boldsymbol{\beta}^*) + (\|\boldsymbol{\lambda} \odot \boldsymbol{\beta}\|_1 - \|\boldsymbol{\lambda} \odot \boldsymbol{\beta}^*\|_1) \leq C\lambda^2 s.$$

Plugging (E.21) into the left-hand side of (E.22), we immediately obtain

$$\frac{\rho_-}{2}\|\boldsymbol{\beta}-\boldsymbol{\beta}^*\|_2^2 \leq C\lambda^2 s \underbrace{-\langle \nabla \mathcal{L}(\boldsymbol{\beta}^*), \boldsymbol{\beta}-\boldsymbol{\beta}^* \rangle}_{\text{I}} + \underbrace{(\|\boldsymbol{\lambda} \odot \boldsymbol{\beta}^*\|_1 - \|\boldsymbol{\lambda} \odot \boldsymbol{\beta}\|_1)}_{\text{II}}.$$

Following a similar argument in the proof of Lemma B.7 in the appendix, we have

$$\text{I} \leq \|(\boldsymbol{\beta}-\boldsymbol{\beta}^*)_{\mathcal{E}^c}\|_1 \|\nabla \mathcal{L}(\boldsymbol{\beta}^*)\|_\infty + \|(\boldsymbol{\beta}-\boldsymbol{\beta}^*)_{\mathcal{E}}\|_1 \|\nabla \mathcal{L}(\boldsymbol{\beta}^*)\|_\infty,$$
$$\text{II} \leq \lambda\|(\boldsymbol{\beta}-\boldsymbol{\beta}^*)_{\mathcal{E}}\|_1 - \lambda/2\|(\boldsymbol{\beta}-\boldsymbol{\beta}^*)_{\mathcal{E}^c}\|_1.$$

Therefore, we have

$$\frac{\rho_-}{2}\|\boldsymbol{\beta}-\boldsymbol{\beta}^*\|_2^2 + (\lambda/2 - \|\nabla \mathcal{L}(\boldsymbol{\beta}^*)\|_\infty)\|(\boldsymbol{\beta}-\boldsymbol{\beta}^*)_{\mathcal{E}^c}\|_1$$
$$\leq (\lambda + \|\nabla \mathcal{L}(\boldsymbol{\beta}^*)\|_\infty)\|(\boldsymbol{\beta}-\boldsymbol{\beta}^*)_{\mathcal{E}}\|_1 + C\lambda^2 s.$$

The proof is finished by noticing that $\|\nabla \mathcal{L}(\boldsymbol{\beta}^*)\|_\infty \leq \lambda/4$. □

**Lemma E.12.** Assume Assumption 4.1 holds. Take $\mathcal{E}$ such that $S \subseteq \mathcal{E}$ and $|\mathcal{E}| \leq 2s$. Let $\|\nabla \mathcal{L}(\boldsymbol{\beta}^*)\|_\infty + \varepsilon \leq \lambda/4$ and $\|\boldsymbol{\lambda}_{\mathcal{E}^c}\|_{\min} \geq \lambda/2$. If $\boldsymbol{\beta} \in B_2(r, \boldsymbol{\beta}^*)$ satisfies $\|\boldsymbol{\beta}_{S^c}\|_0 \leq \widetilde{s}$ and $F(\boldsymbol{\beta}, \boldsymbol{\lambda}) - F(\boldsymbol{\beta}^*, \boldsymbol{\lambda}) \leq C\lambda^2 s$, then we must have

$$\|\boldsymbol{\beta} - \boldsymbol{\beta}^*\|_2 \leq C'\lambda\sqrt{s},$$
$$\langle \nabla \mathcal{L}(\boldsymbol{\beta}) - \nabla \mathcal{L}(\boldsymbol{\beta}^*), \boldsymbol{\beta} - \boldsymbol{\beta}^* \rangle \leq C'^2 \rho_+(2s + \widetilde{s}, r)\lambda^2 s,$$

where $C' = \max\{2\sqrt{C/\rho_-(2s+\widetilde{s}, r)}, 5\sqrt{2}/\rho_-(2s+\widetilde{s}, r)\}$.

PROOF. We omit the arguments in $\rho_-(k, r)$ and $\rho_+(k, r)$ when they are clear form the context. Directly applying Lemma E.11, it follows that

$$\frac{\rho_-}{2}\|\boldsymbol{\beta} - \boldsymbol{\beta}^*\|_2^2 \leq \frac{5\lambda}{4}\|(\boldsymbol{\beta} - \boldsymbol{\beta}^*)_\mathcal{E}\|_1 + C\lambda^2 s.$$

To further bound the right-hand side of the inequality above, we discuss two cases regarding the magnitude of $\|(\boldsymbol{\beta} - \boldsymbol{\beta}^*)_\mathcal{E}\|_1$ with respect to $\lambda s$:

- If $5\lambda\|(\boldsymbol{\beta} - \boldsymbol{\beta}^*)_\mathcal{E}\|_1/4 \leq C\lambda^2 s$, we have

(E.23) $\quad \dfrac{\rho_-}{2}\|\boldsymbol{\beta} - \boldsymbol{\beta}^*\|_2^2 \leq 2C\lambda^2 s,$ and thus $\|\boldsymbol{\beta} - \boldsymbol{\beta}^*\|_2 \leq 2\sqrt{\dfrac{C}{\rho_-}}\lambda\sqrt{s}.$

- If $5\lambda\|(\boldsymbol{\beta} - \boldsymbol{\beta}^*)_\mathcal{E}\|_1/4 > C\lambda^2 s$, we have

$$\frac{\rho_-}{2}\|\boldsymbol{\beta} - \boldsymbol{\beta}^*\|_2^2 \leq 5\lambda\|(\boldsymbol{\beta} - \boldsymbol{\beta}^*)_\mathcal{E}\|_1/2 \leq 5\lambda\sqrt{2s}\|\boldsymbol{\beta} - \boldsymbol{\beta}^*\|_2/2,$$

which further yields

(E.24) $$\|\boldsymbol{\beta} - \boldsymbol{\beta}^*\|_2 \leq \frac{5\sqrt{2}}{\rho_-}\lambda\sqrt{s}.$$

Combining (E.23) and (E.24), we obtain

$$\|\boldsymbol{\beta} - \boldsymbol{\beta}^*\|_2^2 \leq \max\left\{2\sqrt{\frac{C}{\rho_-}}, \frac{5\sqrt{2}}{\rho_-}\right\}\lambda\sqrt{s} = C'\lambda\sqrt{s},$$

where $C' = \max\{2\sqrt{C/\rho_-}, 5\sqrt{2}/\rho_-\}$. Using Proposition B.3, we obtain

$$D_\mathcal{L}^s(\boldsymbol{\beta}, \boldsymbol{\beta}^*) = \langle \mathcal{L}(\boldsymbol{\beta}) - \mathcal{L}(\boldsymbol{\beta}^*), \boldsymbol{\beta} - \boldsymbol{\beta}^* \rangle \leq \rho_+\|\boldsymbol{\beta} - \boldsymbol{\beta}^*\|_2^2 \leq C'^2\rho_+\lambda^2 s.$$

This completes the proof. $\square$

**Lemma E.13.** Assume Assumption 4.1 holds. Take $\mathcal{E}$ such that $S \subseteq \mathcal{E}$ with $|\mathcal{E}| \leq 2s$. Let Let $\|\nabla \mathcal{L}(\boldsymbol{\beta}^*)\|_\infty + \varepsilon \leq \lambda/4$ and $\|\boldsymbol{\lambda}_{\mathcal{E}^c}\|_{\min} \geq \lambda/2$. Let $\boldsymbol{\beta} \in B_2(r, \boldsymbol{\beta}^*)$ satisfy $\|\boldsymbol{\beta}_{\mathcal{E}^c}\|_0 \leq \widetilde{s}$ and $F(\boldsymbol{\beta}, \boldsymbol{\lambda}) - F(\boldsymbol{\beta}^*, \boldsymbol{\lambda}) \leq C\lambda^2 s$. Let $C_0 = 80\gamma_u \rho^*/\rho_* \max\{\sqrt{2C\rho_*}, 5\} + 64(\rho^*/\rho_*)^2 \max\{4C\rho_*, 50\}$. If $\widetilde{s} \geq C_0 s$, then the one-step LAMM algorithm produces a $(2s + \widetilde{s})$-sparse solution: $\|(T_{\boldsymbol{\lambda}, \phi_t}(\boldsymbol{\beta}))_{\mathcal{E}^c}\|_0 \leq \widetilde{s}$.

PROOF. For simplicity, we write $\bar{\boldsymbol{\beta}} = \boldsymbol{\beta} - \phi_t^{-1} \nabla \mathcal{L}(\boldsymbol{\beta})$. To show that $\|(S(\bar{\boldsymbol{\beta}}, \phi_t^{-1}\boldsymbol{\lambda}))_{\mathcal{E}^c}\|_0 \leq \widetilde{s}$, it suffices to prove that, for any $j \in \mathcal{E}^c$, the total number of $\beta_j$'s such that $\bar{\beta}_j > \lambda_j/\phi_t$ is no more than $\widetilde{s}$. We first write $\bar{\boldsymbol{\beta}}$ as

$$\bar{\boldsymbol{\beta}} = \boldsymbol{\beta} - \frac{1}{\phi_t} \nabla \mathcal{L}(\boldsymbol{\beta}) = \boldsymbol{\beta} - \frac{1}{\phi_t} \nabla \mathcal{L}(\boldsymbol{\beta}^*) + \frac{1}{\phi_t} \nabla \mathcal{L}(\boldsymbol{\beta}^*) - \frac{1}{\phi_t} \nabla \mathcal{L}(\boldsymbol{\beta}).$$

Define $S_n = \{j \in \mathcal{E}^c : (\boldsymbol{\beta} - \phi_t^{-1} \nabla \mathcal{L}(\boldsymbol{\beta}))_j = \lambda_j/\phi_t\}$, and notice that $\{j : (T_{\boldsymbol{\lambda}, \phi_t}(\boldsymbol{\beta}))_j \neq 0\} \subseteq S_n$, thus it suffices to show $|S_n| \leq \widetilde{s}$. We further define $S_n^1$, $S_n^2$ and $S_n^3$ as:

(E.25) $$S_n^1 \equiv \left\{j \in \mathcal{E}^c : |\boldsymbol{\beta}_j| \geq \frac{1}{4} \cdot \frac{\lambda_j}{\phi_t}\right\},$$

(E.26) $$S_n^2 \equiv \left\{j \in \mathcal{E}^c : |\nabla \mathcal{L}(\boldsymbol{\beta}^*)_j/\phi_t| > \frac{1}{2} \cdot \frac{\lambda_j}{\phi_t}\right\},$$

(E.27) $$S_n^3 \equiv \left\{j \in \mathcal{E}^c : \left|\left(\frac{\nabla \mathcal{L}(\boldsymbol{\beta}) - \nabla \mathcal{L}(\boldsymbol{\beta}^*)}{\phi_t}\right)_j\right| > \frac{1}{4} \cdot \frac{\lambda_j}{\phi_t}\right\}.$$

We immediately have $S_n \subseteq S_n^1 \cup S_n^2 \cup S_n^3$. It suffices to prove that $|S_n^1| + |S_n^2| + |S_n^3| \leq \widetilde{s}$. The assumption that $\|\nabla \mathcal{L}(\boldsymbol{\beta}^*)\|_\infty + \varepsilon \leq \lambda/4$ implies $S_n^2 = \emptyset$ and thus $|S_n^2| = 0$. In what follows, we bound $|S_n^1|$ and $|S_n^3|$, respectively.

**Bound for $|S_n^1|$:**
For $\forall j \in \mathcal{E}^c$, we have $\beta_j^\circ = 0$. Using Markov inequality, we obtain

$$|S_n^1| = \left|\left\{j \in \mathcal{E}^c : |\boldsymbol{\beta}_j| \geq \frac{1}{4} \cdot \frac{\lambda_j}{\phi_t}\right\}\right| \leq \sum_{j \in \mathcal{E}^c} \frac{4\phi_t}{\lambda_j} |\beta_j - \beta_j^*|.$$

Because $\|\boldsymbol{\lambda}_{\mathcal{E}^c}\|_{\min} \geq \lambda/2$, we have

$$|S_n^1| \leq \sum_{j \in \mathcal{E}^c} \frac{8\phi_t}{\lambda} |\beta_j - \beta_j^\circ| \leq \frac{8\phi_t}{\lambda} \|(\boldsymbol{\beta} - \boldsymbol{\beta}^*)_{\mathcal{E}^c}\|_1.$$

It remains to bound $\|(\boldsymbol{\beta} - \boldsymbol{\beta}^*)_{\mathcal{E}^c}\|_1$. A similar argument in Lemma E.11 implies

$$\frac{1}{4} \|(\boldsymbol{\beta} - \boldsymbol{\beta}^*)_{\mathcal{E}^c}\|_1 \leq \frac{5\lambda}{4} \|(\boldsymbol{\beta} - \boldsymbol{\beta}^*)_\mathcal{E}\|_1 + C\lambda^2 s.$$

Therefore, $\boldsymbol{\beta} - \boldsymbol{\beta}^*$ falls in the approximate $\ell_1$ cone:

$$\|(\boldsymbol{\beta} - \boldsymbol{\beta}^*)_{\mathcal{E}^c}\|_1 \leq 5\lambda\|(\boldsymbol{\beta} - \boldsymbol{\beta}^*)_{\mathcal{E}}\|_1 + 4C\lambda s \leq 5\sqrt{2}C'\lambda s + 4C\lambda s,$$

where $C' = \max\{2\sqrt{C/\rho_-}, 5\sqrt{2}/\rho_-\}$ and the last inequality is due to Lemma E.12. Let $C'' = \max\{10\sqrt{2C/\rho_-}, 50/\rho_-\} + 4C$, then we have

$$|S_n^1| \leq 8\phi_t C'' s \leq 8C'' \gamma_u \rho_+ s,$$

where we use the fact $\phi_t \leq \gamma_u \rho_+$ in the last inequality.

**Bound for $|S_n^3|$:**
Consider an arbitrary subset $S' \subseteq S_n^3$ with size $s' = |S'| \leq \widetilde{s}$. Let us further consider a $d$-dimensional sign vector $\mathbf{u}$ such that $\|\mathbf{u}\|_\infty = 1$ and $\|\mathbf{u}\|_0 = s'$. There exists some $\mathbf{u}$ such that

$$\frac{1}{8}\lambda s' \leq \sum_{j \in \mathcal{E}^c} \frac{1}{4}\lambda_j |u_j| \leq \mathbf{u}^T \{\nabla \mathcal{L}(\boldsymbol{\beta}) - \nabla \mathcal{L}(\boldsymbol{\beta}^*)\}.$$

By the Mean Value theorem, there exists some $\gamma \in [0, 1]$ such that $\nabla \mathcal{L}(\boldsymbol{\beta}) - \nabla \mathcal{L}(\boldsymbol{\beta}^*) = [\nabla^2 \mathcal{L}(\gamma \boldsymbol{\beta} + (1-\gamma)\boldsymbol{\beta}^*)](\boldsymbol{\beta} - \boldsymbol{\beta}^*)$. Let $\mathbf{H} \equiv [\nabla^2 \mathcal{L}(\gamma \boldsymbol{\beta} + (1-\gamma)\boldsymbol{\beta}^*)]$. Writing $\mathbf{u}^T(\nabla \mathcal{L}(\boldsymbol{\beta}) - \nabla \mathcal{L}(\boldsymbol{\beta}^*))$ as $\langle \mathbf{H}^{1/2}\mathbf{u}, \mathbf{H}^{1/2}(\boldsymbol{\beta} - \boldsymbol{\beta}^*)\rangle$ and applying the Hölder inequality, we obtain

$$(\text{E.28}) \quad \lambda s'/8 \leq \|\mathbf{H}^{1/2}\mathbf{u}\|_2 \|\mathbf{H}^{1/2}(\boldsymbol{\beta} - \boldsymbol{\beta}^*)\|_2 \leq \sqrt{\rho_+(s',r)s'} \underbrace{\|\mathbf{H}^{1/2}(\boldsymbol{\beta} - \boldsymbol{\beta}^*)\|_2}_{\text{I}}.$$

To bound term I, we apply Lemma E.12 and obtain that

$$\text{I} = \|\mathbf{H}^{1/2}(\boldsymbol{\beta} - \boldsymbol{\beta}^*)\|_2 \leq C'\sqrt{\rho_+(2s + \widetilde{s}, r)}\lambda\sqrt{s},$$

where $C' = \max\{2\sqrt{C/\rho_-}, 5\sqrt{2}/\rho_-\}$. Plugging the above inequality into (E.28), we obtain

$$\lambda s'/8 \leq \sqrt{\rho_+(s',r)}\sqrt{s'} \times C'\sqrt{\rho_+(2s + \widetilde{2}, r)}\lambda\sqrt{s}.$$

Taking squares of both sides yields

$$s' \leq 64\rho_+(s',r)C'^2 \rho_+(2s + \widetilde{s}, r)s \leq 64\rho_+(\widetilde{s}, r)C'^2 \rho_+(2s + \widetilde{s}, r)s < \widetilde{s}$$

where the last inequality is due to the assumption. Since $s' = |S'|$ achieves the maximum possible value such that $s' \leq \widetilde{s}$ for any subset $S'$ of $S_n^3$ and the above inequality shows that $s' < \widetilde{s}$, we must have

$$S' = \{j : |(\nabla \mathcal{L}(\boldsymbol{\beta}) - \nabla \mathcal{L}(\boldsymbol{\beta}^*))_j| \geq \lambda_j/4\}.$$

Finally, combining bounds for $|S_n^1|$, $|S_n^2|$ and $|S_n^3|$, we obtain

$$\|(T_{\boldsymbol{\lambda},\phi_t}(\boldsymbol{\beta}))_{E^c}\|_0 \leq 8C''\gamma_u\rho_+ s + 64\rho_+(\widetilde{s},r)C'^2\rho_+(2s+\widetilde{s},r)s \leq \widetilde{s}.$$

□

**Lemma E.14.** Assume the same conditions in Theorem 4.7 hold. The solution sequence $\{\boldsymbol{\beta}^{(\ell,k)}\}_{k=0}^{\infty}$ always satisfies that

$$F(\boldsymbol{\beta}^{(\ell,k)}, \boldsymbol{\lambda}^{(\ell-1)}) - F(\boldsymbol{\beta}^*, \boldsymbol{\lambda}^{(\ell-1)}) \leq C\lambda^2 s,$$
$$\|(\boldsymbol{\beta}^{(\ell,k)})_{\mathcal{E}_\ell^c}\|_0 \leq \widetilde{s}, \text{ and } \|\boldsymbol{\beta}^{(\ell,k)} - \boldsymbol{\beta}^*\|_2 \leq C'\lambda\sqrt{s}.$$

for $\ell \geq 2, k \geq 0$, where $C = 15/(4\rho_*)$ and $C' = 5\sqrt{2}/\rho_*$.

PROOF. We omit the argument $\boldsymbol{\lambda}$ in $F(\boldsymbol{\beta},\boldsymbol{\lambda})$, for notation simplicity. We prove the theorem by mathematical induction on $(\ell,k)$.

**Base case**: For the $\ell$th tightening step, the stopping criterion implies that $\omega_{\boldsymbol{\lambda}^{(\ell-2)}}(\boldsymbol{\beta}^{(\ell,0)}) \leq \varepsilon_t$. On the other hand, the suboptimality condition for the 1st iteration in the $\ell$th step can be written as

$$\omega_{\boldsymbol{\lambda}^{(\ell-1)}}(\boldsymbol{\beta}^{(\ell,0)}) = \min_{\boldsymbol{\xi} \in \partial \|\boldsymbol{\beta}^{(\ell,0)}\|_1} \left\{ \|\nabla \mathcal{L}(\boldsymbol{\beta}^{(\ell,0)}) + \boldsymbol{\lambda}^{(\ell-1)} \odot \boldsymbol{\xi}\|_\infty \right\}$$

which, together with the triangle inequality, yields

$$\omega_{\boldsymbol{\lambda}^{(\ell-1)}}(\boldsymbol{\beta}^{(\ell,0)}) \leq \min_{\boldsymbol{\xi}} \left\{ \|\nabla \mathcal{L}(\boldsymbol{\beta}^{(\ell,0)}) + \boldsymbol{\lambda}^{(\ell-2)} \odot \boldsymbol{\xi}\|_\infty + \underbrace{\|(\boldsymbol{\lambda}^{(\ell-1)} - \boldsymbol{\lambda}^{(\ell-2)}) \odot \boldsymbol{\xi}\|_\infty}_{\text{I}} \right\}.$$

For the second term I in the right hand side, we have

$$\|(\boldsymbol{\lambda}^{(\ell-1)} - \boldsymbol{\lambda}^{(\ell-2)}) \odot \boldsymbol{\xi}\|_\infty \leq \|\boldsymbol{\lambda}^{(\ell-1)} - \boldsymbol{\lambda}^{(\ell-2)}\|_\infty \leq \lambda/8.$$

Using the fact that $\varepsilon \leq \lambda/8$, we obtain

$$\omega_{\boldsymbol{\lambda}^{(\ell-1)}}(\boldsymbol{\beta}^{(\ell,0)}) \leq \min_{\boldsymbol{\xi}} \left\{ \|\nabla \mathcal{L}(\boldsymbol{\beta}^{(\ell,0)}) + \boldsymbol{\lambda}^{(\ell-2)} \odot \boldsymbol{\xi}\|_\infty + \lambda/8 \right\} \leq \lambda/4,$$

Thus the initialization satisfies that

$$\|(\boldsymbol{\beta}^{(\ell,0)})_{\mathcal{E}_\ell^c}\|_0 \leq \widetilde{s}, \ \omega_{\boldsymbol{\lambda}^{(\ell-1)}}(\boldsymbol{\beta}^{(\ell,0)}) \leq \lambda/4, \text{ and } \phi_t \leq \gamma_u\rho_+(2s+2\widetilde{s},r).$$

Therefore, using Lemma E.10, we obtain

$$F(\boldsymbol{\beta}^{(\ell,0)}) - F(\boldsymbol{\beta}^*) \leq C\lambda^2 s, \text{ where } C = 15/(4\rho_*).$$

Therefore, directly applying Lemma E.12 results

$$\|\boldsymbol{\beta}^{(\ell,0)} - \boldsymbol{\beta}^*\|_2 \leq C'\lambda\sqrt{s},$$

where $C' = 5\sqrt{2}\rho_*$.

**Induction step**: Suppose that, at the $(k-1)$-th iteration of the LAMM method in the $\ell$-th step, we have

$$\|(\boldsymbol{\beta}^{(\ell,k-1)})_{\mathcal{E}_\ell^c}\|_0 \leq \widetilde{s}, \ \phi \leq \gamma_u\rho_+, \ \text{and} \ F(\boldsymbol{\beta}^{(\ell,k-1)}) - F(\boldsymbol{\beta}^*) \leq C\lambda^2 s.$$

Then according to Lemma E.13, we have that the solution to the LAMM method at the $k$th iteration is $(2s+\widetilde{s})$-sparse: $\boldsymbol{\beta}^{(\ell,k)} = T_{\boldsymbol{\lambda}^{(\ell-1)},\phi_t}(\boldsymbol{\beta}^{(\ell,k-1)})$ satisfies $\|(\boldsymbol{\beta}^{(\ell,k)})_{\mathcal{E}_\ell^c}\|_0 \leq \widetilde{s}$. Thus Lemma E.9 implies

$$F(\boldsymbol{\beta}^{(\ell,k)}) \leq F(\boldsymbol{\beta}^{(\ell,k-1)}) - \frac{\phi_t}{2}\|\boldsymbol{\beta}^{(\ell,k)} - \boldsymbol{\beta}^{(\ell,k-1)}\|.$$

which implies that

$$F(\boldsymbol{\beta}^{(\ell,k)}) - F(\boldsymbol{\beta}^*) \leq F(\boldsymbol{\beta}^{(\ell,k-1)}) - F(\boldsymbol{\beta}^*) - \frac{\phi_t}{2}\|\boldsymbol{\beta}^{(\ell,k)} - \boldsymbol{\beta}^{(\ell,k-1)}\|_2^2 \leq C\lambda^2 s.$$

Therefore we have the induction holds at the $k$th iteration:

$$\|(\boldsymbol{\beta}^{(\ell,k)})_{\mathcal{E}_\ell^c}\|_0 \leq \widetilde{s}, \phi_t \leq \gamma_u\rho_+(2s+2\widetilde{s}), \ \text{and} \ F(\boldsymbol{\beta}^{(\ell,k)}) - F(\boldsymbol{\beta}^*) \leq C\lambda^2 s.$$

Using Lemma E.12, for $C'$ defined as before, we obtain

$$\|\boldsymbol{\beta}^{(\ell,k)} - \boldsymbol{\beta}^*\|_2 \leq C'\lambda\sqrt{s}.$$

We complete induction on $k$. For $\ell$, the proof is similar. $\square$

## APPENDIX F: PRELIMINARY LEMMAS

In this section, we collect several preliminary lemmas.

**Lemma F.1** ($\ell_1$ Cone Property For Approximate Solution). Let $\mathcal{E}$ such that $S \subseteq \mathcal{E}$. If $\|\nabla\mathcal{L}(\boldsymbol{\beta}^*)\|_\infty + \varepsilon \leq \|\boldsymbol{\lambda}_\mathcal{E}\|_{\min}$, we must have

$$\|(\widetilde{\boldsymbol{\beta}} - \boldsymbol{\beta}^*)_{\mathcal{E}^c}\|_1 \leq \frac{\|\boldsymbol{\lambda}\|_\infty + \|\nabla\mathcal{L}(\boldsymbol{\beta}^*)\|_\infty + \varepsilon}{\|\boldsymbol{\lambda}_\mathcal{E}\|_{\min} - (\|\nabla\mathcal{L}(\boldsymbol{\beta}^*)\|_\infty + \varepsilon)}\|(\widetilde{\boldsymbol{\beta}} - \boldsymbol{\beta}^*)_\mathcal{E}\|_1.$$

PROOF OF LEMMA F.1. For any $\boldsymbol{\xi} \in \partial \|\widetilde{\boldsymbol{\beta}}\|_1$, let $\mathbf{u} = \nabla \mathcal{L}(\widetilde{\boldsymbol{\beta}}) + \boldsymbol{\lambda} \odot \boldsymbol{\xi}$. By the Mean Value theory, there exists a $\gamma \in [0, 1]$, such that $\nabla \mathcal{L}(\widetilde{\boldsymbol{\beta}}) - \nabla \mathcal{L}(\boldsymbol{\beta}^*) = \left[\nabla^2 \mathcal{L}(\gamma \boldsymbol{\beta}^* + (1-\gamma)\widetilde{\boldsymbol{\beta}})\right](\widetilde{\boldsymbol{\beta}} - \boldsymbol{\beta}^*)$. Write $\mathbf{H} = \nabla^2 \mathcal{L}(\gamma \boldsymbol{\beta}^* + (1-\gamma)\widetilde{\boldsymbol{\beta}})$. Then we have

$$(\text{F.1}) \quad \langle \nabla \mathcal{L}(\widetilde{\boldsymbol{\beta}}) + \boldsymbol{\lambda} \odot \boldsymbol{\xi}, \widetilde{\boldsymbol{\beta}} - \boldsymbol{\beta}^* \rangle = \langle \nabla \mathcal{L}(\boldsymbol{\beta}^*) + H(\widetilde{\boldsymbol{\beta}} - \boldsymbol{\beta}^*), \widetilde{\boldsymbol{\beta}} - \boldsymbol{\beta}^* \rangle$$
$$\leq \|\mathbf{u}\|_\infty \|\widetilde{\boldsymbol{\beta}} - \boldsymbol{\beta}\|_1.$$

Using the fact $(\widetilde{\boldsymbol{\beta}} - \boldsymbol{\beta}^*)^T \mathbf{H} (\widetilde{\boldsymbol{\beta}} - \boldsymbol{\beta}^*) \geq 0$, we have

$$0 \leq \|\mathbf{u}\|_\infty \|\widetilde{\boldsymbol{\beta}} - \boldsymbol{\beta}^*\|_1 - \underbrace{\langle \nabla \mathcal{L}(\boldsymbol{\beta}^*), \widetilde{\boldsymbol{\beta}} - \boldsymbol{\beta}^* \rangle}_{\text{I}} - \underbrace{\langle \boldsymbol{\lambda} \odot \boldsymbol{\xi}, \widetilde{\boldsymbol{\beta}} - \boldsymbol{\beta}^* \rangle}_{\text{II}}.$$

Using a similar argument in the proof of Lemma B.7, we have $\text{I} \geq -\|\nabla \mathcal{L}(\boldsymbol{\beta}^*)\|_\infty \|\widetilde{\boldsymbol{\beta}} - \boldsymbol{\beta}\|_1$, and

$$\text{II} = \langle \boldsymbol{\lambda} \odot \boldsymbol{\xi}, \widetilde{\boldsymbol{\beta}} - \boldsymbol{\beta}^* \rangle = \langle (\boldsymbol{\lambda} \odot \boldsymbol{\xi})_{\mathcal{E}^c}, (\widetilde{\boldsymbol{\beta}} - \boldsymbol{\beta}^*)_{\mathcal{E}^c} \rangle + \langle (\boldsymbol{\lambda} \odot \boldsymbol{\xi})_{\mathcal{E}}, (\widetilde{\boldsymbol{\beta}} - \boldsymbol{\beta}^*)_{\mathcal{E}} \rangle$$
$$\geq \|\boldsymbol{\lambda}_{\mathcal{E}^c}\|_{\min} \|(\widetilde{\boldsymbol{\beta}} - \boldsymbol{\beta}^*)_{\mathcal{E}^c}\|_1 - \|\boldsymbol{\lambda}_{\mathcal{E}}\|_\infty \|(\widetilde{\boldsymbol{\beta}} - \boldsymbol{\beta}^*)_{\mathcal{E}}\|_1.$$

Plugging the above bounds into (F.1) and taking inf with respect to $\boldsymbol{\xi} \in \partial \|\widetilde{\boldsymbol{\beta}}\|_1$ yields

$$0 \leq -(\|\boldsymbol{\lambda}_{\mathcal{E}^c}\|_{\min} - (\|\nabla \mathcal{L}(\boldsymbol{\beta}^*)\|_\infty + \omega_{\boldsymbol{\lambda}}(\widetilde{\boldsymbol{\beta}})))\|(\widetilde{\boldsymbol{\beta}} - \boldsymbol{\beta}^*)_{\mathcal{E}^c}\|_1$$
$$+ (\|\boldsymbol{\lambda}_{\mathcal{E}^c}\|_{\min} + \|\nabla \mathcal{L}(\boldsymbol{\beta}^*)\|_\infty + \omega_{\boldsymbol{\lambda}}(\widetilde{\boldsymbol{\beta}}))\|(\widetilde{\boldsymbol{\beta}} - \boldsymbol{\beta}^*)_{\mathcal{E}}\|_1,$$

or equivalently

$$\|(\widetilde{\boldsymbol{\beta}} - \boldsymbol{\beta}^*)_{\mathcal{E}^c}\|_1 \leq \frac{\lambda + \|\nabla \mathcal{L}(\boldsymbol{\beta}^*)\|_\infty + \omega_{\boldsymbol{\lambda}}(\widetilde{\boldsymbol{\beta}})}{\|\boldsymbol{\lambda}_{\mathcal{E}^c}\|_{\min} - (\|\nabla \mathcal{L}(\boldsymbol{\beta}^*)\|_\infty + \omega_{\boldsymbol{\lambda}}(\widetilde{\boldsymbol{\beta}}))} \|(\widetilde{\boldsymbol{\beta}} - \boldsymbol{\beta}^*)_{\mathcal{E}}\|_1$$

Using the stopping criterion, i.e. $\omega_{\boldsymbol{\lambda}}(\widetilde{\boldsymbol{\beta}}) \leq \varepsilon$, we have that

$$\|(\widetilde{\boldsymbol{\beta}} - \boldsymbol{\beta}^*)_{\mathcal{E}^c}\|_1 \leq \frac{\lambda + \|\nabla \mathcal{L}(\boldsymbol{\beta}^*)\|_\infty + \varepsilon}{\|\boldsymbol{\lambda}_{\mathcal{E}^c}\|_{\min} - (\|\nabla \mathcal{L}(\boldsymbol{\beta}^*)\|_\infty + \varepsilon)} \|(\widetilde{\boldsymbol{\beta}} - \boldsymbol{\beta}^*)_{\mathcal{E}}\|_1$$

Therefore we proved the desired result. □

**Lemma F.2.** Let $D_{\mathcal{L}}(\boldsymbol{\beta}_1, \boldsymbol{\beta}_2) = \mathcal{L}(\boldsymbol{\beta}_1) - \mathcal{L}(\boldsymbol{\beta}_2) - \langle \mathcal{L}(\boldsymbol{\beta}_2), \boldsymbol{\beta}_1 - \boldsymbol{\beta}_2 \rangle$ and $D^s_{\mathcal{L}}(\boldsymbol{\beta}_1, \boldsymbol{\beta}_2) = D_{\mathcal{L}}(\boldsymbol{\beta}_1, \boldsymbol{\beta}_2) + D_{\mathcal{L}}(\boldsymbol{\beta}_2, \boldsymbol{\beta}_1)$. For $\boldsymbol{\beta}(t) = \boldsymbol{\beta}^* + t(\boldsymbol{\beta} - \boldsymbol{\beta}^*)$ with $t \in (0, 1]$, we have that

$$D^s_{\mathcal{L}}(\boldsymbol{\beta}(t), \boldsymbol{\beta}^*) \leq t D^s_{\mathcal{L}}(\boldsymbol{\beta}, \boldsymbol{\beta}^*).$$

PROOF OF LEMMA F.2. Let $Q(t) = D_{\mathcal{L}}(\boldsymbol{\beta}(t), \boldsymbol{\beta}^*) = \mathcal{L}(\boldsymbol{\beta}(t)) - \mathcal{L}(\boldsymbol{\beta}^*) - \langle \nabla \mathcal{L}(\boldsymbol{\beta}^*), \boldsymbol{\beta}(t) - \boldsymbol{\beta}^* \rangle$. Since the derivative of $\mathcal{L}(\boldsymbol{\beta}(t))$ with respect to $t$ is $\langle \nabla \mathcal{L}(\boldsymbol{\beta}(t)), \boldsymbol{\beta} - \boldsymbol{\beta}^* \rangle$, it follows that

$$Q'(t) = \langle \nabla \mathcal{L}(\boldsymbol{\beta}(t)) - \nabla \mathcal{L}(\boldsymbol{\beta}^*), \boldsymbol{\beta} - \boldsymbol{\beta}^* \rangle.$$

Therefore, the symmetric Bregman divergence $D^s_{\mathcal{L}}(\boldsymbol{\beta}(t) - \boldsymbol{\beta}^*)$ can be written as

$$D^s_{\mathcal{L}}(\widetilde{\boldsymbol{\beta}}(t) - \boldsymbol{\beta}^*) = \langle \nabla \mathcal{L}(\widetilde{\boldsymbol{\beta}}(t)) - \nabla \mathcal{L}(\boldsymbol{\beta}^*), t(\boldsymbol{\beta} - \boldsymbol{\beta}^*) \rangle = tQ'(t) \text{ for } 0 < t \leq 1.$$

Plugging $t = 1$ in the equation above, we have $Q'(1) = D^s_{\mathcal{L}}(\boldsymbol{\beta}, \boldsymbol{\beta}^*)$ as a special case. If we assume that $Q(t)$ is convex, then $Q'(t)$ is non-decreasing and thus

$$D^s_{\mathcal{L}}(\boldsymbol{\beta}(t), \boldsymbol{\beta}^*) = tQ'(t) \leq tQ'(1) = tD^s_{\mathcal{L}}(\boldsymbol{\beta}, \boldsymbol{\beta}^*).$$

It remains to show the convexity of $Q(t)$ with respect to $t$; or equivalently, the convexity of $\mathcal{L}(\boldsymbol{\beta}(t))$ and $\langle \nabla \mathcal{L}(\boldsymbol{\beta}^*), \boldsymbol{\beta}^* - \boldsymbol{\beta}(t) \rangle$, respectively. First, we have the fact that $\boldsymbol{\beta}(t)$ is linear in $t$, that is, $\boldsymbol{\beta}(\alpha_1 t_1 + \alpha_2 t_2) = \alpha_1 \boldsymbol{\beta}(t_1) + \alpha_2 \boldsymbol{\beta}(t_2)$, for $t_1, t_2 \in [0, 1]$ and $\alpha_1, \alpha_2 \geq 0$ such that $\alpha_1 + \alpha_2 = 1$. Then the convexity of $\mathcal{L}(\boldsymbol{\beta}(t))$ follows from this linearity property of $\boldsymbol{\beta}(t)$ and the convexity of the Huber loss. For the second term, the convexity directly follows from the bi-linearity of the inner product function. This finishes the proof. □

The following lemma is taken from [2] and describes a general concentration for quadratic forms in sub-Gaussian random variables.

**Lemma F.3.** (Hanson-Wright Inequality, [2]).
Let $\mathbf{v} = (v_1, \ldots, v_d) \in \mathbb{R}^d$ be a random vector with independent components $v_i$ such that $v_i \sim \text{sub-Gaussian}(0, \sigma^2)$. Let $\mathbf{A}$ be an $n \times n$ matrix. Then, for every $t \geq 0$,

$$\mathbb{P}\Big(|\mathbf{v}^T \mathbf{A} \mathbf{v} - \mathbb{E} \mathbf{v}^T \mathbf{A} \mathbf{v}| > t\Big) \leq 2 \exp\bigg(- C_h \min\Big\{\frac{t^2}{\sigma^4 \|\mathbf{A}\|_F^2}, \frac{t}{\sigma^2 \|\mathbf{A}\|_2}\Big\}\bigg),$$

where $C_h$ is a universal constant, not depending on $\mathbf{A}, \mathbf{v}$ and $n$.

## REFERENCES

[1] AGARWAL, A., NEGAHBAN, S. and WAINWRIGHT, M. J. (2012). Fast global convergence rates of gradient methods for high-dimensional statistical recovery. *The Annals of Statistics* **40** 2452–2482.

[2] RUDELSON, M. and VERSHYNIN, R. (2013). Hanson-wright inequality and sub-gaussian concentration. *arXiv preprint arXiv:1306.2872* .

DEPARTMENT OF OPERATIONS RESEARCH
AND FINANCIAL ENGINEERING
PRINCETON UNIVERSITY
PRINCETON, NJ 08544;
SCHOOL OF DATA SCIENCE
FUDAN UNIVERSITY
SHANGHAI, CHINA
E-MAIL: jqfan@princeton.edu

DEPARTMENT OF OPERATIONS RESEARCH
AND FINANCIAL ENGINEERING
PRINCETON UNIVERSITY
PRINCETON, NJ 08544
E-MAIL: hanliu@princeton.edu
qsun.ustc@gmail.com

TENCENT AI LAB
SHENNAN AVE, NANSHAN DISTRICT
SHEN ZHEN, GUANGDONG, CHINA;
SCHOOL OF DATA SCIENCE
FUDAN UNIVERSITY
SHANGHAI, CHINA
E-MAIL: tongzhang@tongzhang-ml.org